\documentstyle[amssymb,12pt]{article}

\input amssym.def
\input amssym

\newtheorem{thm}{Theorem}[section]
\newtheorem{lem}[thm]{Lemma}
\newtheorem{prop}[thm]{Proposition}
\newtheorem{cor}[thm]{Corollary}
\newtheorem{fact}[thm]{Fact}

\newtheorem{defn}[thm]{Definition}

\newtheorem{nrmk}[thm]{Remark}

\newcommand{\qed}{\hfill $\Box$ \vspace{.5cm}}
\newcommand{\pf}{{\bf Proof. }}

\title {Covers of groups definable in o-minimal structures}

\author {M\'{a}rio J. Edmundo \thanks{Supported by the EPSRC grant GR/M66332.}
\\The Mathematical Institute
\\24-29 St Giles
\\OX1 3LB  Oxford, U.K
\\edmundo@maths.ox.ac.uk
\\
\\}
\date{March 25, 2001}

\newcommand{\into}{\longrightarrow}

\renewcommand{\tilde}{\widetilde}
\renewcommand{\bar}{\overline}

\newcommand{\NN}{\mathbb{N}}
\newcommand{\ZZ}{\mathbb{Z}}
\newcommand{\QQ}{\mathbb{Q}}
\newcommand{\RR}{\mathbb{R}}

\newcommand{\M}{\mbox{$\cal M$}}
\newcommand{\N}{\mbox{$\cal N$}}

\begin{document}

\maketitle
\begin{abstract}
We develop in this paper the theory of covers for Hausdorff properly 
$\bigvee $-definable 
manifolds with definable choice in an o-minimal structure $\N$. 
In particular, we show that given an $\N$-definably connected
$\N$-definable group $G$ we have 
$1\rightarrow \pi _1(G)\rightarrow \tilde{G}\stackrel{p}\rightarrow 
G\rightarrow 1$ in the category of strictly properly $\bigvee
$-definable groups with strictly properly $\bigvee $-definable
homomorphisms, where $\pi _1(G)$ is the o-minimal fundamental group of
$G$.
\end{abstract}

\newpage

\begin{section}{Introduction}\label{section introduction}

Throughout this paper, $\N$ will be an o-minimal structure and definable
will mean $\N$-definable (possibly with parameters). We will assume
the readers familiarity with the basic theory of o-minimal structures
(see for example \cite{vdd}).

In section \ref{section manifolds}, we will introduce
several categories - properly $\bigvee $-definable manifolds with 
(properly) $\bigvee $-definable maps and properly $\bigvee $-definable 
manifolds with strictly (properly) $\bigvee $-definable maps - 
generalising the category of definable manifolds with definable maps. 
Properly $\bigvee $-definable subsets of properly $\bigvee $-definable 
manifolds, which will play an important role, are
introduced and several notions such as properly $\bigvee $-definable
cell decomposition, properly $\bigvee $-definable completeness,
properly $\bigvee $-definable connectedness, dimensions and generic 
points are studied. These notions, generalise corresponding notions
for definable subsets of definable manifolds and are shown to be
invariant under the model theoretic operations of taking elementary 
extensions, elementary substructures, reducts and expansions.  

In section \ref{section fundamental group}, we introduce the o-minimal 
fundamental group functor for the category of properly 
$\bigvee $-definable manifolds with definable choice together with 
continuous strictly properly $\bigvee $-definable maps. The
construction of this functor is based on $k$-cells ($k=0,1,2$) with an 
orientation, and we prove that this functor satisfies all the relevant 
properties that one should expect for the fundamental group functor. 
Namely, the o-minimal Tietze and Seifert-van Kampen theorems are
proved for  properly $\bigvee $-definably complete, locally finite
properly $\bigvee $-definable manifolds with
definable choice, generalising results from \cite{bo} on o-minimal 
fundamental groups of definable sets in o-minimal expansions of real 
closed fields. Moreover, the o-minimal fundamental group of a properly
$\bigvee $-definably complete, properly $\bigvee $-definable manifold with
definable choice is proved to be 
invariant under the model theoretic operations of 
taking elementary extensions, elementary substructures, reducts and expansions.
In section \ref{section covering spaces}, all the
theory of strictly properly $\bigvee $-definable covering spaces is 
developed in the category of properly $\bigvee $-definable manifolds 
with definable choice together with continuous strictly properly 
$\bigvee $-definable maps. 

Finally, in section \ref{section  strictly properly bigvee definable 
groups} we apply our previous theory to $\bigvee $-definable groups -
which we prefer to call in this paper, strictly 
properly $\bigvee $-definable groups mainly for two reasons: first 
because we do not assume that $\N$ is $\aleph _1$-saturated (as in 
\cite{pst2}, where part of the theory of $\bigvee $-definable groups
is developed) and secondly because, as we show throughout the paper, the 
distinction between properly $\bigvee $-definable subgroups and 
$\bigvee $-definable subgroups is a very important one, with a better 
theory in the first case. We start section 
\ref{section  strictly properly bigvee definable groups}  by 
developing the basic theory of 
strictly properly $\bigvee $-definable groups: existence of a
properly $\bigvee $-definable manifold structure making the group 
operations and strictly properly $\bigvee $-definable homomorphisms 
continuous (this is already implicit in \cite{pst2}); DCC for strictly 
properly $\bigvee $-definable subgroups; existence of infinite
strictly properly $\bigvee $-definable abelian subgroups of infinite 
strictly properly $\bigvee $-definable groups; existence of strictly 
properly $\bigvee $-definable quotient of a strictly properly 
$\bigvee $-definable group by a strictly properly $\bigvee $-definable 
normal subgroup and existence of a corresponding strictly properly 
$\bigvee $-definable section; centerless strictly properly 
$\bigvee $-definable groups with no strictly properly $\bigvee
$-definable normal subgroups of positive dimension are shown to be the 
$\bigvee $-definable open and closed subgroups of definably semisimple 
definable groups generated by open definable subsets and,
the classification of 
solvable (and therefore by the above, of all) strictly properly 
$\bigvee $-definable groups is reduced to the classification of
properly $\bigvee $-definably complete, strictly properly 
$\bigvee $-definable solvable groups; existence of definable choice is 
proved for strictly properly $\bigvee $-definable groups
and the theory of strictly properly $\bigvee $-definable
coverings of strictly properly $\bigvee $-definable groups is presented.

There are two problems that we do not handle in this paper. The first 
is the classification of properly $\bigvee $-definably 
complete, properly $\bigvee $-definably connected, strictly properly 
$\bigvee $-definable solvable groups. We consider this problem in more 
detail in \cite{e3}. The second problem is the computation of the 
o-minimal fundamental group of a definable group (or even more
generally of a strictly properly $\bigvee $-definable group). We show
here that such groups are abelian and finitely generated, and in
\cite{e2} where we compute the o-minimal fundamental groups of groups 
definable in o-minimal expansions of real closed fields, we reduce as
well this problem to the problem of computing the o-minimal
fundamental groups of definably compact, definable abelian
groups. We show here that the o-minimal fundamental group $\pi _1(G)$ of a 
definably compact, definable abelian group $G$ is a torsion-free finitely 
generated abelian group. We conjecture that the rank of $\pi _1(G)$
the equals the dimension of $G$. 
This problem however, can only be solved using a general cohomology
theory for general o-minimal structures which we hope to develop in \cite{e4}.

\end{section}

\begin{section}{Properly $\bigvee $-definable manifolds}
\label{section manifolds}

\begin{subsection}{Properly $\bigvee $-definable manifolds}
\label{subsection manifolds}

\begin{defn}\label{defn nspaces}
{\em  
A {\it properly $\bigvee $-definable manifold  (over $A$) of dimension $m$}
where $A\subseteq N$ is such that $|A|<\aleph _0$, is
a triple $\mathbf{X}$ $:=(X,(X_i,\phi _i)_{i\in I})$ where $\{X_i:i\in
I\}$ is cover of the set $X$ with $|I|<\aleph _1 $ such that for each
$i\in I$, if $I_i:=\{j\in I:X_i\cap X_j\neq \emptyset \}$ then
we have injective maps $\phi _i:X_i\into N^m$ such that
$\phi _i(X_i)$ is an open definably connected definable set
(over $A$), for all $j\in I_i$, $\phi _i(X_i\cap X_j)$ is an
definable (over $A$) open subset of $\phi_i(X_i)$ and the map
$\phi _{ij}:\phi _i(X_i\cap X_j)\into \phi _j(X_i\cap X_j)$ given by
$\phi _{ij}:=\phi _j\circ \phi _i^{-1}$ is a definable homeomorphism (over
$A$). 

A properly $\bigvee $-definable manifold  (over $A$) of dimension $m$ will be
called a 
{\it locally finite properly $\bigvee $-definable manifold  (over $A$) 
of dimension $m$} if $|I_i|<\aleph _0$; a properly $\bigvee $-definable
manifold  (over $A$) of dimension $m$ will be
called a {\it definable manifold  (over $A$) of dimension $m$} if
$|I|<\aleph _0$. 
}  
\end{defn}

For the rest of the section, let $\mathbf{X}$$=(X,(X_i,\phi _i)_{i\in
I})$ and $\mathbf{Y}$ $=(Y,(Y_j,\psi _j)_{j\in J})$ be properly 
$\bigvee $-definable manifolds of dimension $m$ and $n$ respectively and
defined over $A_X$ and $A_Y$ respectively. Note that,
$\mathbf{X\times Y}$ $=(X\times Y,(X_i\times Y_j, 
(\phi _i,\psi _j))_{i\in I,j\in J})$ is then a properly $\bigvee $-definable
manifold of dimension $m\times n$ and defined over $A_X\cup A_Y$.

\begin{defn}\label{defn subsets}
{\em
Given $Z\subseteq X$ let $I^Z:=\{i\in I:Z\cap X_i\neq \emptyset
\}$ (necessarily $|I^Z|<\aleph _1$) and  for each $i\in I^Z$, let 
$Z_i:=Z\cap X_i$. Let $A_X\subseteq B\subseteq N$ be such that
$|B|<\aleph _0$. We say that $Z$ is a 
{\it properly $\bigvee $-definable subset of $X$  (over $B$)} if 
for each $i\in I^Z$, the set $\phi _i(Z_i)$ is a definable subset
of $\phi _i(X_i)$ (over $B$); we say that $Z$ is a {\it
definable subset of $X$ (over $B$)} if $Z$ is a  properly 
$\bigvee $-definable subset of $X$ (over $B$) and there is a finite
subset $I_0$ of $I$ such that $Z\subseteq \cup \{X_i:i\in I_0\}$. 
Finally, we
say that $Z$ is a {\it $\bigvee $-definable subset of $X$  (over $B$)}
if $Z=\cup \{Z^s:s\in S\}$ where $|S|<\aleph _1$ and for each $s\in
S$, $Z^s$ is  is a definable subset of $X$ (over $B$).
}
\end{defn}

Of course, a definable subset of $X$ is a properly 
$\bigvee $-definable subset of $X$, a properly $\bigvee $-definable subset
of $X$ is a $\bigvee $-definable subset of $X$ and $X$ is
always a properly $\bigvee $-definable subset of $X$. By a definable (resp., 
properly $\bigvee $-definable and $\bigvee $-definable) subset $U$ of
a $\bigvee $-definable subset $Z$ of $X$ we mean a definable (resp., 
properly $\bigvee $-definable and $\bigvee $-definable) subset $U$ of
$X$ which is a subset of $Z$. Definable subsets of $X$ have a very
well behaved theory, induced from the theory of definable sets in
$N^m$. The collection of properly $\bigvee
$-definable subsets of $X$ which is closed under finite unions, finite 
intersections and under taking complements, but not under the
projection maps (in fact, general $\bigvee $-definable subsets
typically occur in this way -see the examples below), will also have
a well behaved theory as we shall show throughout the paper. However,
general $\bigvee $-definable subsets do not have an interesting theory 
when $\N$ is not $\aleph _1$-saturated, consider for
example in $\N$$=(\RR$$,<)$ the $\bigvee $-definable subset $\QQ$ of
$\RR$. 
Fortunately, general $\bigvee $-definable subsets will not play 
an essential role in the construction of o-minimal fundamental groups
and covers, so we will work most of the time without the 
$\aleph _1$-saturation condition.

\begin{defn}
{\em
Let $Z$ be a properly $\bigvee $-definable subset of $X$  (over
$B$). A map $f:Z\subseteq X\into Y$ is a {\it properly $\bigvee
$-definable map  (over $B$)} if
its graph $\Gamma (f)$ is a properly $\bigvee $-definable subset of
$X\times Y$ (over $B$). Similarly, we say that a properly
$\bigvee $-definable map $f:Z\into Y$  is a {\it definable map  (over
$B$)} if its graph $\Gamma (f)$ is a definable subset of $X\times Y$
(over $B$). And finally if $Z$ is a $\bigvee
$-definable subset of $X$ (over $B$) and $f:Z\into Y$ is a map, 
we say that $f$ is a {\it $\bigvee $-definable map  (over
$B$)} if $Z=\cup \{Z^s:s\in S\}$ with $|S|<\aleph _1$, 
where for each $s\in S$, $Z^s$ is a definable subset of $X$ (over $B$) and 
$f_{|Z^s}$ is a definable map (over $B$).
}
\end{defn}

Note that: $(1)$  $f:Z\subseteq X\into Y$ is a properly $\bigvee
$-definable map  (over $B$) iff for all $i\in I^Z$, and for all 
$j\in J^{f(Z_i)}$, the map $\psi _j\circ f\circ \phi _{i|}^{-1}: 
\phi _i(Z_i)\into \psi _j(Y_j)$ is a definable map (over $B$); $(2)$
if $f:Z\into Y$ is a properly $\bigvee $-definable map  (over $B$) then
the image $f(U)$ of a definable subset $U$ of $Z$ (over $B$) is not 
necessarily a definable subset of $Y$ but a it is a properly 
$\bigvee $-definable subset of $Y$ (over $B$) and therefore, the 
image $f(U)$ of a (properly) $\bigvee $-definable subset $U$ of $Z$ 
(over $B$) is in general a $\bigvee $-definable subset of $Y$ (over
$B$) (e.g., take $\N$$=(\QQ$$,<)$, $I=\NN$ and $X=\cup \{X_i: i\in
I\}$  with $X_i=\{i\}$, $J=\{1\}$ and $Y=Y_1=\QQ$, and
$f(n)=\frac{1}{n}$); $(3)$ if $f:Z\into Y$ is a properly $\bigvee $-
definable map  (over $B$) then the inverse image $f^{-1}(W)$ of a
definable subset $W$ of $f(Z)$ (over $B$) is not necessarily a
definable subset of $Z$ but a it is a properly $\bigvee $-definable
subset of $Z$ (over $B$) and therefore, the inverse image $f^{-1}(W)$
of a (properly) $\bigvee $-definable subset $W$ of $f(Z)$ (over $B$)
is in general a $\bigvee $-definable subset of $Z$ (over $B$) (e.g.,
take $\N$$=(\QQ$$,<)$, $I=\{1\}$, $X=X_1=\QQ$ $^{\geq 0}$, $J=\NN$ and 
$Y=\cup \{Y_j:j\in J\}$ with $Y_j=\{j\}$ and $f(x)=2k$ if $x\in (\frac{k}{2},
\frac{k+1}{2})$ and $f(\frac{k}{2})=2k+1$, where $k\in \NN$ and take 
$W=2\NN$$+1$) and $(4)$ if $f:Z\subseteq X\into Y$ is a properly $\bigvee
$-definable map  (over $B$) then, we have 
$Z=\cup \{Z^s:s\in S\}$ with $|S|<\aleph _1$, 
where for each $s\in S$, $Z^s$ is a definable subset of $X$ (over $B$) and 
$f_{|Z^s}$ is a definable map (over $B$).

\begin{defn}
{\em
Let $Z$ be a properly $\bigvee $-definable subset of $X$  (over
$B$). A properly $\bigvee $-definable map $f:Z\subseteq X\into Y$
(over $B$) is a {\it strictly properly $\bigvee $-definable map  (over
$B$)} if for all $i\in I^Z$, $f(Z_i)$ is a definable
subset of $Y$. Finally, any $\bigvee $-definable map $f:Z\subseteq X\into Y$
(over $B$) is a {\it strictly $\bigvee $-definable map (over $B$)}
since we can write $Z=\cup \{Z^s:s\in S\}$ with $|S|<\aleph _1$, 
where for each $s\in S$, $Z^s$ is a definable subset of $X$ (over $B$) and 
$f_{|Z^s}$ is a definable map (over $B$) such that $f(Z_s)$ is a
definable subset of $Y$.
}
\end{defn}

$\mathbf{X}$ can be made into a topological space: the basis for the
topology is the collection of {\it open definable subsets of $X$}
i.e., definable subsets $U$ of $X$  such that for all $i\in I^U$, 
$\phi _i(U_i)$ is an open definable subset of $\phi _i(X_i)$. 
We will often identify two properly $\bigvee $-definable manifolds 
${\mathbf X}$ and ${\mathbf Y}$ if $X=Y$ and the identity map $1_X:X\into Y$ is
a strictly (properly) $\bigvee $-definable homeomorphism.
$\mathbf{Y}$ is a {\it strictly (properly) $\bigvee $-definable
submanifold} of $\mathbf{X}$ if $Y$ is a (properly) $\bigvee
$-definable subset of $X$, and the inclusion map $Y\into Y\subseteq X$
is a strictly (properly) $\bigvee $-definable homeomorphism onto its image.
The strictly (properly) $\bigvee$-definable submanifolds
$\mathbf{Y}$  of $N^n$ are called {\it strictly (properly) 
$\bigvee $-definable affine manifolds}.

\end{subsection}

\begin{subsection}{Properly $\bigvee $-definable cell decomposition}
\label{subsection cell decomposition}

In this subsection ${\mathbf X}$ will be a locally finite
properly $\bigvee $-definable 
manifold. There are many geometric properties of definable sets and
definable maps in the o-minimal structure $\N$. However, 
two of the most powerful results that we will be using throughout this
paper are the monotonicity theorem for 
definable one variable functions and the $C^p$-cell decomposition 
theorem for definable sets and definable maps. We will now explain the
$C^p$-cell decomposition theorem (here $p=0$ if $\N$ is not an
expansion of a (real closed) field) in order to introduce as well 
the notions of properly $\bigvee $-definable $C^p$-cell decomposition and
o-minimal dimension of properly $\bigvee $-definable subsets of $X$.

\begin{defn}\label{defn cells and dimension}
{\em
$C^p$-cells and o-minimal dimension are defined inductively as
follows: $(i)$ the unique non empty definable subset of $N^0$ is a 
$C^p$-cell of dimension zero, a point in $N^1$ is a $C^p$-cell of
dimension zero and an open interval in $N^1$ is a $C^p$-cell of
dimension one; $(ii)$ a $C^p$-cell in $N^{l+1}$ of dimension $k$
(resp., $k+1$) is a definable set of the form $\Gamma
(f)$ (the graph of $f$) where $f:C\into N$ is a $C^p$-definable
function and $C$ is a $C^p$-cell in $N^l$ of dimension $k$ (resp., of
the form $(f,g)_C:=\{(x,y)\in C\times N:f(x)<y<g(x)\}$ where
$f,g:C\into N$ are definable $C^p$-function with $-\infty \leq f<g\leq
+\infty $ and $C$ is a $C^p$-cell in $N^l$ of dimension $k$.  The
Euler characteristic $E(C)$ of a $C^p$-cell $C$ of dimension $k$ is 
defined to be $(-1)^k$.
}
\end{defn}

\begin{defn}\label{defn cell decomposition}
{\em
A $C^p$-cell decomposition of $N^m$ is a special kind of partition of
$N^m$ into finitely many $C^p$-cells: a partition of $N^1$ into
finitely many disjoint $C^p$-cells of dimension zero and one is a
$C^p$-cell decomposition of $N^m$ and, a partition of $N^{k+1}$ into
finitely many disjoint $C^p$-cells $C_1, \dots , C_m$ is a $C^p$-cell 
decomposition of $N^{k+1}$ if $\pi (C_1), \dots , \pi (C_m)$ is a
$C^p$-cell decomposition of $N^k$ (where $\pi:N^{k+1}\into N^k$ is the
projection map onto the first $k$ coordinates). Let $A_1,\dots
A_k\subseteq A\subseteq N^m$ be definable sets. A $C^p$-cell
decomposition of $A$ compatible with $A_1,\dots A_k$  is a finite
collection $C_1, \dots , C_l$ of $C^p$ partitioning $A$ obtained from
a $C^p$-cell decomposition of $N^m$ such that for every $(i, j)\in
\{1, \dots ,k\}\times \{1, \dots , l\}$ if $C_j\cap A_i\neq \emptyset $ 
then $C_j\subseteq A_i$. 
}
\end{defn}

\begin{fact}\label{fact cell decomposition}
\cite{vdd}
Given definable sets $A_1,\dots ,A_k$ $\subseteq A\subseteq N^m$ there
is a $C^p$-cell decomposition of $A$ compatible with $A_1,\dots A_k$
and, for every definable function $f:A\into N$, $A\subseteq N^m$,
there is a $C^p$-cell decomposition of $A$, such that each restriction 
$f_{|C}:C\into N$ is $C^p$ for each cell $C\subseteq A$ of the
$C^p$-cell decomposition.
\end{fact}

The o-minimal dimension $dim(A)$ and Euler characteristic $E(A)$ of a
definable set $A$ are defined by $dim(A)=\max \{dim (C):C\in
\mathcal{C}$$\}$ and $E(A)=\sum _{C\in \mathcal{C}}$$E(C)$ where
$\mathcal{C}$ is some (equivalently any) $C^p$-cell decomposition of
$A$. These notions are well behaved under the usual
set theoretic operations on definable sets, are invariant under
definable bijections and given a definable family of 
definable sets, the set of parameters whose fibre in the family has a 
fixed dimension (resp., Euler characteristic) is also a definable set.
The cell decomposition theorem is also used to
show that every definable set has only finitely many definably
connected components, and given a definable family of definable sets
there is a uniform bound on the number of definably connected
components of the fibres in the family.

\begin{nrmk}\label{nrmk homeomorphism a cells}
{\em
Let $A, B\subseteq N^m$ be definable sets and $\phi :A\into B$ a
definable homeomorphism. If $C_1, \dots , C_n$ is a cell decomposition
of $A$ then, $\phi (C_1), \dots , \phi (C_n)$ is a cell decomposition of $B$.
}
\end{nrmk}

\begin{defn}\label{defn properly bigvee definable cell decomposition}
{\em
Let $A_1, \dots , A_n, B, Z$ be properly $\bigvee $-definable subsets
of $X$. Let $I=\{1, 2, \dots \}$ be an enumeration
of $I$. Define inductively $(X'_i, M_i, N_i)$ for $i\in I$ by:
$X'_1=X_1$, $M_1$ is a cell decomposition of $\phi _1(X_1)$ compatible
with the definable subsets $\phi _1(X_1\cap X_j)$, $\phi _1(X_1\cap
A_l)$ and $\phi _1(X_1\cap X_j\cap A_l)$ for all $l\in \{1, \dots ,
n\}$ and all $j\in I_1$, and $N_1=M_1$ (recall that for $i\in I$,
$I_i:=\{j\in I:X_i\cap X_j\neq \emptyset \}$ is finite); let
$X'_{i+1}:=X_{i+1}\setminus \cup \{C:\exists r\in \{1, \dots , i\}, \,\,\phi
_r(C)\in N_i\}$ and $M_{i+1}$ is a cell decomposition of $\phi
_{i+1}(X'_{i+1})$ compatible with the definable sets $\phi
_{i+1}(X'_{i+1}\cap X_j)$, $\phi _{i+1}(X'_{i+1}\cap A_l)$ and $\phi
_{i+1}(X'_{i+1}\cap X_j\cap A_l)$ for all $l\in \{1, \dots , n\}$ and
all $j\in I_{i+1}$ and such that $N_{i+1}$ which is equal to $M_{i+1}$
together with all the cell in $N_j$ for 
$j\in I_{i+1}\cap \{1, \dots ,i\}$ is a cell decomposition of 
$\phi _{i+1}(X_{i+1})$. 

We define a properly $\bigvee $-definable cell decomposition of $X$
compatible with $A_1, \dots ,$$ A_n$ to be  a sequence $K=\{C:\phi _i
(C)\in N_i\,\,\, for\,\,\, some\,\,\, i\in I\}$ some $\{N_i: i\in I\}$
like above. By a properly $\bigvee $-definable cell decomposition of
$B$ compatible with $A_1, \dots , A_n$ we mean a sequence
$K_B=\{C\subseteq B: C\in K \}$ for some properly $\bigvee $-definable
cell decomposition $K$ of $X$ compatible with $A_1, \dots , A_n,
B$. Note that if $C\in K$ then, for all $j\in I$ if $C\cap X_j\neq
\emptyset$ then $C\subseteq X_j$ and $\phi _j(C)$ is a $k$-cell in
$\phi _j(X_j)$. If $C\subseteq Z$ we say that $C$ is a $k$-cell (of
$K$) in $Z$. 
}
\end{defn}

From fact \ref{fact cell decomposition} and definition 
\ref{defn properly bigvee definable cell decomposition} we get:

\begin{fact}\label{fact properly bigvee definable cell decomposition}
Given properly $\bigvee $-definable subsets $A_1,\dots ,A_k$ 
$\subseteq A\subseteq X$ there is a $C^p$-properly $\bigvee
$-definable cell decomposition of $A$ compatible with $A_1,\dots A_k$
and, for every strictly properly $\bigvee $-definable map $f:A\into N$, 
there is a $C^p$-properly $\bigvee $-definable
cell decomposition of $A$, such that each restriction 
$f_{|C}:C\into N$ is $C^p$ for each cell $C\subseteq A$ of the
$C^p$-properly $\bigvee $-definable cell decomposition.
\end{fact}

There is no $\bigvee $-definable cell decomposition of general 
$\bigvee $-definable subsets of $X$ and there is no corresponding 
$\bigvee $-definable cell decomposition theorem. 

If ${\mathbf X}$ is a properly $\bigvee $-definable but not locally
finite properly 
$\bigvee $-definable manifold, then there is no cell decomposition
theorem for general (properly) $\bigvee $-definable subsets of $X$,
however if $Z$ is a properly $\bigvee $-definable subset of $X$ for
which there is a subset $I'$ of $I^Z$ such that: $Z\subseteq \cup 
\{X_i:i\in I'\}$ and for all $i\in I'$, the set $\{j\in I':X_i\cap X_j
\neq \emptyset \}$ is finite, then since $Z$ is a properly 
$\bigvee $-definable subset of the locally finite 
properly $\bigvee $-definable
manifold ${\mathbf X'}$$=(X', (X'_i,\phi '_i)_{i\in I'})$ where
$X'=\cup \{X'_i:i\in I'\}$ and for each $i\in I'$, $X'_i=X_i$ and
$\phi '_i=\phi _i$,  $Z$ will have a properly $\bigvee $-definable
cell decomposition relative to ${\mathbf X'}$. Under these conditions, 
we will say that $Z$ is a {\it properly $\bigvee $-definable subset of 
$X$ with properly $\bigvee $-definable cell decomposition}. This fact,
will allow us to talk above notions  in properly $\bigvee $-definable 
manifolds that involve the properly $\bigvee $-definable cell
decomposition.

\medskip
When $\N$ expands a real closed field, then we can use the definable
triangulation theorem instead of the cell decomposition theorem. Below
we include the definition of properly $\bigvee $-definable
triangulation of $X$ compatible with finitely many properly $\bigvee
$-definable subsets. Similarly, all the notions that we define
using cell decomposition have an analogue obtained by using the
definable triangulation theorem, and moreover the two versions are compatible.

\begin{defn}\label{defn triangulaton}
{\em
Let $S_1,\dots , S_k$ $\subseteq S\subseteq N^m$ be definable sets. A
{\it definable triangulation} in $N^m$ of $S$ compatible with $S_1,\dots ,
S_k$ is a pair $(\Phi ,K)$ consisting of a complex $K$ in $N^m$ and a
definable homeomorphism $\Phi :S\into |K|$ such that each $S_i$ is a
union of elements of $\Phi ^{-1}(K)$. We say that $(\Phi ,K)$ is a
{\it stratified definable triangulation} of $S$ compatible with 
$S_1,\dots , S_k$ if: $m=0$ or $m>0$ and there is a stratified
definable triangulation $(\Psi ,L)$ of $\pi (S)$ compatible with $\pi
(S_1),\dots , \pi (S_k)$ (where $\pi :N^m\into N^{m-1}$ is the
projection onto the first $m-1$ coordinates) such that $\pi
_{|Vert(K)}:K\into L$ is a simplicial map and the diagram 
\[
\begin{array}{clcr}
\,\,S\,\,\,\stackrel{\Phi }{\rightarrow}\,\,\,|K|\\
\,\,\,\,\,\,\,\,\,\,{\downarrow}^{\pi }\,\,\,\,\,\,\,\,\,\,\,\,\,\, 
{\downarrow}^{\pi }\,\,\,\,\,\,\,\,\,\\
\pi (S)\stackrel{\Psi }{\rightarrow}|L|\\
\end{array}
\]
commutes.
We say that $(\Phi ,K)$ is a {\it quasi-stratified definable
triangulation} of $S$ compatible with $S_1,\dots , S_k$ if there is a
linear bijection $\alpha :N^m\into N^m$ such that $(\alpha \Phi \alpha
^{-1}, \alpha K)$ is a stratified definable triangulation of $\alpha
(S)$ compatible with $\alpha (S_1),\dots , \alpha (S_k)$.
}
\end{defn}

\begin{fact}\label{fact definable triangulation theorem}
\cite{vdd}
Let $S_1,\dots , S_k$ $\subseteq S\subseteq N^m$ be definable
sets. Then, there is a definable triangulation of $S$ compatible with 
$S_1,\dots , S_k$. Moreover, if $S$ is bounded then, there is a
quasi-stratified definable triangulation of $S$ compatible with 
$S_1,\dots , S_k$.
\end{fact}

\begin{defn}\label{defn properly bigvee definable triangulation}
{\em
Suppose that $\N$ expands a real closed field and let $A_1, \dots ,$
$A_n, B, Z$ be properly $\bigvee $-definable subsets of $X$.  Let
$I=\{1, 2, \dots \}$ be an enumeration of $I$. Define inductively $(X'_i,
(\Psi _i, M_i), (\Phi _i, N_i))$ for $i\in I$ by:  $X'_1=X_1$, $(\Psi
_1, M_1)$ is a definable triangulation of $\phi _1(X_1)$ compatible
with the definable subsets $\phi _1(X_1\cap X_j)$, $\phi _1(X_1\cap
A_l)$ and $\phi _1(X_1\cap X_j\cap A_l)$ for all $l\in \{1, \dots ,
n\}$ and all $j\in I_1$, and $(\Phi _1, N_1)=(\Psi _1, M_1)$; let
$X'_{i+1}:=X_{i+1}\setminus \cup \{C:\exists r\in \{1, \dots ,
i\},\,\, \phi _r(C)\in N_i\}$ and be 
$(\Psi _{i+1}, M_{i+1})$ be a definable triangulation of
$\phi _{i+1}(X'_{i+1})$ compatible with the definable sets $\phi
_{i+1}(X'_{i+1}\cap X_j)$, $\phi _{i+1}(X'_{i+1}\cap A_l)$ and $\phi
_{i+1}(X'_{i+1}\cap X_j\cap A_l)$ for all $l\in \{1, \dots , n\}$ and
all $j\in I_{i+1}$ and such that $(\Phi _{i+1}, N_{i+1})$ which is
equal to $(\Psi _{i+1}, M_{i+1})$ together with all the $(\Phi _{j|},
N_j)$ for $j\in I_{i+1}\cap \{1, \dots , i\}$ is a definable
triangulation of $\phi _{i+1}(X_{i+1})$. 

By a properly $\bigvee $-definable triangulation of $X$ compatible
with $A_1, \dots , A_n$ we mean a sequence $(\Phi , K)=\{(\Phi _i,
N_i): i\in I\}$ some $\{(\Phi _i, N_i): i\in I\}$ like above.
By a properly $\bigvee $-definable triangulation of $B$ compatible
with $A_1, \dots , A_n$ we mean a sequence $(\Phi _B, K_B)=\{(\Lambda
, L)\in (\Phi , K): |L|\subseteq B\}$ for some properly $\bigvee
$-definable triangulation $(\Phi , K)$ of $X$ compatible with $A_1, 
\dots , A_n, B$.  Note that for each $C\in K$ such that $\phi _i(C)\in
K_i$ for some $i\in I$, $C$ is definably homeomorphic to a
$k$-simplex. If $C\subseteq Z$ we say that $C$ is a $k$-simplex of
($(\Phi , K)$) in $Z$. 
}
\end{defn}

\end{subsection}

\begin{subsection}{Properly $\bigvee $-definable completeness}
\label{subsection completeness}

In this subsection ${\mathbf X}$ will be a properly $\bigvee $-definable 
manifold.
Recall that, by \cite{ps} a Hausdorff definable manifold ${\mathbf X}$
is called {\it definably compact} if for every definable continuous
map $\sigma:(a, b)\subseteq N\cup \{-\infty , +\infty \}\into X$ the
limits $\lim _{t\into a^{+}}\sigma (t)$ and $\lim _{t\into
b^{-}}\sigma (t)$ exist in $X$. Here we will introduce the notion of 
{\it properly $\bigvee $-definable completeness} for properly 
$\bigvee $-definable subsets of $X$, 
which coincides with the notion of definable compactness on Hausdorff
definable 
manifolds. The notion of properly $\bigvee $-definable completeness is 
invariant under taking elementary extensions, elementary substructures 
of $\N$ (containing the parameters over which ${\mathbf X}$ is
defined) and under taking expansions of $\N$ and reducts of $\N$ on
which ${\mathbf X}$ is defined.  In particular this shows that the
same holds for the notion of definable compactness of definable
manifolds, solving in this way a question raised (and solved for the
special case of affine definable manifolds) in \cite{ps}.  

\medskip
For a properly $\bigvee $-definable subset $Z$ of $X$, $\bar{Z}$
denotes the (topological) closure of $Z$ in $X$ and $\bar{Z}\setminus
Z$ its boundary in $X$. Both $\bar{Z}$ and $\bar{Z}\setminus Z$ are
properly $\bigvee $-definable subsets of $X$.
Note that, if $C\subseteq X_i$ is a cell in
$X$ then $\bar{C}$ is a disjoint union of cells in $X$ and
$\bar{\phi _i(C)}\subseteq N^m$ is a disjoint union of cells in
$N^m$. Below it we be convenient to work with cell in $[-\infty ,
+\infty ]^m$. These cells are defined in the same way as the cells in
$N^m$ with the convention that the new definable sets (and maps) are
the sets (and maps) definable in the structure obtained from $\N$ by
adding the constant symbols $-\infty $ and $+\infty $. Its easy to see
that the cell decomposition theorem also holds for this new structure
and a cell in $N^m$ is also a cell in $[-\infty , +\infty ]$.

\begin{defn}\label{defn definably complete}
{\em
Let $Z$ be  a properly $\bigvee $-definable subset of $X$.
We say that a cell $C$ in $Z$ is {\it definably complete (in $Z$)} if
for some (equivalently for all) $i\in I$ such that $C\subseteq X_i$,
for every cell of $\bar{\phi _i(C)}$ in $[-\infty , +\infty ]^m$ there
is a unique cell in $Z$ contained in $\bar{C}$ of the same dimension
and such that the incidence relations among the cells of 
$\bar{\phi _i(C)}$ in $[-\infty , +\infty ]^m$ are preserved under this
correspondence. We say that $Z$ is 
{\it properly $\bigvee $-definably complete} if for each $l\in I^Z$,
there is a subset $I^l$ of $I^Z$ such 
that  ${\mathbf X}^l$$=(X^l, (X^l_i,\phi ^l_i)_{i\in I^l})$ where $X^l=\cup 
\{X^l_i:i\in I^l\}$ and for each $i\in I^l$, $X^l_i=X_i$ and $\phi
^l_i=\phi _i$ is a locally finite
properly $\bigvee $-definable manifold and $Z^l:=\{Z_i:i\in I^l\}$ is
a properly $\bigvee $-definably  subset of $X^l$ with a properly 
$\bigvee $-definable cell decomposition $K_l$ in $X^l$
such that every cell $C$ of $K_l$ in $Z_l$ is definably complete in $Z^l$.
}
\end{defn}

\begin{thm}\label{thm definably complete}
If ${\mathbf X}$ is a Hausdorff properly $\bigvee $-definable manifold and
$Z$ is  a properly $\bigvee $-definable subset of $X$, then 
$Z$ is properly $\bigvee $-definably complete  iff for every definable 
continuous map $\sigma :(a,b)\subseteq
N\cup \{-\infty , +\infty \}\into Z$ the limits $\lim _{t\into
a^{+}}\sigma (t)$ and $\lim _{t\into b^{-}}\sigma (t)$ exist in
$Z$. In particular, the notion of properly $\bigvee $-definably
complete does not depend on the properly $\bigvee $-definable cell
decomposition $K_l$ of $Z_l$ ($l\in I^Z$) and is invariant under
taking elementary
extensions, elementary substructures of $\N$ (containing the
parameters over which $Z$ is defined) and under taking
expansions of $\N$ and reducts of $\N$ on which $Z$ is defined.
\end{thm}
 
\pf
By o-minimality, its enough to show that a cell $C\subseteq X_i$ 
in $X$ is definably complete in $X$ iff for every definable
continuous map $\sigma :(a,b)\subseteq N\cup \{-\infty , +\infty
\}\into C$ the limit $\lim _{t\into b^{-}}\sigma (t)$ exist in $\bar{C}$. 

By o-minimality, $\lim _{t\into b^{-}}\phi _i\circ \sigma (t)$ exists
in $\bar{\phi _i(C)}\subseteq [-\infty , +\infty ]^m$ (this notation
means the closure in $[-\infty , +\infty ]^m$) and so there is a cell
$B\subseteq \bar{\phi _i(C)}\subseteq [-\infty , +\infty ]^m$
containing this limit. If $B\subseteq \phi _i(X_i)$ then we are done,
otherwise there is a cell $D\subseteq \bar{C}$ in $K$ which
corresponds to $B$ and there is $j\in I$ such that $D\subseteq X_j$,
$C\subseteq X_i\cap X_j$ and $\lim _{t\into b^{-}}\phi _j\circ \sigma
(t) \in \phi _j(D)$. 

Suppose that $dim\phi _i(C)=k$. Then there is a open $k$-cell
$B\subseteq N^k$ definably homeomorphic to $\phi _i(C)$ and by
considering the definable sets $B\cap \pi _{\Sigma }^{-1}(x)$ for
$\Sigma =(i_1, \dots , i_{k-1})$, $\pi _{\Sigma }:N^k\into
N^{k-1}$ the projection onto the $\Sigma $-coordinates
and $x\in \pi _{\Sigma }(B)$, we get a family of continuous maps
$\sigma _{\Sigma , x, r}: (a_{\Sigma , x, r}, b_{\Sigma , x, r})\into
C$ such that $\cup \{\sigma _{\Sigma , x, r}(a_{\Sigma , x, r},
b_{\Sigma , x, r}): \Sigma , \,\,\, x\in \pi _{\Sigma }(B),$$ \,\,\,
r=1,\dots , R_{\Sigma , x}\}$$=C$. Using the fact that all these
definable continuous maps have limits in $X$ and the fact that the
definable continuous maps $\phi _i\circ \sigma _{\Sigma , x, r}$ have
limits in $\bar{\phi _i(C)}\subseteq [-\infty , +\infty ]^m$ the
result follows.        
\qed 

Note that, a properly $\bigvee $-definable subset $Z$ of $X$ is properly 
$\bigvee $-definably complete iff it is a (countable) union of
definably complete (i.e., definably compact) definable subsets of
$X$. We can use this to extend this notion to $\bigvee $-definable
subsets: we say that a $\bigvee $-definable subset $Z$ of $X$ is 
{\it $\bigvee $-definably complete} iff it is a (countable) union of 
definably complete definable subsets of $X$. However, 
the $\bigvee $-definable analogue of theorem \ref{thm definably
complete} only holds if
we assume  $\aleph _1$-saturation (consider again the case $\N$$=(\QQ$$,<)$). 

\begin{defn}\label{defn properly bigvee def compact}
{\em
A (properly) $\bigvee $-definable subset $Z$ of $X$ is said to be 
{\it (properly) $\bigvee $-definably compact} if $Z$ is (properly) 
$\bigvee $-definably complete and ${\mathbf X}$ is Hausdorff.
}
\end{defn}

\medskip
{\it We will from now on, assume that every properly 
$\bigvee $-definable manifold is Hausdorff.}
 
\end{subsection}

\begin{subsection}{Properly $\bigvee $-definable connectedness}
\label{subsection connectedness}

Here, ${\mathbf X}$ will be a properly $\bigvee $-definable manifold
and $Z$ a subset of $X$.

\begin{defn}\label{defn connected}
{\em 
If $Z$ is a properly $\bigvee $-definable subset of $X$. We say that
$Z$ is {\it properly $\bigvee $-definably connected} if there do not
exist two disjoint (relatively) open properly $\bigvee $-definable
subsets of $Z$ whose union is $Z$. Or equivalently, if there is no
open and closed proper properly $\bigvee $-definable subset of $Z$. 
}
\end{defn}

There is no interesting notion of $\bigvee $-definable connectedness
for general $\bigvee $-definable subsets of $X$: in $\N$$=(\QQ$$,<)$, $\QQ$ is
as a definable set definably connected, but is not a ``$\bigvee
$-definably connected'' set since it is a disjoint union of the open
$\bigvee $-definable subsets $(-\infty ,\sqrt{2})$ and $(\sqrt{2},
+\infty )$; in an $\aleph _1$-saturated extension
$\N$$=(N,<)$ of $(\RR$$,<)$, the convex hull $R$ of $\RR$ is a
proper, open and closed $\bigvee $-definable subset of the definable (and
therefore also $\bigvee $-definable) set $N$. Note also, that the
complement of $R$ in $N$ is not $\bigvee $-definable.

\begin{lem}\label{lem properly bigvee definable connected}
Let $Z$ be properly $\bigvee $-definable subset of $X$. Then $Z$ is a 
countable union of properly $\bigvee $-definably connected properly 
$\bigvee $-definable subsets.
\end{lem}

\pf
For each $i\in I^Z$, we have $Z_i=\cup \{Z_{i,j}:j=1,\dots ,m_i\}$
where for each $j\in \{1,\dots ,m_i\}$, $Z_{i,j}$ is a (closed and
open) definably connected component of $Z_i$. Let $S:=\{(i,j):i\in
I^Z,\,\,\, j=1,\dots ,m_i\}$. Clearly (by Zorn's lemma) there are
disjoint subsets $S_k$ of $S$ with $k\in K$ such that $S=\cup \{S_k:
k\in K\}$ and for each $k\in K$, $Z^k:=\cup \{Z_s:s\in S_k\}$ is a
properly $\bigvee $-definably connected component of $Z$. 
\qed

By o-minimility, if $\sigma :(a,b)\subseteq N\into X$ (with $a<b$) if a 
definable continuous injective map, then a total definable ordering 
$\leq _{\sigma }$ in $\sigma (a,b)$ is necessarily the ordering
induced by either the ordering of $(a,b)$ or the reverse ordering of 
$(a,b)$. If $\Sigma =(\sigma (a,b), \leq _{\sigma })$ is as above
and $\leq _{\sigma }$ is induced by the ordering of
(resp., the reverse ordering) of $(a,b)$ , we define the initial point 
$\inf \Sigma $ of $\Sigma $ to be the point $\sigma (a):=\lim_{t\into
a^+}\sigma (t)$ (resp., the point $\sigma (b):=\lim _{t\into
b^-}\sigma (t)$) and the final point of $\Sigma $ denoted $\sup \Sigma
$ to be the point $\sigma (b)$ (resp., $\sigma (a)$).

\begin{defn}\label{defn basic definable paths}
{\em
By a {\it basic definable path in $Z\subseteq X$}  we mean either the {\it
constant path $\epsilon _x$ at  $x\in Z$} or a definable totally
ordered set $\Sigma =(\Sigma , \leq _{\Sigma })$ contained in $Z$ such that
$\Sigma =\sigma ((a,b))$ for some continuous definable injective map
$\sigma :(a,b)\subseteq N\into X$ (with $-\infty <a<b<+\infty $) 
and $\inf \Sigma $, 
$\sup \Sigma \in Z$. We call $\sigma :(a,b)\into X$ a {\it definable 
parametrisation of $\Sigma$};  the {\it support of
$\Sigma $} denoted by $|\Sigma |$ is the closed set $\{\inf \Sigma
\}\cup \Sigma \cup \{\sup \Sigma \}$; the {\it inverse} of the basic 
definable path $\Sigma $ is the basic definable path $\Sigma ^{-1}:
=(\Sigma ,\leq _{\Sigma ^{-1}})$ where for all $s,t\in \Sigma $, 
$s\leq _{\Sigma ^{-1}}t$ iff $t\leq _{\Sigma }s$. 

For convenience,  if $\Sigma $ is the basic definable path $\epsilon
_x$, we call $\sigma :\emptyset \into X$ the definable parametrisation
of $\Sigma $; $x$ the support of $\Sigma$ (denoted as
well $|\Sigma |$); $x$ the initial and final point of $\Sigma$
(denoted as before $\inf \Sigma $ and $\sup \Sigma $); and $\Sigma
^{-1}=\Sigma $.
}
\end{defn}

\begin{defn}\label{defn definable paths}
{\em
A {\it definable path in $Z\subseteq X$ from $x$ to $y$} is a sequence $\Sigma 
=\Sigma _1\cdot \cdots \cdot \Sigma _l$ of basic definable paths 
$\Sigma _i$'s in $Z\subseteq X$ such that $\inf \Sigma _1=x$,
 $\sup \Sigma _l=y$ and for each $i=1,\dots , l-1$ we have  
$\sup \Sigma _i=\inf \Sigma _{i+1}$ . $x$ and $y$ are also called the 
initial and the final point of $\Sigma $ respectively (we use the
notation $x=\inf \Sigma $ and $y=\sup \Sigma $); if $\inf \Sigma =\sup 
\Sigma =z$ we say that $\Sigma $ is a {\it definable loop in
$Z\subseteq X$ at $z$}; the 
sequence $\{\sigma _i:i=1,\dots ,l\}$ where $\sigma _i$ is a definable 
parametrisation of $\Sigma _i$ is called a definable parametrisation
(of length $l$) of $\Sigma $; the {\it support} of $\Sigma $ 
denoted $|\Sigma |$ is $\cup \{|\Sigma _i|:i=1,\dots ,l\}$; the {\it
inverse} of $\Sigma $ is the definable path $\Sigma ^{-1}:= \Sigma
_l^{-1}\cdot \cdots \cdot \Sigma _1^{-1}$ in $Z\subseteq X$. 
We also have a natural
definition of the {\it product} $\Sigma \cdot \Gamma $ of two
definable paths $\Sigma $ and $\Gamma $ such that $\sup \Sigma =\inf \Gamma $.
}
\end{defn}

Finally, given two definable paths $\Sigma =\Sigma _1\cdot \cdots
\cdot \Sigma _l$ and $\Gamma =\Gamma _1\cdot \cdots \cdot \Gamma _k$
where $\Sigma _i$'s and $\Gamma _j$'s are basic definable paths, we
define $\Sigma \simeq \Gamma $ iff  $|\Sigma |=|\Gamma |$ and for all
$i=1,\dots ,l$ and $j=1,\dots ,k$ we have $(\Sigma _i\cap \Gamma _j,\leq
_{\Sigma _i})=$ $(\Sigma _i\cap \Gamma _j,\leq _{\Gamma _j})$. This is
an equivalence relation, and we will some times identify definable
paths under this equivalence relation. If $f:Z\subseteq X\into Y$ is a 
strictly properly $\bigvee $-definable map and $\Sigma $ is a
definable path in 
$Z\subseteq X$ from $x$ to $y$, then o-minimality implies that up to $\simeq $,
there a unique definable path $f\circ \Sigma $ in $f(Z)\subseteq Y$ 
from $f(x)$ to $f(y)$ with $|f\circ \Sigma |=f(|\Sigma |)$.

\begin{nrmk}\label{nrmk definable paths}
{\em
Of course if $\N$ expands an ordered group then, for every definable
path $\Sigma $ in $Z\subseteq X$ from $x$ to $y$, we have $|\Sigma |=\tau
([a,b])$ for some definable continuous map
$\tau :[a,b]\into X$. If in addition, $\N$ expands an ordered 
field then we can take $[a,b]=[0,1]$.
}
\end{nrmk}

\begin{defn}
{\em 
We say that $Z\subseteq X$ is {\it definably path connected}  if for
any two points 
$y,z\in Z$ there is a definable path in $Z$ with initial point $y$ and
final point  $z$. $Z\subseteq X$ is {\it locally definably path
connected} if every point in $Z$ has a definable open neighbourhood
(in $Z$) which is definably path connected. 
}
\end{defn}

Note that a properly $\bigvee $-definably connected properly 
$\bigvee $-definable subset $Z$ of $X$ is not necessarily definably 
path connected: take $\N$$=(\QQ$$,<)$, $X=\QQ$$^2$ and 
$Z=\{(x,y)\in \QQ$$^2:0<y<x\}\cup \{(0,0), (1,1)\}\subseteq X$; 
However, because the support of a definable path in $X$ is a definably 
connected definable subset of $X$, a definably path connected properly 
$\bigvee $-definable subset $Z$ of $X$ is properly 
$\bigvee $-definably connected. 

We now introduce two definitions that will play an important role in this
paper. First recall that if, $U\subseteq X$ and $Z\subseteq X$ be a definable 
subset of $X$ and a properly $\bigvee $-definable subset of $X$ respectively.
A {\it definable family of definable subsets of $U$} is simply a 
definable subset $F$ of $U\times Y$. Note that $W:=\pi (F)$ 
(where $\pi :X\times Y\into Y$ is the projection onto $Y$) is a 
definable subset of $Y$ and for each $w\in W$ the
fibre $F_w:=\{u\in U: (u,w)\in F\}$ is a definable 
subset of $U$. We use the notation $(F_w)_{w\in W}$.

\begin{defn}\label{defn bigvee definable with choice}
{\em
We say that a (properly) $\bigvee $-definable subset $Z$ of $X$ has 
{\it definable choice} if for every
definable family $(F_w)_{w\in W}$ of definable subsets of $Z$ there is
a definable map $s:W\into Z$ such that for every $w\in W$, $s(w)\in
F_w$. If moreover, for every
definable family $(F_w)_{w\in W}$ of definable subsets of $Z$ there is
a definable map $s:W\into Z$ such that for every $w, v\in W$, $s(w)\in
F_w$ and $F_w=F_v$ iff $s(w)=s(v)$ then we say that $Z$ has {\it
strong definable choice}. Note that, these notions are invariant under
taking elementary extensions or elementary substructures of $\N$ which
contain the parameter over which $Z$ is defined.
}
\end{defn}

The next result was proved in \cite{e1} in the definable case. 
The proof there can easily be adapted to get the following

\begin{fact}\label{fact choice on complete}
If $Z$ is a properly $\bigvee $-definably complete properly 
$\bigvee $-definable subset of $X$ and $(F_w)_{w\in W}$ is a definable 
family of definable closed subsets of $Z$, then there is a definable
map $s:W\into Z$ such that for all $u,v\in W$, $s(u)\in F_u$ and
$F_u=F_v$ iff $s(u)=s(v)$. 
\end{fact}

\begin{cor}\label{cor properly bigvee def compact and maps}
If $f:X\into Y$ is a strictly properly $\bigvee $-definable continuous
map and $Z$ is a definably compact definable subset of $X$, then
$f(Z)$ is a definably compact definable subset of $f(X)$. 
\end{cor}

\pf
If $\alpha :(a,b)\into f(Z)$ is a definable continuous map, then by
fact \ref{fact choice on complete}, there is a definable map 
$s:(a,b)\into Z$, such that for each $x\in (a,b)$, 
$s(x)\in f^{-1}(\alpha (x))$. Now since $Z$ is definably compact the
limit $\lim _{x\into a^+}s(x)$ (resp., $\lim _{x\into b^-}s(x)$) exist
in $Z$ and, since $f$ is continuous and $\alpha (x)=f(s(x))$, the
limit $\lim _{x\into a^+}\alpha (x)$ (resp., 
$\lim _{x\into b^-}\alpha (x)$) exist in $f(Z)$.
\qed

A {\it properly $\bigvee $-definable family of properly $\bigvee $-definable 
subsets of $Z$} is
a properly $\bigvee $-definable subset $F\subseteq Z\times Y$.
Note that, for each $(i,j)\in (I\times J)^F$, $F_{(i,j)}$ is a definable
subset of $F$ and $W_{(i,j)}:=\pi (F_{(i,j)})$ is a 
definable subset of $Y$. Consider the $\bigvee $-definable subset
$W:=\cup \{W_{(i,j)}:(i,j)\in (I\times J)^F\}$ of $Y$. Then for each 
$w\in W$ the fibre $F_w:=\{u\in Z: (u,w)\in F\}$ is a properly
$\bigvee $-definable subset of $Z$. We use the notation $(F_w)_{w\in W}$. 
If for each $w\in W$, $F_w$ is a definable subset of $Z$, then we say
that $(F_w)_{w\in W}$ is a {\it properly $\bigvee $-definable family
of definable subsets of $Z$}. Finally, if $Z$ is a
$\bigvee $-definable subset of $X$, by a {\it $\bigvee
$-definable family of $\bigvee $-definable (resp., properly $\bigvee
$-definable and definable) subsets of $Z$} we mean a $\bigvee
$-definable subset $F$ of $Z\times Y$ such that if
$w\in W$ where $W:=\pi (F)$
then $F_w:=\{u\in Z: (u,w)\in F\}$ is a $\bigvee $-definable (resp.,
properly $\bigvee $-definable and definable) subset of $Z$.

\begin{defn}\label{defn system of definable paths}
{\em
Suppose that $Z\subseteq X$ is a (properly) $\bigvee $-definable subset.
A {\it (properly) $\bigvee $-definable system of definable paths in $Z$} is 
a (properly) $\bigvee $-definable family $\Gamma \subseteq Z\times
Z\times Z$ such that for each $(u,v)\in Z\times Z$, the fibre $\Gamma
_{u,v}$ of $\Gamma $ at $(u,v)$ is a definable path in $Z$ from $u$ to $v$.
We will often use the notation $\{ \Gamma _{u,v}:u,v\in Z\}$ for a 
(properly) $\bigvee $-definable system of definable paths in $Z$.
}
\end{defn}

\begin{lem}\label{lem connected path connected}
Let $U\subseteq N^k$ be an open definable subset with definable choice and
let $B$ be a definable subset of $U$. Then 
$B$ is definably connected iff $B$ has a definable system of definable 
paths iff $B$ is definably path connected.
\end{lem}

\pf
So suppose that $B$ is definably connected.
We prove the result by induction on $k$. The case $k=1$ is
obvious. Suppose that the result holds in $N^{k-1}$. We now prove it
in $N^k$. First suppose that $B$ is a cell, and without loss of
generality we can assume that $B$ is an open cell in $U\subseteq N^k$ (for
otherwise, $dimB<k$ and $B$ is definably homeomorphic to an open cell
in $V\subseteq N^{dimB}$) with $V$ an open definable set definably
homomorphic to $U$. Then we have $B=(f,g)_C$ where $C$ is the projection
of $B$ in $N^k$ and 
$f,g:C\into N\cup \{-\infty ,+\infty \}$ are continuous definable 
functions such that $f<g$. Since $C$ is definably connected, it has a 
definable system $\{\Gamma ^C_{u,v}:u,v\in C\}$ of definable paths in
$C$. By definable choice and cell decomposition, there is a cell 
decomposition $C_1, \dots ,C_m$ of $C$ such that for $i=1,\dots , m-1$ 
either $C_i$ intersects the closure of $C_{i+1}$ or $C_{i+1}$
intersects the closure of $C_i$, together with definable continuous 
injective maps $\rho _i:C_i \into B$ such that $\pi \circ \rho
_i=1_{C_i}$ where $\pi :N^k\into N^{k-1}$ is the projection onto the
first $k-1$ coordinates. The cell decomposition $C_1,\dots ,C_m$ of
$C$ induces a cell decomposition $B_1,\dots ,B_m$ of $B$ such that for 
$i=1,\dots , m-1$ either $B_i$ intersects the closure of $B_{i+1}$ or 
$B_{i+1}$ intersects the closure of $B_i$. Taking the products of 
"vertical paths" together the definable path in the definable system 
$\{\rho_i\circ \Gamma ^C_{\rho _i(u),\rho _i(v)}:u,v\in C_i\}$ of
definable paths in $\rho _i(C_i)$, its clear that each $B_i$ has a
definable system of definable paths.  On the other hand, by definable
choice in $U$, there is a
definable path in $B_i\cup B_{i+1}$ connecting a point of $B_i$ with a
point of $B_{i+1}$ therefore,  $B$ has a definable system of definable paths.  

If $B$ is not a cell, then $B$ is a union of cells $C_1,\dots ,C_k$
where for each $i<k$, either $C_i$ intersects the closure of
$C_{i+1}$, or $C_{i+1}$ intersects the closure of $C_i$ and so, by definable
choice in $U$, there
is a definable path in $C_i\cup C_{i+1}$ connecting a point of $C_i$
with a point of $C_{i+1}$. Since the result holds for each $C_i$, it
also holds for $B$. 

Clearly, if $B$  has a definable system of definable paths, $B$ is
definably path connected. The above argument and cell decomposition
show that if $B$ is definably path connected then $B$ is definably connected.  
\qed

\begin{prop}\label{prop uniform path connected}
Suppose that $Z$ is a properly $\bigvee $-definable subset of $X$ with
definable choice. Then $Z$ is properly $\bigvee
$-definably connected iff $Z$ has a properly $\bigvee $-definable 
system of definable paths iff $Z$ is definably path connected. 
\end{prop}

\pf
For each $i\in I^Z$, $Z_i=\cup \{Z_{i,1}, \dots Z_{i,m_i}\}$ where
for each $j\in \{1,\dots ,m_i\}$, $Z_{i,j}$ is a definably connected
component of $Z_i$. Let $S:=\{(i,j):i\in I^Z,\,\, j=1,\dots ,m_i\}$. Then by
lemma \ref{lem connected path connected}, for each $s\in S$,
there is a definable system $\{ \Gamma ^s_{x,y}:x,y\in Z_s\}$
of definable paths in $Z_s$. Let $z\in Z$ and for each 
$s\in S$, let $z_s\in Z_s$. To finish its enough to find definable
paths $\Gamma _s$'s from $z$ to $z_s$. Let $R$ be the subset of $S$ of
all $s$'s for which $\Gamma _s$ exists. By the above, both 
$\bigcup _{s\in R}Z_s$ and $\bigcup _{s\in S\setminus R}Z_s$ are
disjoint open properly $\bigvee $-definable subsets of $Z$ whose union
is $Z$. Therefore, $Z$ is properly $\bigvee $-definably connected iff 
$R=S$ iff $Z$ has a properly $\bigvee $-definable system of definable
paths iff $Z$ is definably path connected.
\qed

\end{subsection}

\begin{subsection}{Dimension of $\bigvee $-definable subsets}
\label{subsection dimension of vee definable subsets}

The results in this subsection are obtained by easy modifications
of similar results from \cite{p1} for definable sets and definable maps.

\begin{lem}\label{lem dimension}
Let $Z\subseteq X$ be a properly $\bigvee $-definable subset of $X$
(over $B$) and define the dimension $dimZ$ of $Z$ to be
$\max \{dimZ_i:i\in I^Z\}$. Then the
following holds: (1) for $k\leq m$, $dimZ\geq k$ iff there is $i\in
I^Z$ such that $dimZ_i\geq k$ iff some
projection of $\phi _i(Z_i)$ onto $N^k$ has interior in $N^k$ iff
there is a definable equivalence relation on $Z_i$ (over $B$) infinitely many
equivalence classes of which have dimension $\geq k-1$ iff there is a
definable subset $U_i$ of $Z_i$ (over $B$) such that $dimU_i\geq
k-1$ and $U_i$ has no interior in $Z$ iff there is a properly 
$\bigvee $-definable subset $U$ of $Z$ such that $dimU\geq k-1$ and
$U$ has no interior in $Z$ iff there is a properly $\bigvee
$-definable equivalence relation on $Z$ with properly $\bigvee
$-definable classes and with infinitely many classes of dimension 
$\geq k-1$ on some $Z_i$;
(2) if $f:Z\into Y$ is a (strictly) properly $\bigvee $-definable injective
map then for each $i\in I^Z$, $dimZ_i=dimf(Z_i)$; 
(3) if $\{Z^s:s\in S\}$ with $|S|<\aleph _0$ is a collection of 
properly $\bigvee $-definable subsets of $X$ then, 
$dim(\cup \{Z^s:s\in S\})$$=\max \{dimZ^s:s\in S\}$ and (4) 
if $(F_w)_{w\in W}$ is a properly $\bigvee $-definable family of
(properly) $\bigvee $-definable subsets of
$Z$ (over $B$) then, the set $\{w\in W: dimF_w=k\}$ is a
$\bigvee $-definable subset of $W$ (over $B$). 
\end{lem}

\pf
$(1)$  The first equivalence follows from the definition; the second 
equivalence is lemma 1.4 in \cite{p1}; the third equivalence, is 
proposition 1.8 in \cite{p1}; the fourth equivalence is proposition
1.9 in \cite{p1}; the fifth equivalence follows from the previous ones
by taking $U:=\cup \{U_i:i\in I'\}$ where $I'$ is the subset of $I^Z$
of all $i$'s such that $dimZ_i\geq k$. Now suppose that there is a
properly $\bigvee $-definable subset $U$ of $Z$ such that $dimU\geq
k-1$ and $U$ has no interior in $Z$. By the above, we can assume
without loss of generality that for each $i\in I^U$, $dimU_i\geq k-1$
and $U_i$ has no interior in $Z$. Consider the following properly 
$\bigvee $-definable equivalence relation on $Z$ given by:
$x\thicksim y$ iff  $x,y\in \cup \{Z_i:i\in I^U\}$ and there are 
$\{i_1,\dots ,i_m\}\subseteq I^U$ and $\{x=x_1,x_2,\dots x_m=y\}$
$\subseteq \cup \{Z_i:i\in I^U\}$ such that for each $j=1,\dots ,m-1$, 
we have $x_j\thicksim _{i_j}x_j$ where $\thicksim _{i_j}$ is the 
definable equivalence relation on $Z_{i_j}$ given by the fourth 
equivalence, or otherwise.   Then, $\thicksim $ is a properly 
$\bigvee $-definable equivalence relation on $Z$ with properly 
$\bigvee $-definable classes such that for some $i\in I^Z$,
$\thicksim $ has infinitely many classes of dimension $\geq k-1$ on
some $Z_i$. The converse is immediate, by definition of $dimZ$ and previous 
equivalences.

Finally, $(2)$ and $(3)$ are lemma 1.5 in \cite{p1}, and $(4)$ follows
from lemma 1.6 in \cite{p1}. 
\qed

Recall that, given $A\subseteq N$, $a\in N$ is in the (model
theoretic) {\it algebraic closure of $A$} denoted $a\in acl(A)$, if
$a$ lies in a finite set definable over $A$, equivalently, $a$ is in
the {\it definable closure of $A$} denoted $a\in dcl(A)$ i.e., $\{a\}$
is definable over $A$. We obtain in this way an operator $acl(-):
\mathcal{P}$$(N)\into \mathcal{P}$$(N)$, where $\mathcal{P}$$(N)$ is
the set of all subsets of $N$. By the monotonicity theorem and uniform 
bounds, we get (see \cite{PiS1}) that $\N$ is a {\it geometric
structure} i.e., for any formula $\phi (x,y)$ there is $n\in \NN$ such
that for any $b\in N^l$, either $\phi (x,b)$ has less than $n$
solutions in $N^m$ or it has infinitely many, and $(N,acl(-))$ is a 
{\it pregeometry} (which means that $acl(-):\mathcal{P}$$(N)\into 
\mathcal{P}$$(N)$ satisfies the following: $(i)$ if $A\subseteq N$ then 
$A\subseteq acl(A)$, $(ii)$ if $A\subseteq N$ and $a\in acl(A)$ then, 
there is a finite $B\subseteq A$ such that $a\in acl(B)$, $(iii)$ if 
$A\subseteq N$ then $acl(acl(A))=acl(A)$ and $(iv)$ if $a,b\in N$, 
$A\subseteq N$ and $b\in acl(A,a)$ then either $b\in acl(A)$ or 
$a\in acl(A,b)$ (exchange principle). 

\begin{lem}\label{lem another dim}
For $B\subseteq N$ and $a\in
N^m$ let $dim(a/B):=\min \{|a'|:a'\subseteq a,\,\,\, a\subseteq 
acl(B\cup a')\}$. If $Z\subseteq X$ is a properly $\bigvee $-definable
subset of $X$ (over $B$) then 
$dim(Z)=\max \{dim(\phi _i(a)/B): a\in Z_i,\,\, i\in I^Z\}$.
\end{lem}

\pf
This follows from the definition of $dimZ$ and the corresponding
result for definable sets, see lemma 1.4 in \cite{p1}.
\qed

This approach gives a rise to a good notion of dimension for general 
$\bigvee $-definable subsets of $X$ when $\N$ is $\aleph _1$-saturated and
which is coherent with the definition of dimension for properly 
$\bigvee $-definable subsets (for details see \cite{p1} where the
definable case is treated): 

\begin{lem}\label{lem dimension again}
Suppose that $\N$ is $\aleph _1$-saturated and let $Z\subseteq X$ be a 
$\bigvee $-definable subset of $X$ (over $B$). Let
$dimZ$$:=\max \{dim(\phi _i(a)/B): a\in Z_i,\,\, i\in I^Z\}$. Then the
following holds: (1) for $k\leq m$, $dimZ\geq k$ iff there is a
definable subset $W$ of $Z$ such that $dimW\geq k$ iff there are $i\in
I^W$ such that $dimW_i\geq k$ iff some
projection of $\phi _i(W_i)$ onto $N^k$ has interior in $N^k$ iff
there is a definable equivalence relation on $W_i$ (over $B$) infinitely many
equivalence classes of which have dimension $\geq k-1$ iff there is a
definable subset $U_i$ of $W_i$ (over $B$) such that $dimU_i\geq
k-1$ and $U_i$ has no interior in $Z$ iff there is a 
$\bigvee $-definable subset $U$ of $Z$ such that $dimU\geq k-1$ and
$U$ has no interior in $Z$ iff there is a $\bigvee $-definable
equivalence relation on $Z$ with $\bigvee $-definable classes and 
with infinitely many classes of dimension $\geq k-1$ on some 
definable subset $W$ of $Z$;
(2) if $f:Z\into Y$ is a (strictly) $\bigvee $-definable injective
map then $dimZ=dimf(Z)$; 
(3) if $\{Z^s:s\in S\}$ with $|S|<\aleph _1$ is a collection of 
$\bigvee $-definable subsets of $X$ then, 
$dim(\cup \{Z^s:s\in S\})$$=\max \{dimZ^s:s\in S\}$ and (4)
if $(F_w)_{w\in W}$ is a $\bigvee $-definable family of $\bigvee
$-definable subsets of
$Z$ (over $B$) then, the set $\{w\in W: dimF_w=k\}$ is a
$\bigvee $-definable subset of $W$ (over $B$). 
\end{lem}

\pf
$(2)$ and $(3)$ are lemma 1.5 in \cite{p1}, and $(4)$ follows from 
lemma 1.6 in \cite{p1}. Using this and the fact that $Z=\cup
\{Z^s:s\in S\}$ with $|S|<\aleph _1$ and for each $s\in S$, $Z^s$ is a 
properly $\bigvee $-definable subset of $X$, $(1)$ follows from the 
corresponding result for properly $\bigvee $-definable subsets of $X$.
\qed

Given properly $\bigvee $-definable subsets $V\subseteq Z\subseteq X$, we say
that $V$ is {\it large in $Z$} if $dim(Z\setminus V)<dim Z$.
Lemma 1.6 in \cite{p1} implies the following result:

\begin{lem}\label{lem dim of fibers}
Let $Z\subseteq X$ be a properly $\bigvee $-definable subset of $X$
(over $B$). If $(F_w)_{w\in W}$ is a properly $\bigvee $-definable family of
properly $\bigvee $-definable subsets of
$Z$ (over $B$) then the set $\{w\in W:F_w \,\, is \,\, large \,\,
in\,\, Z\}$ is a $\bigvee $-definable subset of $W$ (over $B$).
\end{lem}

If $Z\subseteq X$ is a (properly) 
$\bigvee $-definable subset of $X$ (over $B$) and $a\in Z$ then, we
say that $a$ is a {\it generic point of $Z$ over $B$} if 
$dim(\phi _i(a)/B)=dim(Z)$ for some (for every) $i\in I^Z$ such
that $a\in Z_i$. Note that if $\N$ is $\aleph _1$-saturated, then
generic points of $Z$ over $B$ exist. 
The following is an easy consequence of the definitions:

\begin{lem}\label{lem large and generics}
Let $V\subseteq Z\subseteq X$ be properly $\bigvee $-definable
subsets. Then $V$ is large in $Z$ iff for every $B$  over which $V$,
$Z$ and $X$ are defined, every generic point of $Z$ over $B$ is in $V$.
\end{lem}

In general, given $\bigvee $-definable subsets $V\subseteq Z\subseteq
X$, we will say that $V$ is large in $Z$ iff for every $B$  over which $V$,
$Z$ and $X$ are defined, every generic point of $Z$ over $B$ is in
$V$. It easy to see that with this definition, we also get the
analogue of lemma \ref{lem dim of fibers} for $\bigvee $-definable subsets.

\end{subsection}

\end{section}

\begin{section}{The fundamental group}
\label{section fundamental group}

Throughout this section, $\mathbf{X}$$=(X,(X_i,\phi _i)_{i\in I})$ and 
$\mathbf{Y}$$=(Y,(Y_j,\psi _j)_{j\in J})$ will be properly $\bigvee $-
definably connected properly $\bigvee $-definable manifolds with definable
choice and of dimension $m$ and $n$ respectively. Note that, although
we work in this section in the category of properly $\bigvee $-definable 
manifolds with strictly (properly) $\bigvee $-definable
continuous maps, all the results we present here also hold in the
category of properly $\bigvee $-definable topological spaces with
corresponding strictly properly $\bigvee $-definable maps, where by a
properly $\bigvee $-definable topological space we mean a properly
$\bigvee $-definable subset of a properly $\bigvee $-definable
manifold with the induced topology. In fact, as one can easily check,
these results, except those from subsection \ref{subsection the
o-minimal tietze theorem} which use the properly $\bigvee $-definable
cell decomposition theorem, also hold under the assumption that $\N$
is $\aleph _1$-saturated, in the category of $\bigvee $-definable
topological spaces with corresponding strictly $\bigvee $-definable
continuous maps, where by a $\bigvee $-definable topological space we mean a 
$\bigvee $-definable subset of a properly $\bigvee $-definable
manifold with the induced topology.

\begin{subsection}{The fundamental group}
\label{subsection the fundamental group}

\begin{defn}\label{defn ordered numbering}
{\em
A $k$-cell ($k=0,1,2$) with an ordered numbering of the $0$-cells in
$X$ is a pair $(H, \zeta _H)$ where $H$ is a definably complete
$k$-cell in $X$ and $\zeta _H$ is a numbering of the $0$-cell of
$\bar{H}$ (the "corners of $\bar{H}$") which induces an orientation of
the boundary $\bar{H}\setminus H$ of $H$. Note that $|\zeta _H|$ is
the number of $0$-cells in $\bar{H}$.
}
\end{defn} 

Given $n\in \NN$, the function $[n]:\{1,\dots , 2n\}$
$\into \{1,\dots, n\}$ is defined by: $[n](m)=m$ if $m\leq n$ and 
$[n](m)=m-n$ otherwise.

\begin{defn}\label{special basic definable paths}
{\em
Let $(H, \zeta _H)$ be a $k$-cell ($k=0,1,2$) in $X_i$ with an ordered
numbering of the $0$-cells. The {\it special basic definable paths of
$(H,\zeta _H)$} are the basic definable paths $h_{i,j}$ in
$\bar{H}\setminus H$ joining the $0$-cell $i$ to the $0$-cell $j$
where $i=1,\dots , |\zeta _H|$ and $j=i$ or $j=[|\zeta _H|](i+1)$. The
{\it special definable paths of $(H, \zeta _H)$} are the following
definable paths: $(1)$ special basic definable path of $(H, \zeta
_H)$; $(2)$ $h_{k,k}\cdot h_{k,j}$ and $h_{k,j}\cdot h_{j,j}$; $(3)$
$h_{k, [|\zeta |](k+1)}\cdot h_{[|\zeta _H|](k+1), [|\zeta
_H|](k+2)}\cdot \cdots \cdot h_{[|\zeta _H|](k+j-1), [|\zeta
_H|](k+j)}$ for $k,j\in \{1,\dots , |\zeta _H|\}$; and $(4)$ $\Sigma
^{-1}$ where $\Sigma $ is a special definable path of $(H,\zeta _H)$.
}
\end{defn}
 
We say that two definable paths $\Sigma =\Sigma _1\cdot \cdots \cdot
\Sigma _l$ and $\Gamma =\Gamma _1\cdot \cdots \cdot \Gamma _l$ in $X$
are {\it adjacent} if for all $j=1,\dots , l$ $\Sigma _j$ and $\Gamma
_j$ are definable paths with  $|\Sigma _j|, |\Gamma _j|\subseteq
\bar{X_{j_i}}$. 

\begin{defn}\label{defn  definable pre-homotopy of paths}
{\em
Given two definable paths $\Sigma $ and $\Gamma $ in $X$, we define
$\Sigma \thicksim ^0\Gamma $ iff one of the following holds: $(1)$ there is 
$(H, \zeta _H)$ a $k$-cell ($k=0,1,2$) in $X_i$ with an ordered
numbering of the $0$-cells and there are special definable paths
$\Sigma '$ and $\Gamma '$ of $(H, \zeta _H)$ such that $\Sigma \simeq 
\Sigma '$, $\Gamma \simeq \Gamma '$, $\inf \Sigma ' =\inf \Gamma '$
and $\sup \Sigma '=\sup \Gamma '$, in this case we call $H$ a {\it
basic definable pre-homotopy} and we write $\Sigma \thicksim ^0_H\Gamma $;
$(2)$ there are two adjacent definable paths $\Sigma '=\Sigma _1\cdot
\cdots \cdot \Sigma _l$ and $\Gamma '=\Gamma _1\cdot \cdots \cdot
\Gamma _l$ in $X$ such that $\Sigma \simeq \Sigma '$ and $\Gamma
\simeq \Gamma '$, and there is a sequence $H=\{H_j:j=1,\dots, l\}$ of
basic definable pre-homotopies $\Sigma _j \thicksim ^0_{H_j}\Gamma _j$. 

If $(1)$ or $(2)$ holds, we write $\Sigma \thicksim ^0_H\Gamma $ and say that
$H$ is a {\it definable pre-homotopy} of $\Sigma $ and $\Gamma $.
}
\end{defn}

Finally we are ready to the define the notion of definable homotopy
between definable paths in $X$. Note that, by the definition of $\thicksim ^0$,
if $\Sigma \thicksim ^0\Gamma $ then we have $\inf \Sigma =\inf \Gamma $ and
$\sup \Sigma =\sup \Gamma $. 

\begin{defn}\label{defn definable homotopy of paths}
{\em
Given two definable paths $\Sigma $ and $\Gamma $ in $X$ we define 
$\Sigma \thicksim \Gamma $ iff one of the following holds: $(1)$ there is a
sequence of definable paths $\{\Sigma _j:j=1,\dots ,l\}$ such that  
$\Sigma _1\simeq \Sigma $, $\Sigma _l\simeq \Gamma $ and there is a sequence 
$H=\{H_j:j=1,\dots, l-1\}$ of definable pre-homotopies $\Sigma _j
\thicksim ^0_{H_j}\Sigma _{j+1}$, in this case we write 
$\Sigma \thicksim _H\Gamma $ and
say that $H$ is a {\it definable homotopy with fixed endpoints} of
$\Sigma $ and $\Gamma $; $(2)$ there are definable paths $\Lambda _1$
and $\Lambda _2$ such that 
$\Sigma \thicksim _H\Lambda _1\cdot \Gamma \cdot \Lambda
_2$ for some definable homotopy with fixed endpoints $H$.

If $(1)$ or $(2)$ holds, we write $\Sigma \thicksim _H\Gamma $ and say that $H$
is a {\it definable homotopy} of $\Sigma $ and $\Gamma $. 
}
\end{defn}

\begin{defn}\label{defn refiniment}
{\em
Suppose that $\Sigma \thicksim _B\Gamma $ with $B=\{B_i: i=1,\dots , n\}$, 
$\Sigma _i\thicksim _{B_i}\Sigma _{i+1}$ for $i=1, \dots , m-1$ and $\Sigma
\simeq \Sigma _1$ and $\Gamma \simeq \Sigma _m$ is a definable
homotopy. The {\it support} $|B|$ of the definable
homotopy $B$ is the union of the closure of all $k$-cells ($k=0,1,2$)
in $B$. Let $K$ be a finite collection of properly $\bigvee
$-definable sets in $X$. A {\it refinement $B'$ of $B$ respecting $K$}
is a definable homotopy $\Sigma \thicksim _{B'}\Gamma $ obtained by taking a
cell decomposition of $|B|$ compatible with the definable sets: 
$|\Sigma _i|$, all the cells in the closure of $B_i$ and $B_i\cap N$
for all $i\in \{1, \dots , n\}$ and all definable sets $N$ of $K$, and 
modifying the definable homotopy $B$ by adding the new $k$-cells
($k=0,1,2$), modifying the definable paths $\Sigma _i$ and enlarge the
list. We will often identify two definable homotopies $B$ and $B'$ under
the equivalence relation $B\simeq B'$ iff $B$ and $B'$ have a common 
refinement.
}
\end{defn}

If $\Sigma \thicksim _H\Gamma $ and $\inf \Sigma =\inf \Gamma $ we say that
$\Sigma $ and $\Gamma $ are definably homotopic with fixed initial
point denoted $\Sigma \thicksim _H\Gamma (rel\{\inf \})$. Similarly we define
definably homotopic with fixed final point and definably homotopic
with fixed end points and call $H$ a definable homotopy with fixed
final point, end points respectively and write $\Sigma \thicksim _H\Gamma (rel 
\{\sup \})$ and $\Sigma \thicksim _H\Gamma (rel \{\inf ,\sup \})$. As before, 
we sometimes omit the subscript $H$ in $\thicksim _H$. More generally and
similarly, given properly $\bigvee $-definable subsets $A$ and $B$ of
$X$, we can define the notion of relative definable homotopy of
definable paths $\Sigma $ and $\Gamma $ with initial point in $A$ and
final point in $B$, denoted $\Sigma \thicksim _H\Gamma (rel \{A,B\})$. It
follows easily from the definitions that all these relations are in
fact equivalence relations, we denote the set of equivalence classes
by $\pi _1(X,A,B)$ and the equivalence class of $\Sigma $ is denoted
by $[\Sigma ]$. As usual, we write $\pi _1(X,A)$ for $\pi _1(X,A,A)$
and $\pi _1(X,x)$ for $\pi _1(X,\{x\})$. 

\begin{lem}\label{lem fundamental group}
Let $\Sigma $ and $\Gamma $ be two definable paths in $X$ and let
$f:X\into Y$ be a strictly properly $\bigvee $-definable continuous
map. 
If $\Sigma \thicksim \Gamma $ then $f\circ \Sigma \thicksim f\circ \Gamma $. 
\end{lem}

\pf
The proof is of course by induction on the definition of $\thicksim $. There
result will therefore follow, if we prove it in the case where there
is $(H,\zeta _H)$ a $k$-cell ($k=0,1,2$) in $X_i$ with an ordered
numbering of the $0$-cells, and $\Sigma $ and $\Gamma $ are special
definable paths of $(H,\zeta _H)$ such that $\Sigma \thicksim ^0_H\Gamma $.

Let $J'$ be a finite subset of $J$ such that $f(H)\subseteq \{Y_j:j\in
J'\}$. If $k=0,1$ the result is obvious, so we can assume that $k=2$.
We now prove the result by induction on $|J'|$. 

Suppose that $|J'|=1$. Then using the fact that there is a definable 
homeomorphism $\rho :(g,h)_{(a,b)}\subseteq N^2 \into H$ where 
$g,h:(a,b)\into N$ are definable continuous maps such 
that $g<h$, and using cell decomposition its easy to see that there is
a cell decomposition $C$ of $f(H)$ such that $f\circ \Sigma \thicksim _Cf\circ 
\Gamma $: Let $K$ be a cell decomposition of $f(H)$, then there is a
cell decomposition $K'$ of $H$ such that for each cell $B\in K'$,
$f(B)$ is a cell in $f(H)$. Take $C$ to be the cell decomposition of
$f(H)$ compatible with $K$ and each cell $f(B)$ where $B\in K'$.  

For the case $|J'|>1$, we first show that there is a cell
decomposition $C=\{C_1,\dots ,C_q\}$ of $H$ such that $\Sigma
\thicksim _C\Gamma $ and for each $i=1, \dots , q$ there is 
$j\in J'$ such
that $C_i\subseteq Y_j$. Let $K=\{K_1,\dots , K_p\}$ be a cell 
decomposition of $H$ compatible with the definable sets $f^{-1}(f(H)
\cap Y_j)$ where $j\in J'$. Using $\rho $ and o-minimality we can
refine the cell decomposition $K$ if necessary, by adding points
$a=a_1<a_2<\cdots <a_t=b$ to $[a,b]$ to get a cell decomposition $C$
such that: each $0$-cell in $C$ is a $0$-cell of some $[\rho \circ
g(a_i), \rho \circ h(a_i)]$; each $1$-cell in $C$ is either a
subinterval of some $[\rho \circ g(a_i), \rho \circ h(a_i)]$ or the
graph of a definable continuous injective function $\rho \circ
k:(a_i,a_{i+1})\into H$; and each a $2$-cell in $C$ is of the form 
$\rho ((k,k')_{(a_i,a_{i+1})})$ for some 
definable continuous injective functions $k, k':(a_i, a_{i+1})\into 
\rho ^{-1}(H)$ such that $k<k'$. Its clear now that this (refined)
cell decomposition $C$ of $H$ satisfies the claim.

Therefore, we have reduced the case $|J'|>1$ to the case where, 
$\Sigma \thicksim _C\Gamma $ and each definable pre-homotopy $\Sigma
_i\thicksim ^0_{C_i}\Gamma _i$ which occurs here is such that $f(C_i)\subseteq
Y_j$ for some $j$. But then, induction and the case $|J'|=1$ proves
the result for $|J'|>1$. 
\qed

Given $A_1,\dots ,A_l\subseteq X'$ and  $B_1,\dots ,B_l\subseteq Y'$ 
properly $\bigvee $-definable subsets of $X$ and $Y$ respectively and
a continuous strictly properly $\bigvee $-definable map $f:X'\into Y'$ 
such that $f(A_i)\subseteq B_i$, we write $f:(X',A_1,\dots ,A_l)
\into (Y',B_1,\dots ,B_l)$ (if $l=1$, $A_1=\{x\}$ and $B_1=\{y\}$ we
write $f:(X',x)\into (Y',y)$). 

\begin{thm}\label{thm fundamental group}
$\pi _1$ is a covariant functor from the category of pointed properly 
$\bigvee $-definable manifolds into the category of groups. $\pi _1(X,x)$ is a
group called the definable fundamental group of $X$ at $x$
with the product defined by $[\Sigma ][\Gamma ]:=[\Sigma \cdot \Gamma ]$
and given strictly properly $\bigvee $-definable map $f:(X,x)\into
(Y,y)$, $\pi _1$$(f)$ (which is denoted $f_*$) is defined by
$f_*([\Sigma ])=[f\circ \Sigma ]$. Moreover,  $\pi _1(X\times Y, (x,y))\simeq 
\pi _1(X,x)\times \pi _1(Y,y)$ and, if there is a definable
path in $X$ from $x_0$ to $x_1$, then $\pi _1(X,x_0)$ $\simeq \pi _1(X, x_1)$.
\end{thm}

\pf
It follows easily from the definition of definable homotopy of
definable paths that $\pi _1(X,x)$ is in fact a well defined group
with identity $[\epsilon _x]$ and the inverse $[\Sigma ]^{-1}$ of
$[\Sigma ]$ given by $[\Sigma ^{-1}]$. The fact that $f_*:\pi
_1(X,x)\into \pi _1(Y,y)$ is well defined follows from lemma \ref{lem
fundamental group}, and $f_*$ is a homomorphism since $f\circ (\Sigma 
\cdot \Gamma )\simeq (f\circ \Sigma )\cdot (f\circ \Gamma )$.

Let $\Theta $ be a definable path from $x_0$ to $x_1$. Its easy to see
that the map $\Lambda :\pi _1(X,x_0)\into \pi _1(X,x_1)$ given by 
$\Lambda ([\Sigma ]):=[\Theta \cdot \Sigma \cdot \Theta ^{-1}]$ is a 
well defined isomorphism. 

The isomorphism $\pi _1(X\times Y, (x,y))\simeq 
\pi _1(X,x)\times \pi _1(Y,y)$  is easy to verify: let $\pi ^X:X\times
Y\into X$ and $\pi ^Y:X\times Y\into Y$ be the natural projections,
and let $i^X:X\into X\times Y$ and $i^Y:Y\into X\times Y$ be given by
$i^X(u):=(u,y)$ and $i^Y(v):=(x,v)$ respectively; let $\beta :\pi
_1((X\times Y),(x,y))\into \pi _1(X,x)\times \pi _1(Y,y)$ be given by
$\beta ([\Sigma ]):=(\pi ^X_*([\Sigma ]),\pi ^Y_*([\Sigma ]))$ and let
$\alpha :\pi _1(X,x)\times \pi _1(Y,y)\into \pi _1((X\times Y,(x,y))$
be given by $\alpha ([\Sigma ],[\Gamma ]):=i^X_*([\Sigma ])\cdot
i^Y_*([\Gamma ])$. Then clearly, $\beta \circ \alpha =
1_{\pi _1(X,x)\times \pi _1(Y,y)}$ and therefore, $\beta $ is
surjective. On the other hand, if $\beta ([\Sigma ])=\beta ([\Gamma
])$ and $\pi ^X\circ \Sigma \thicksim _H\pi ^X\circ \Gamma $ and $\pi
^Y\circ \Sigma \thicksim _K\pi ^Y\circ \Gamma $, then by considering
the definable set $(\pi ^X)^{-1}(|H|)\cap (\pi ^Y)^{-1}(|K|)$ and cell
decomposition, its easy to construct a definable homotopy $\Sigma
\thicksim \Gamma $, and therefore, $\beta $ is injective. 
\qed

\begin{defn}
{\em
The properly $\bigvee $-definably connected, properly $\bigvee $-definable
set $X$ is called {\it definably simply connected} if 
$\pi _1(X,x)=0$ for some (equivalently for all) $x\in X$.
}
\end{defn}
\end{subsection}

\begin{subsection}{Homotopy type}
\label{subsection homotopy type}

\begin{defn}\label{defn homotopy type}
{\em
Let $f,g:(X', A_1,\dots , A_l)\into (Y', B_1,\dots , B_l)$ be two
continuous strictly properly $\bigvee $-definable maps. We say that
$f$ and $g$ are {\it strictly  properly $\bigvee $-definably
pre-homotopic} if there is a continuous strictly properly $\bigvee
$-definable map $H:(X'\times |\Gamma |, A_1\times |\Gamma |,\dots ,
A_l\times |\Gamma |)\into (Y', B_1,\dots , B_l)$ such that $H(x,\inf
\Gamma )=f(x)$ and $H(x,\sup \Gamma )=g(x)$ for all $x\in X'$, where 
$\Gamma $ is a definable path in some $\bigvee $-definable manifold 
${\mathbf Z}$. $H$ is called a {\it strictly properly $\bigvee
$-definable pre-homotopy} between $f$ and $g$. We say that  $f$ and
$g$ are {\it strictly  properly $\bigvee $-definably homotopic} if
there is a sequence $f=h_0,h_1,\dots , h_m=g$ of continuous strictly 
properly $\bigvee $-definable maps from $(X',A_1,\dots ,A_l)$ into 
$(Y',B_1,\dots ,B_l)$ and a sequence $H=\{H_i:i=1,\dots , m\}$ of 
strictly properly $\bigvee $-definable pre-homotopies $H_i$ between 
$h_{i-1}$ and $h_i$.  $H$ is called a {\it strictly properly $\bigvee
$-definable homotopy} between $f$ and $g$.
}
\end{defn}

Strictly properly $\bigvee $-definable homotopy between continuous
strictly properly $\bigvee $-definable maps is an equivalence relation 
compatible with composition in the set of all such continuous
strictly properly $\bigvee $-definable maps. We denote by $[f]$ the 
equivalence class of $f:(X', A_1,\dots , A_l)\into (Y', B_1,\dots ,
B_l)$ and by $[(X', A_1,\dots , A_l),(Y', B_1,\dots , B_l)]$ the set
of all such classes. We say that
$(X', A_1,\dots , A_l)$ and $(Y', B_1,\dots ,B_l)$ 
have the same {\it strictly properly $\bigvee
$-definable homotopy type} if there are continuous strictly properly 
$\bigvee $-definable maps $f:(X', A_1,\dots , A_l)\into (Y', B_1,\dots
, B_l)$ and $g:(Y', B_1,\dots , B_l)\into (X', A_1,\dots , A_l)$ such
that $[g\circ f]=[1_{X'}]$ and $[f\circ g]=[1_{Y'}]$. This is an 
equivalence relation and we denote by $[(X', A_1,\dots , A_l)]$ the 
equivalence class of $(X', A_1,\dots , A_l)$. 

\begin{fact}\label{fact homotopy type}
Let $f,g:(X,x)\into (Y,y)$ and $f_i:(X,x)\into (Y,\{y_0,y_1\})$
($i=0,1$) with $f_i(x)=y_i$ be continuous strictly properly 
$\bigvee $-definable maps. 
If $[f]=[g]$, then $f_*=g_*$ and if $H$ is a strictly properly 
$\bigvee $-definable homotopy between $f_0$ and $f_1$ then, there is a 
definable path $\Gamma $ in $Y$ from $y_0$ to $y_1$ which determines a 
homomorphism $\Gamma _{\#}:\pi _1(Y,y_0)\into \pi _1(Y,y_1)$ such that 
$\Gamma _{\#}\circ f_{0*}=f_{1*}$. In particular,  if $f:(X,x)\into
(Y,y)$ is a strictly properly $\bigvee $-definable homotopy
equivalence then, the induced homomorphism $f_*:\pi _1(X,x)\into \pi
_1(Y,y)$ is an isomorphism.
\end{fact}

Of course if $\N$ expands a real closed field, then we can take in 
definition \ref{defn homotopy type} $|\Gamma |$ to be
$[0,1]$. Moreover we also have the following result:

\begin{prop}\label{prop two definitions of homotopy}
Suppose that $\N$ is an expansion of a real closed field. Then 
$\pi _1(X, A, B)=[([0,1],\{0\}, \{1\}), (X, A, B)]$.
\end{prop}

\pf
Let $\Sigma $ and $\Gamma $ be two definable paths in $X$.We need to
show that, $\Sigma \thicksim \Gamma $ iff there is a definable continuous
function $K:([0,1]\times [0,1], \{0\}\times [0,1], \{1\}\times [0,1])
\into (X, A, B)$ such that $\Sigma \simeq K_0$ and $\Gamma \simeq K_1$ 
where for $i=0,1$ $K_i$ is the definable path in $X$ parametrised by 
$t\into K(i,t)$. The implication from right to left follows from lemma 
\ref{lem fundamental group}. By the transitivity of the two notions of 
definable homotopy and by the inductive definition of $\thicksim $, to 
show the 
implication from left to right, its enough to assume that there is 
$(H,\zeta _H)$ a $k$-cell ($k=0,1,2$) in $X_i$ with an ordered
numbering of the $0$-cells, and $\Sigma $ and $\Gamma $ are special
definable paths of $(H,\zeta _H)$ such that $\Sigma \thicksim ^0_H\Gamma $. The
cases $k=0,1$ are easy to prove, for the case $k=2$ use first the fact
that the closure of a $2$-cell $(g,h)_{(a,b)}\subseteq N^2$ 
is definable homeomorphic to $[0,1]
\times [0,1]$. 
\qed

\end{subsection}

\begin{subsection}{The o-minimal Tietze theorem}
\label{subsection the o-minimal tietze theorem}

In this subsection, we generalise results proved in \cite{bo} for
definable sets in an o-minimal expansion of a real closed field. {\it
We assume here that $X$ is properly $\bigvee $-definably complete with
a properly $\bigvee $-definable cell decomposition}.

\begin{defn}\label{defn combinatorial pi}
{\em
Let $K$ be a properly $\bigvee $-definable cell decomposition of $X$
and $v$ a $0$-cell of $K$. An edge path of $K$ is a sequence $\sigma
=u_1, u_2, \dots , u_l$ of $0$-cells of $K$ such that for all $i=1,
\dots , l-1$, $u_i, u_{i+1}$ are $0$-cells of a $1$-cell of $K$. If 
$\gamma =w_1, w_2,\dots ,w_k$ is another edge path of $K$ and $u_l,
w_1$ are $0$-cells of a $1$-cell of $K$ then the concatenation $\sigma 
\cdot \gamma $ is also an edge path of $K$. $\sigma $ is an edge loop
of $K$ at $v$ if $v=u_1=u_l$. Note that an edge path $\sigma $ of $K$ 
determines uniquely a definable path $\Sigma $ in $X$ and the
concatenation of two edge paths of $K$ corresponds to the product of
the corresponding definable paths. $E(K,v)$ is the group under the
operation $[\sigma ][\gamma ]:=[\sigma \cdot \gamma ]$ of classes
$[\sigma ]$ of edge loops $\sigma $ of $K$ at $v$ under the
equivalence relation: $\sigma \thicksim \gamma $ iff 
$\Sigma \thicksim _H\Gamma $ where
every $k$-cell ($k=0,1,2$) in $H$ is a $k$-cell of $K$.
}
\end{defn}

Let $K$ be as above. A tree in $K$ is a collection of $1$-cells and 
$0$-cells  which is a tree; a maximal tree in $K$ necessarily contains 
all the $0$-cells of $K$.  A special (basic) definable path of $K$ is
a special (basic) definable path of some $k$-cell ($k=0,1,2$) $H$ of $K$.

\begin{lem}\label{lem combinatorial pi}
Let $K$ be a properly $\bigvee $-definable cell decomposition of $X$,
$v$ a $0$-cell of $K$ and $T$ a maximal tree in $K$. Then $E(K,v)$ is 
isomorphic to the group $G(K,T)$ generated by the special basic
definable paths of $K$ with relations: (1) $\Sigma =1$ if $\Sigma $
is a special basic definable paths of $K$ contained in $T$; and (2) 
$\Sigma =\Gamma $ if $\Sigma $ and $\Gamma $ are special definable
paths of $K$ and $\Sigma \thicksim ^0_H\Gamma $ for some $k$-cell ($k=0,1,2$)
$H$ of $K$.
\end{lem}

\pf
The isomorphism $\bar{\kappa }:G(K,T)\into E(K,v)$ is induced by the
map $\kappa $ defined as follows: Let $h_{a,b}$ be a special basic
definable path of $K$ where $a, b$ are $0$-cells of some $1$-cell of
$K$, then since $T$ is definable path connected, there is a unique
edge path $e$ in $T$ from $v$ to $a$ with no repetition of vertices,
with inverse $e^{-1}$ the unique edge path in $T$ from $a$ to $v$ with
no repetition of vertices. Define $\kappa (h_{a, b}):=e, a, b, a,
e^{-1}$. Let $\sigma $ be an edge loop of $K$ at $v$ then there is an
edge loop $\gamma =u_1, u_2, \dots , u_k$ of $K$ at $v$ with minimal
$k$ such that $[\sigma ] = [\gamma ]$. Then $[\kappa (h_{u_1,
u_2}\cdot h_{u_2, u_3}\cdots h_{u_{k-1}, u_k})]=[\sigma ]$. Similarly, 
its easy to see that $\bar{\kappa }$ is injective.
\qed

\begin{nrmk}\label{nrmk paths in a closed cell}.
{\em
Let $K$ be a properly $\bigvee $-definable cell decomposition of $X$
and let $C$ be a $l$-cell of $K$ in $X$. Then by an easy induction
argument on $l$ one can prove that: for any definable path 
$\Sigma $ in $\bar{C}$ there is a definable homotopy 
$\Sigma \thicksim _H\Gamma $ where $\Gamma $ is a definable path contained
in $\bar{C}$ obtained from edge paths of $K$ contained in $\bar{C}$; 
moreover, for any definable paths $\Sigma $ and $\Gamma $ obtained
from edge paths of $K$ contained in $\bar{C}$ (with the same
endpoints) there is a definable homotopy $\Sigma \thicksim _H\Gamma $ (with
fixed endpoints) such that every $k$-cell ($k=0,1,2$) of $H$ is a
$k$-cell of $K$.
}
\end{nrmk}

\begin{thm}\label{thm tietze}
(Tietze theorem).
Let $K$ be a properly $\bigvee $-definable cell decomposition of $X$
and $T$ a maximal tree in $K$. Then, $\pi _1(X,x)\simeq G(K,T)$. In 
particular, $\pi _1(X,x)$ is invariant under taking elementary
extensions, elementary substructures of $\N$ (containing the
parameters over which ${\mathbf X}$ is defined) and under taking
expansions of $\N$ and reducts of $\N$ on which ${\mathbf X}$ is
defined and $X$ has definable choice.  
\end{thm}

\pf
By theorem \ref{thm fundamental group} and lemma \ref{lem
combinatorial pi} its enough to show that: $(i)$ every definable loop
in $X$ at $v$ (a fixed $0$-cell of $K$) is definably homotopic to a
definable loop in $X$ at $v$ obtained from an edge loop of $K$ at $v$
and $(ii)$ if $\Sigma $ and $\Gamma $ are two definable loops in $X$ at
$v$ obtained from edge loops $\sigma $ and $\gamma $ of $K$ at $v$
such that $\Sigma $ and $\Gamma $ are definably homotopic then,
$\Sigma \thicksim _H\Gamma $ where every $k$-cell ($k=0,1,2$) in $H$ is a
$k$-cell of $K$.

Let $\Sigma $ be a definable loop in $X$ at $v$. Then clearly, there
are definable paths $\Sigma _1, \dots ,\Sigma _k$ and there are cells
$C_1, C_2, \dots , C_m$ of $K$ such that: $(1)$ $\Sigma \simeq \Sigma
_1\cdot \Sigma _2\cdots \Sigma _k$; $(2)$ for each $i=1, 2, \dots , k$ 
there is $l(i)\in \{1, \dots , m\}$ such that $|\Sigma _i |\subseteq 
C_{l(i)}$  and $(3)$ for each $l\in \{1, \dots , m-1\}$, $C_l$ is in
the closure of $C_{l+1}$ or $C_{l+1}$ is in the closure of $C_l$. Its
clear that, we can enlarge if necessary the list $\Sigma _1, \dots , 
\Sigma _k$  and the list $C_1, C_2, \dots , C_m$ preserving properties 
$(2)$ and  $(3)$ above, so that we have:
$(4)$ for every $i\in \{1, \dots , m\}$ every cell of $K$ in
$\bar{C_i}$ is in $\{C_1, \dots , C_m\}$; $(5)$ $\Sigma \thicksim \Gamma
_1\cdot \Gamma _2\cdots \Gamma _n$ where for each $j=1, \dots , n$ the 
definable path $\Gamma _j$ is such that the following holds: $(5a)$ $\Gamma _j$
$=\Sigma _{j_1}\cdot \Sigma _{j_2}\cdots \Sigma _{j_{n(j)}}$ where for 
each $i\in \{j_1, \dots , j_{n(j)}\}$ there is $l(i)\in \{1, \dots ,
m\}$ such that $|\Sigma _i |\subseteq C_{l(i)}$, there is $r(j)\in
\{j_1, \dots , j_{n(j)}\}$ such that the collection $C_{l(j_1)}, \dots
, C_{l(j_n)}$ is contained in the collection of all the cells of $K$
contained in $\bar{C_{r(j)}}$, and $C_{l(j_1)}$ and $C_{l(j_n)}$ are
$0$-cells.  To finish the proof of this part of the theorem, its
enough to show  that each definable path $\Gamma _j$ is definably 
homotopic to a special basic definable path of a $k$-cell ($k=0,1$) in
the list $C_1, \dots , C_m$. But this last claim is remark \ref{nrmk
paths in a closed cell}.

Let $\Sigma $ and $\Gamma $ be definable loops in $X$ at $v$ obtained
from edge loops $\sigma $ and $\gamma $ of $K$ at $v$, and suppose
that $\Sigma \thicksim _B\Gamma $ with $B=\{B_i: i=1,\dots , n\}$, $\Sigma
_i\thicksim _{B_i}\Sigma _{i+1}$ for $i=1, \dots , m-1$ and $\Sigma \simeq
\Sigma _1$ and $\Gamma \simeq \Sigma _m$. By taking if necessary a
refinement of $B$ respecting $K$, we can assume that, for each $i\in
\{1, \dots , n\}$ there are cells $C_{1^i}, \dots , C_{m^i}$ of $K$
satisfying condition $(4)$, each definable path $\Sigma _i$ satisfies
condition $(5)$ (with $\Sigma $ substituted by $\Sigma _i$, each
$\Gamma _j$ substituted by $\Gamma _{j^i}$, $n$ by $n^i$ and $j$ by
$j^i$) and $(6)$ for each $u\in \{1^i, \dots , m^i\}$ there is
$v(u)\in \{1^{i+1}, \dots , m^{i+1}\}$ such that $\bar{C_u}\cap
\bar{C_{v(u)}}\neq \emptyset $. Now the theorem follows from remark
\ref{nrmk paths in a closed cell}. 
\qed

\begin{thm}\label{thm van kampen}
(van Kampen theorem).
Let $X_0, X_1, X_2$ be closed properly $\bigvee $-definably connected,
properly $\bigvee $-definable subsets of $X$ with $X=X_1\cup X_2$ and
$X_0=X_1\cap X_2$. Let $x\in X_0$. Then for any group $G$ and any
homomorphisms $h_{\alpha }:\pi _1(X_{\alpha }, x)\into G$ ($\alpha
=0,1,2$) such that $h_0=h_{\alpha }\circ i_{\alpha }$ where $i_{\alpha
}:\pi _1(X_0, x)\into \pi _1(X_{\alpha }, x)$ is the homomorphism
induced by the inclusion, there is a unique homomorphism $h:\pi _1(X,
x)\into G$ such that $h_{\alpha }=h\circ j_{\alpha }$, where
$j_{\alpha }:\pi _1(X_{\alpha }, x)\into \pi _1(X, x)$ is the
homomorphism induced by the inclusion.
\end{thm} 

\pf
Take  a properly $\bigvee $-definable cell decomposition $K$ of $X$
compatible with $X_0, X_1$ and $X_2$. $K$ determines properly
$\bigvee $-definable cells decompositions $K_0, K_1$ and $K_2$
respectively. Take a maximal tree $T_0$ in $K_0$ and extend it to
maximal trees $T_1$ and $T_2$ in $K_1$ and $K_2$ respectively. Then,
$T:=T_1\cup T_2$ is a maximal tree in $K$. By theorem \ref{thm tietze}
we have $\pi _1(X, x)=G(K,T)$ and $\pi _1(X_{\alpha }, x)=G(K_{\alpha
}, T_{\alpha })$. The results follows since $i_{\alpha }(\pi _1(X_0,
x))$ is generated by the special basic definable paths of
$X_0\setminus T_0$ in $X_{\alpha }$. 
\qed

When $\N$ expands a real closed field, then we can use the definable
triangulation theorem instead of the cell decomposition theorem.

\begin{defn}\label{defn tietze}
{\em
Suppose that $\N$ is an o-minimal expansion of a real closed field and
let $(\Phi , K)$ be a properly $\bigvee $-definable triangulation of
$X$, $v$ a vertice of $K$ and $T$ a maximal tree in $K$ ($T$ contains
all the vertices of $K$). An edge path of $K$ is a sequence $\sigma
=u_1, u_2, \dots , u_l$ of vertices of $K$ such that for all $i=1,
\dots , l-1$, $u_i, u_{i+1}$ are vertices of an edge of $K$. If
$\gamma =w_1, w_2, \dots , w_k$ is another edge path of $K$ and $u_l,
w_1$ are vertices of an edge of $K$ then the concatenation $\sigma
\cdot \gamma $ is also an edge path of $K$. $\sigma $ is an edge loop
of $K$ at $v$ if $v=u_1=u_l$. $E(K,v)$ is the group under the
operation $[\sigma ][\gamma ]:=[\sigma \cdot \gamma ]$ of classes
$[\sigma ]$ of edge loops $\sigma $ of $K$ at $v$ under the
equivalence relation: $v, u, a, b, c,\dots , w, v \thicksim  v, u, a, c, \dots
, w, v$ iff $abc$ spans a $k$-simplex of $K$ where $k=0,1,2$. $E(K,v)$
is isomorphic to the group $G(K,T)$ generated by the $1$-simplexes of
$K$, denoted $g_{a,b}$ for each edge $a,b$ and with relations: $(1)$
$g_{a,b}=1$ if $a,b$ spans a simplex of $L$ and $(2)$
$g_{a,b}g_{b,c}=g_{a,c}$ if $a,b,c$ spans a simplex of $K$.
}
\end{defn}

\begin{cor}\label{cor tietze}
(Tietze theorem).
Suppose that $\N$ is an o-minimal expansion of a real closed field and
let $(\Phi , K)$ be a properly $\bigvee $-definable triangulation of
$X$ and $T$ a maximal tree in $K$. Then, $\pi _1(X,x)\simeq
G(K,T)$. In particular, $\pi _1(X,x)$ is invariant under taking
elementary extensions, elementary substructures of $\N$ (containing
the parameters over which ${\mathbf X}$ is defined) and under taking
expansions of $\N$ and reducts of $\N$ on which ${\mathbf X}$ is
defined and $X$ has definable choice.  
\end{cor}

\end{subsection}
\end{section}

\begin{section}{$\bigvee$-definable covering spaces}
\label{section covering spaces}

For the rest of this section an less otherwise stated,
$\mathbf{X}$$=(X,(X_i,\phi _i)_{i\in I})$,$\mathbf{Y}$$=(Y,(Y_j,\psi
_j)_{j\in J})$ and $\mathbf{Z}$$=(Z,(Z_k,\tau _k)_{k\in K})$ will be
properly $\bigvee $-definably connected, properly 
$\bigvee $-definable manifolds with definable choice. As in the last
section, the results of this section also hold in the category of
properly $\bigvee $-definable topological spaces with strictly
properly $\bigvee $-definable continuous maps.

\begin{subsection}{Strictly properly $\bigvee $-definable covering spaces}
\label{subsection covering spaces}

\begin{defn}\label{defn covering space}
{\em 
$(Y,p,X)$ is called a {\it strictly properly $\bigvee $-definable covering
space} if the strictly properly $\bigvee $-definable map $p:Y\into X$ is
continuous, surjective and there is a cover $\{U_l: l\in L\}$ of $X$
with $|L|<\aleph _1 $ such that: $(1)$ for each $l\in L$, 
$U_l$ is an open definably
connected definable subset of $X$; $(2)$ for each $i\in I$, there a
finite subset $L_i$ of $L$ such that $X_i\subseteq \cup \{U_l:l\in
L_i\}$; and $(3)$ for each
$l\in L$, $p^{-1}(U_l)$ is a disjoint union of open definable subsets
of $Y$, each of which is mapped homeomorphically by $p$ onto
$U_l$. We say that $p:Y\into X$ is a {\it strictly properly 
$\bigvee $-definable covering map} and $\{U_l:l\in L\}$ is called a 
{\it $p$-admissible family of definable neighbourhoods}. 
}
\end{defn}

If $(Y, p,X)$ is a strictly properly $\bigvee $-definable covering space then:
$(1)$ for every definable subset $Z$ of $X$, $p^{-1}(Z)$ is a properly
$\bigvee $-definable subset of $Y$, and in particular, for each $x\in
X$, $p^{-1}(x)$ is a properly $\bigvee $-definable discrete subset of $Y$ and
therefore, $|p^{-1}(x)|<\aleph _1$; $(2)$ if for each $l\in L$,
$p^{-1}(U_l)$ is a disjoint union of the open definable subsets
$V_{l,s}$ of $Y$ with $s\in S_l$, each of which is mapped
homeomorphically by $p$ onto $U_l$ then, since for each $j\in J$, the
map $p_{|Y_j}:Y_j\into X$ is definable, there is a finite subset $Q_j$
of $\{(l,s):l\in L,\,\,\, s\in S_l\}$ such that $Y_j\subseteq 
\cup \{V_q:q\in Q_j\}$; $(3)$ let $i\in I$, then $p^{-1}(X_i)$ is an
open properly $\bigvee $-definable subset of $Y$ and considering the
open definable subsets $V_{i,s}$ of $Y$ where $s\in \cup \{S_l:l\in
L^i\}$ and $L^i$ is a finite subset of $L$ such that $X_i\subseteq
\cup \{U_l:i\in L^i\}$, we see that there are open definable subsets 
$Y_{i,r}$ of $p^{-1}(X_i)$ for $r\in R_i$ with $|R_i|<\aleph _1$ such 
that: $(i)$ $p^{-1}(X_i)=\cup \{Y_{i,r}:r\in R_i\}$; $(ii)$   for each 
$r\in R_i$, the set $\{s\in R_i:Y_{i,s}\cap Y_{i,r}\neq \emptyset \}$
is finite; and $(iii)$ for each $r\in R_i$, $p_{|Y_{i,r}}:Y_{i,r}\into
X_i$ is a definable surjective map. The collection $\{Y_{i,r}:i\in
I,\,\, r\in R_i\}$ of open definable subsets of $Y$ determine in a
natural way a $\bigvee $-definable manifold ${\mathbf Y'}$ which we
can identify with ${\mathbf Y}$ (they are strictly properly 
$\bigvee $-definably isomorphic), and under this identification 
(which we will be considering from now on), its easy to see that 
${\mathbf X}$ is a properly $\bigvee $-definable manifold iff
${\mathbf Y}$ is a properly $\bigvee $-definable manifold; $(4)$ the  properly
$\bigvee $-definable covering map $p:Y\into X$ is an open surjection:
Let $V$ be an open $\bigvee $-definable subset of $Y$ and let $x\in
p(V)$. Let $U$ be a $p$- admissible definable neighbourhood of $x$,
$y\in p^{-1}(x)\cap V$ and let $W$ be the definable sheet over $U$
containing $y$. Then $W\cap V$ is an open definable subset of $V$
containing $y$, $p(W\cap V)$ is an open definable subset of $U$
containing $x$ and therefore, $p(V)$ is open. 

\begin{defn}\label{defn isomorphisms of coverings}
{\em
We define 
a strictly properly $\bigvee $-definable isomorphism between strictly 
properly $\bigvee $-definable covering spaces  
$(Y,p, X)$ and $(Z, q, X)$ to be  a strictly properly
$\bigvee $-definable homeomorphism $\phi :Y\into Z$ such that $q\circ
\phi =p$.  The {\it group of strictly properly $\bigvee $-definable covering
transformations}, $Cov(Y/X)$ is the group of all strictly properly $\bigvee
$-definable homeomorphisms $\phi :Y\into Y$ such that $p\circ \phi =p$. 
}
\end{defn}

Pointed strictly properly $\bigvee $-definable covering space
$((Y,y),p,(X,x))$ and the group $Cov((Y,y)/(X,x))$ of pointed
strictly properly $\bigvee $-definable covering transformations are 
define in the obvious way. 
  
\begin{nrmk}\label{nrmk cells in covers}
{\em
Let $(Y,p,X)$ be a strictly properly $\bigvee $-definable covering
space and let $W$ be a properly $\bigvee $-definable subset of $X$ 
such that $W$ or $p^{-1}(W)$ has a properly $\bigvee $-definable cell 
decomposition.  Then there are properly $\bigvee $-definable cell
decompositions $L$ and $K$ of $p^{-1}(W)$ and $W$ respectively such that for
every cell $C$ of $L$, $p(C)$ is a cell of $K$ and $p_{|C}:C\into
p(C)$ is a definable homeomorphism. In particular, $p^{-1}(W)$ is properly
$\bigvee $-definably complete iff $W$ is properly $\bigvee $-definably
complete. 
}
\end{nrmk}

\begin{lem}\label{lem unique lift} 
Let $p:Y\into X$ be a strictly properly $\bigvee $-definable covering
map. Suppose that $f,g:Z\into Y$ are continuous strictly properly 
$\bigvee $-definable maps such that $p\circ f$$=p\circ g$. If $f(z)=g(z)$ for
some point $z\in Z$, then $f=g$. 
\end{lem}

\pf
Let $W:=\{w\in Z: f(w)=g(w)\}$. Then $W$ and $Z\setminus W$ are
properly $\bigvee $-definable subsets of $Z$ (since $f$ and $g$ are
strictly properly $\bigvee $-definable maps. By properly $\bigvee
$-definable connectedness of $Z$ its enough to show that $W$ is open
and closed in $Z$. Let $w\in W$ and let $U$ be a $p$-admissible
definable open neighbourhood of $p(f(w))=p(g(w))$, and let $V$ be the
definable sheet over $U$ containing $f(w)=g(w)$. Clearly,
$O=f^{-1}(V)\cap g^{-1}(V)\cap Z_k$ for $Z_k$ such that $w\in Z_k$, is
a definable open neighbourhood of $w$ in $Z$. We claim that
$O\subseteq W$. If $v\in O$, then $f(v),g(v)\in V$ and also
$p(f(v))=p(g(v))$, but $p_{|V}$ is a definable homeomorphism and so
$f(v)=g(v)$. Therefore, $W$ is open. Since $Z$ is Hausdorff, $W$ is
closed. This assumption on $Z$ is not necessary: if  $w\in Z\setminus
W$, let $V$ be a definable $p$-admissible open neighbourhood of
$p(f(w))$. If both $f(w)$ and $g(w)$ lie in the same sheet over $V$,
then the argument above shows that $p(f(w))=p(g(w))$.  Therefore,
$f(w)\in S$ and $g(w)\in S'$, where $S$, $S'$ are distinct sheets. But
$U=f^{-1}(S)\cap g^{-1}(S')\cap Z_k$ is a definable open neighbourhood
of $w$ such that $U\subseteq Z\setminus W$. Therefore, $Z\setminus W$
is open as well.
\qed

Recall that some times we identify definable homotopies under the
equivalence $\simeq $ of having a common refinement.
 
\begin{prop}\label{prop lifting paths}
Let $p:Y\into X$ be a strictly properly $\bigvee $-definable covering map. (1)
If $\Gamma $ a definable path in $X$ and $y\in Y$ is such that
$p(y)=\inf \Gamma $, then there is a unique definable path
$\bar{\Gamma }$ in $Y$ such that $y=\inf \bar{\Gamma }$ and $p\circ
\bar{\Gamma } = \Gamma $.

(2) Suppose that $\Gamma \thicksim _H\Sigma $ is a definable homotopy of
definable paths in $X$. Let $\bar{\Gamma }$ be a definable lifting of
$\Gamma $, then there is a unique definable lifting $\bar{H}$ of $H$
(i.e., $p\circ \bar{H}=H$) such that $\bar{\Gamma
}\thicksim _{\bar{H}}\bar{\Sigma }$ where $\bar{\Sigma }$ is a definable
lifting of $\Sigma $. 
\end{prop}

\pf
Let $\{U_l:l\in L\}$ a $p$-admissible family of definable
neighbourhoods. (1) Since for each $i\in I$, there is a finite subset
$L_i$ of $L$ such that $X_i\subseteq \cup \{U_l:l\in L_i\}$ and
$|\Gamma |$ is a definable subset of $X$, there is a finite subset
$I'$ of $I$ such that $|\Gamma |\subseteq \cup \{X_i:i\in I'\}$ and so 
$|\Gamma |\subseteq \cup \{U_l:l\in \cup \{L_i:i\in I'\}\}$. 
Therefore $\Gamma \simeq \Gamma _1\cdot \cdots \cdot \Gamma _n$ 
for some definable paths $\Gamma _j$ ($j=1, \dots , n$) such that for
each $j=1,\dots ,n$ there is $j(l)\in \cup \{L_i:i\in I'\}$ such that 
$|\Gamma _j|\subseteq U_{j(l)}$. Since the result clearly holds for
each $\Gamma _j$ the result holds for $\Gamma $. 
(2) This is proved in a similar way, by taking a refinement of $H$
compatible with $\{U_l:l\in L'\}$ where $L'$ is a finite subset of
$L$ such that $|H|\subset \cup \{U_l:l\in L'\}$. 
\qed

\textbf{Notation:} Referring to proposition \ref{prop lifting paths}, 
we denote by $y*\Gamma $ the final point $\sup \bar{\Gamma }$
of the definable lifting $\bar{\Gamma }$ of $\Gamma $ with initial
point $\inf \bar{\Gamma }=y$.

\begin{cor}\label{cor lifts with the same end point}
Let $((Y,y),p,(X,x))$ be a  strictly properly $\bigvee$-definable covering
space. Then (1) the induced homomorphism $p_*:\pi _1$$(Y,y)\into
\pi_1$$(X,x)$ is injective; (2) if $\Sigma $ is a definable loop at
$x$, then $y=y*\Sigma $ iff $[\Sigma ]\in p_*(\pi _1(Y,y))$ and if
$\Lambda $ and $\Lambda '$ are two definable paths in $X$ from $x$ to
$x'$, then $y*\Lambda =y*\Lambda '$ iff 
$[\Lambda \cdot \Lambda '^{-1}]\in p_*(\pi _1(Y,y))$.
\end{cor}

\pf
(1) Let $\Sigma $ be a definable loop at $y$ such that $p_*([\Sigma
])=1$. And let $H$ be the definable homotopy from $p\circ \Sigma $ to
$\epsilon _x$. By definable homotopy lifting, $H$ lifts to a definable
homotopy $\bar{H}$ from $\Sigma $ to some definable path. By
uniqueness of $\bar{H}$ its easy to see that this definable path is
$\epsilon _y$, and $[\Sigma ]=1$ as required. A similar argument shows (2). 
\qed

\begin{prop}\label{prop size of fibres}
Let $p:(Y,y)\into (X,x)$ be a strictly properly $\bigvee$-definable covering
map.  (1) The map $\rho :\pi _1(X,x)\times p^{-1}(x)\into p^{-1}(x)$
given by $\rho ([\Sigma ],u):= u*\Sigma $ is a transitive action, the
stabiliser of $u\in p^{-1}(x)$ is $p _*(\pi _1(Y,u))$ and
$|p^{-1}(x)|=[\pi _1(X,x):p _*(\pi _1(Y,u))]$.

(2) For all $x_1, x_2\in X$, $|p^{-1}(x_1)|=|p^{-1}(x_2)|$, for all
$u,v\in p^{-1}(x)$, $p _*(\pi _1(Y,u))$ and $p _*(\pi _1(Y,v))$ are
conjugate subgroups of $\pi _1(X,x)$ and if $S$ is a subgroup of $\pi
_1(X,x)$ that is conjugate to $p _*(\pi _1(Y,u))$ then there is $v\in
p^{-1}(x)$ such that $S=p _*(\pi _1(Y,v))$.
\end{prop}

\pf
(1) Note that by proposition \ref{prop lifting paths} (2), $\rho $ is
well defined. We have $\rho ([\epsilon _x],u)=u$ because the lifting
$\bar{\epsilon _x}$ of $\epsilon _x$ at $u$ is $\epsilon _u$. Now
suppose that $[\Lambda ]\in \pi _1(X,x)$. Let $\bar{\Sigma }$ be the
definable lifting of $\Sigma $ at $u$ and let $\bar{\Lambda }$ be the
definable lifting of $\Lambda $ at $u*\Sigma $, then $\bar{\Sigma
}\cdot \bar{\Lambda }$ is the definable lifting of $\Sigma \cdot
\Lambda $ that begins at $u$ and ends at $(u*\Sigma )*\Lambda
$. Therefore, $\rho ([\Sigma ][\Lambda],u)=$$\rho ([\Sigma \cdot
\Lambda ],u)=\rho ([\Lambda ],\rho ([\Sigma ],u))$. Let $u,v\in
p^{-1}(x)$. Since $Y$ is definably path connected, there is a definable
path $\Gamma $ in $Y$ from $u$ to $v$. $p\circ \Gamma $ is a loop at
$x$ whose lifting at $x$ is $\Gamma $. Thus $[p\circ \Gamma ]\in \pi
_1(X,x)$, and $u*(p\circ \Gamma )=v$ i.e., $\rho ([p\circ \Gamma
],u)=v$ and $\rho $ is transitive.

Let $[\Sigma ]\in \pi _1(X,x)_{u}$ (the stabilizer of $u$) and let
$\bar{\Sigma }$ be the definable lifting of $\Sigma $ at $u$. Then
$u=u*\Sigma=\sup \bar{\Sigma }$ and so $[\bar{\Sigma }]\in \pi _1(Y,u)$,
$[\Sigma ]=[p\circ \bar{\Sigma }]\in p_*(\pi _1(Y,u)).$ For the
reverse inclusion, suppose that $[\Sigma ]=[p\circ \Gamma ]$ for some
$[\Gamma ]\in \pi _1(Y,u).$ Then $\bar{\Sigma }=\Gamma $. Therefore, 
$u*\Sigma =\sup \bar{\Sigma }$$=\sup \Gamma =u$, and $[\Sigma ]\in \pi _1(
X,x)_u.$ The rest of (1) follows from the theory of $G$-sets.

(2) We have the following commutative diagram
\[
\begin{array}{clrc}
\pi _1(Y,u_1)\,\,\,\,\stackrel{A}{\rightarrow }\,\,\,\,
\pi _1(Y,u_2)\\
\downarrow ^{p_*}\,\,\,\,\,\,\,\,\,\,\,\,\,\,\,\,\,\,\,\,\,\,\,
\downarrow ^{p_*}\\
\pi _1(X,x_1)\stackrel{a}{\rightarrow }\pi _1(X,x_2) 
\end{array}
\]
where $A([\Gamma ]):=[\Theta ^{-1}\cdot \Gamma \cdot \Theta ]$,
$a([\Sigma ]):=[(p\circ \Theta )^{-1}\cdot \Sigma \cdot (p\circ
\Theta )]$ and $\Theta $ is a definable path from $u_1$ to
$u_2$. Since $A$ and $a$ are isomorphisms and $p_*$ is a monomorphism
it follows from (1) that $|p^{-1}(x_1)|=|p^{-1}(x_2)|$. The same
diagram applied to the case $x_1=x_2=x$ shows that for all $u,v\in
p^{-1}(x)$, $p _*(\pi _1(Y,u))$ and $p _*(\pi _1(Y,v))$ are conjugate
subgroups of $\pi _1(X,x)$.

Suppose now that $S=[\Sigma ^{-1}]p_*(\pi _1(Y,u))[\Sigma ]$ for some 
$[\Sigma ]\in \pi _1(X,x)$. Let $\bar{\Sigma }$ be the definable
lifting of $\Sigma $ at $u$. Note that $v:=u*\Sigma \in
p^{-1}(x)$. Using the commutative diagram we see that $S=p_*(\pi
_1(Y,v)).$ 
\qed

\end{subsection}

\begin{subsection}{Liftings of strictly properly $\bigvee$-definable maps}
\label{subsection liftings}

The results of this subsection are all corollaries of following result
on the possibility of lifting strictly properly $\bigvee $-definable maps.

\begin{prop}\label{prop lifting functions}
Let $p:(Y,y)\into (X,x)$ be a strictly properly $\bigvee $-definable
covering map and let $f:(Z,z)\into (X,x)$ be a continuous strictly
properly $\bigvee $-definable map. Then there is a continuous
strictly properly $\bigvee $-definable map $\bar{f}:(Z,z)\into (Y,y)$
with $p\circ \bar{f}=f$ iff $f_*(\pi _1(Z,z))\subseteq p_*(\pi
_1(Y,y))$. Such strictly properly $\bigvee $-definable lifting
$\bar{f}$, when it exists, is unique. 
\end{prop}

\pf
The necessity is clear and the uniqueness follows from lemma \ref{lem
unique lift}. We will now construct $\bar{f}$. For each $i\in I$
choose $z_i\in Z_i$ such that if $z\in Z_i$ then $z=z_i$, and let 
$\Gamma _i$ be a definable path in $Z$ from $z$ to $z_i$. Given $w\in
Z_i$, let $\Delta _{z,w}$ be the definable path $\Gamma _i\cdot 
\Gamma _{i,z_i, w}$ from $z$ to $w$ where, $\Gamma _{i,z_i,w}$ is
given by proposition \ref{prop uniform path connected}. Let $\Sigma
_{z,w}:=f\circ \Delta _{z,w} $ and put $\bar{f}(w):=y*\Sigma _{z,w}$. 

If $w\in Z_i\cap Z_j$ then we have another definable path $\Delta
'_{z,w}$ from $z$ to $w$ obtained from $Z_j$. $f\circ (\Delta
'_{z,w}\cdot \Delta ^{-1}_{z,w})=\Sigma '_{z,w}\cdot \Sigma
^{-1}_{z,w}$ is a definable loop at $x$. By hypothesis, $[\Sigma
'_{z,w}\cdot \Sigma ^{-1}_{z,w}]\in p_*(\pi _1(Y,y))$ and by corollary
\ref{cor lifts with the same end point} (2), $y*\Sigma _{z,w}=y*\Sigma
'_{z,w}$ and so $\bar{f}$ is well defined. Note that the same argument
shows that $\bar{f}$ does not depend on the choice of the points
$z_i\in Z_i$ or of the definable paths $\Gamma _i$. 

Its easy to see that $\bar{f}$ is a strictly properly $\bigvee
$-definable map.  Let $\{U_l:l\in L\}$ be the $p$-admissible family of
definable open neighbourhoods in $X$. Since $f$ is continuous and
strictly properly $\bigvee $-definable, $\{f^{-1}(U_l):l\in L\}$ is an
open cover of $Z$ such that, for each $Z_i$ there is a finite subset
$S_i$ of $L$ such that $Z_i\subseteq \bigcup _{s\in
S_i}f^{-1}(U_s)$. For each $i$ and each $s\in S_i$ let
$V_{is}:=Z_i\cap f^{-1}(U_s)$ and let $v_{is}\in V_{is}$. Then
$V_{is}$ are open definable sets and $Z_i=\bigcup _{s\in
S_i}V_{is}$. Moreover, by the argument above we can assume that the
properly $\bigvee $-definable system $\Delta _{z,w}$ of definable
paths is, for $w\in V_{is}$ of the form $\Gamma _i\cdot \Gamma
_{i,z_i,v_{is}}\cdot \Gamma _{i,v_{is},w}$. Clearly, we get a properly 
$\bigvee $-definable system $\Lambda _{x,f(w)}:=f\circ \Delta _{z,w}$
of definable paths in $f(Z)$ which is of the form $(f\circ \Gamma
_i)\cdot (f\circ \Gamma _{i,z_i,v_{is}})\cdot (f\circ \Gamma
_{i,v_{is},w})$. From this we see clearly that $\bar{f}_{|V_{is}}$ is
definable and the result follows. 
                                 
We will now show that $\bar{f}$ is continuous at $w$: let $U$ be a
$p$-admissible definable neighbourhood of $f(w)$ and let $V$ be an
open definable set of $p^{-1}(U)$ that is mapped homeomorphically onto
$U$ by $p$ and that contains $\bar{f}(w)$. Choose a definably path
connected definable neighbourhood  $W$ of $w$ so that $f(W)\subseteq
U$. We need to show that $\bar{f}(W)\subseteq V$. For each $w'\in W$
there is a definable path $\Lambda $ from $w$ to $w'$ in $W$, and then
we can use $\Gamma \cdot \Lambda $ as the definable path from $z$ to
$w'$. The definable lifting of $f\circ (\Gamma \cdot \Lambda )=(f\circ
\Gamma )\cdot (f\circ \Lambda )$ is obtained by first lifting $f\circ
\Gamma $ and then lifting $f\circ \Lambda $. Since the latter lifting
stays in $V$, this shows that $\bar{f}(W)\subseteq V$.
\qed

\begin{cor}\label{cor lifting functions}
(1) Let $p:Y\into X$ and $q:Z\into X$ be strictly properly 
$\bigvee $-definable covering maps. There is a strictly 
properly $\bigvee $-definable
isomorphism between the strictly properly $\bigvee $-definable covering spaces
iff $p _*(\pi _1(Y,y))$ and $q _*(\pi _1(Z,z))$ are conjugate
subgroups of $\pi _1(X,x)$ iff $p^{-1}(x)$ and $q^{-1}(x)$ are
isomorphic $\pi _1(X,x)$-sets.

(2) If $p _*(\pi _1(Y,y))$ $\subseteq q _*(\pi _1(Z,z))$, then there
is a unique strictly properly $\bigvee $-definable covering map $r:(Y,y)\into
(Z,z)$ such that $p\circ r=q$. 
\end{cor}

\pf
(1) Suppose that $\phi :(Y,y)\into (Z,z)$ is a strictly properly 
$\bigvee $-definable isomorphism. 
Then $p_*(\pi _1(Y,y))=q_*(\phi (Y),\phi (y))$ which is
conjugate to $q_*(\pi _1(Z,z))$ by proposition \ref{prop size of fibres} 
(2). On the other hand it is clear that $\phi
_{|p^{-1}(x)}:p^{-1}(x)\into q^{-1}(x)$ is a bijection, for $[\Sigma
]\in \pi_1(X,x)$, if $\bar{\Sigma }$ is the definable lifting of
$\Sigma $ at $y$ then $\phi \circ \bar{\Sigma }$ is the definable
lifting of $\Sigma $ at $z$. And so $\phi $ induces an isomorphism
$\bar{\phi }:p^{-1}(x)\into q^{-1}(x)$ of $\pi
_1(X,x)$-sets. Conversely, if $p^{-1}(x)$ and $q^{-1}(x)$ are
isomorphic $\pi _1(X,x)$-sets then by the theory of $G$-sets, the
stabilisers $p_*(\pi _1(Y,y))$ of $y$ and $q_*(\pi _1(Z,z))$ of $z$ in
$\pi _1(X,x)$ are conjugate. By proposition \ref{prop size of fibres}
(2), there is $y'\in p^{-1}(x)$ such that $q_*(\pi _1(Z,z))=p_*(\pi
_1(Y,y'))$. By proposition \ref{prop lifting functions}, there is a
strictly properly $\bigvee $-definable map $\phi :Y\into Z$ (the
strictly properly $\bigvee$-definable lifting of $p:(Y,y')\into
(X,x)$) such that $q\circ \phi=p$. By considering the strictly
properly $\bigvee $-definable lifting of $q:(Z,z)\into (X,x)$ we see
that $\phi $ is a homeomorphism and so a strictly properly $\bigvee $-definable
isomorphism of pointed strictly properly  $\bigvee $-definable pointed
coverings spaces $((Y,y'), p,  (X,x))$ and $((Z,z), q, (X,x))$. 

(2) By proposition \ref{prop lifting functions}, there is a strictly
properly $\bigvee $-definable map $h :(Y,y)\into (Z,z)$ (the strictly
properly $\bigvee $-definable lifting of $p:(Y,y)\into (X,x)$) such
that $q\circ h=p$. It remains to show that $h$ is a strictly properly $\bigvee
$-definable covering map. Let $\{U_l:l\in L\}$ be a $p$-admissible
family of open definably connected definable neighbourhoods in $X$
with $|L|<\aleph _1 $and let $\{V_k:k\in K\}$ be a $q$-admissible
family of open definably connected definable neighbourhoods in $X$
with $|K|<\aleph _1$. For each $l\in L$ there is a finite set
$S_l\subseteq K$ such that $U_l\cap V_s\neq \emptyset $ for all $s\in
S_l$ and $U_l=\bigcup _{s\in S_l}U_l\cap V_s$. Let $\{W_{lk}:l\in L,\, 
k\in S_l\}$ be the family given by $W_{lk}:=U_l\cap V_k$ for all $l\in
L$ and $k\in S_l$. Then $\{W_{lk}:l\in L,\,k\in S_l\}$ is
simultaneously a $p$-admissible and $q$-admissible family of open
definable neighbourhoods in $X$ with $|\{(l,k):l\in L,\, k\in
S_l\}|<\aleph _1$. Since by o-minimality, $W_{lk}$ as only finitely
many definably connected components, we can assume without loss of
generality that each $W_{lk}$ is definably connected.

Let $\{O_m:m\in M\}$ be the family $\{q^{-1}(W_{lk}):l\in L,\, k\in
S_k\}$ and $\{O'_n:n\in N\}$ be the family $\{p^{-1}(W_{lk}):l\in L,\,
k\in S_l\}$. All of $O_m$'s (resp., $O'_n$) are open definably
connected definable neighbourhoods in $Z$ (resp., in $Y$) with $|M\cup
N|<\aleph _1$. We claim that $\{O_m:m\in M\}$ is a $h$-admissible
family of open definably connected definable neighbourhoods in $Z$:
Given $O_m$, there are $l\in L$ and $k\in S_l$ such that
$q_{|O_m}:O_m\into W_{lk}$ is a definable homeomorphism. Now let
$\{O'_{n'}:n'\in N'\subseteq N\}$ be the subfamily of $\{O_n:n\in N\}$
given by $p^{-1}(W_{lk})$. Then clearly $h^{-1}(O_m)=\{O'_{n'}:n'\in
N'\subseteq N\}$ and $h_{|O'_{n'}}:O'_{n'}\into O_m$ is a definable
homeomorphism.
\qed

\begin{cor}\label{cor cov and aut}
Let $p:(Y,y)\into (X,x)$ be a strictly properly $\bigvee $-definable covering
map and consider $p^{-1}(x)$ as a $\pi _1(X,x)$-set. Then we have
canonical isomorphisms
$$Aut(p^{-1}(x))\simeq Cov(Y/X)\simeq N_{\pi }(p_*(\pi
_1(Y,y)))/p_*(\pi _1(Y,y)),$$
where $\pi $ denotes $\pi _1(X,x)$.
\end{cor}

\pf
The prove that $Aut(p^{-1}(x))\simeq Cov(Y/X)$ is contained in the
prove of corollary \ref{cor lifting functions} (1). The rest follows
from the theory of $G$-sets (see lemma 10.26 \cite{ro}).
\qed

\begin{defn}\label{defn regular}
{\em
A strictly properly $\bigvee $-definable covering map $p:(Y,y)\into (X,x)$ is
{\it regular} if $p_*(\pi _1(Y,y))$ is a normal subgroup of $\pi
_1(X,x)$ . A strictly properly $\bigvee $-definable covering space $(Y,p,X)$ is
a {\it universal strictly properly $\bigvee $-definable covering space} of $X$
if $Y$ is definably simply connected.
}
\end{defn}

\begin{lem}\label{lem regular}
A strictly properly $\bigvee $-definable covering map $p:(Y,y)\into
(X,x)$ is {\it regular} iff $Cov(Y/X)$ acts transitively on $p^{-1}(x)$.
\end{lem}

\pf
Suppose that $p:(Y,y)\into (X,x)$ is regular and let $u,v\in
p^{-1}(x)$ Then by proposition \ref{prop size of fibres} (2), $p_*(\pi
_1(Y,u))=p_*(\pi _1(Y,v))$. By proposition \ref{prop lifting
functions}, there is a strictly properly $\bigvee $-definable
homeomorphism $h:(Y,u)\into (Y,v)$ such that $p\circ h=p$; thus $h\in
Cov(Y/X)$ and $h(u)=v$. Conversely, assume that $Cov(Y/X)$ acts
transitively on $p^{-1}(x)$ and let $u,v\in p^{-1}(x)$. Then there is
$h\in Cov(Y/X)$ such that $h(u)=v$. But, $p_*(\pi _1(Y,u))=p_*h_*(\pi
_1(Y,u))=p_*(\pi _1(Y,v)).$ By proposition \ref{prop size of fibres}
(2), $p_*(\pi _1(Y,u))$ is a normal subgroup of $\pi _1(X,x)$.
\qed

\begin{cor}\label{cor regular}
Suppose that $(Y,p,X)$ is a strictly properly $\bigvee $-definable covering
space. Then  $p:(Y,y)\into (X,x)$ is regular iff 
$$Cov(Y/X)\simeq  \pi _1(X,x)/p_*(\pi _1(Y,y))$$ 
and $(Y,p,X)$ is a strictly universal properly $\bigvee $-definable covering
space of $X$ iff 
$$Cov(Y/X)\simeq \pi _1(X,x).$$
\end{cor}

\end{subsection}

\begin{subsection}{Existence of $\bigvee $-definable covering spaces}
\label{susection existence of coverings}
 
Throughout this subsection, $G$ will be a subgroup of $\pi _1(X,x)$ with
$|G|<\aleph _1$ (it will follow from the main theorem of this
subsection that $|\pi _1(X,x)|<\aleph _1$).

\begin{lem}\label{lem nice n space} 
There is a cover $\{U_s:s\in S\}$ of $X$ with $|S|<\aleph _1$ by 
definable open subsets $U_s$'s such that 
every definable loop at $x\in U_s$ is definably
homotopic to the constant path $\epsilon _x$ at $x$ and for each $i\in
I$, there is a finite subset $S_i$ of $S$ with $X_i\subseteq \cup
\{U_s:s\in S_i\}$. We call such
family a good family of open definable neighbourhoods.
\end{lem}

\pf
If $X$ has a properly $\bigvee $-definable cell decomposition $K$, then 
for each $0$-cell $C_s$ ($s\in S$) of
$K$ let $U_s$ be the union of $C_s$ together with all open $k$-cells
$C$ of $K$ such that $C_s\subseteq \bar{C}$. Clearly, $U_s$ is an open
definably path connected definable subset,
which by remark \ref{nrmk paths in a closed cell} has the property of
the lemma. In general, since each $X_i$ has
a properly $\bigvee $-definable cell decomposition, by the argument
above, there is a finite such cover of each $X_i$ and from this, the
result follows.   
\qed

Let $X_G$ be the set of equivalence classes $[\Sigma ]_G$ of definable
paths $\Sigma $ in $X$ with initial point $\inf \Sigma =x$ under the 
following equivalence relation: two such definable paths $\Sigma $ and 
$\Lambda $ are equivalent iff $\sup \Sigma =\sup \Lambda $ and
$[\Sigma \cdot \Lambda ^{-1}]\in G$. Let $x_G:=[\epsilon _x]_G$ and
define $p_G:X_G\into X$ by $p_G([\Sigma ]_G)=\sup \Sigma $.

For each $s\in S$ let $u_s\in U_s$ and $s_G:=\{[\Sigma ]_G: \sup
\Sigma =u_s\}$. Consider the family $\{V_{sk}:s\in S,\, k\in s_G\}$
where each $V_{sk}$ is the set of all $[\Lambda ]_G$ such that there
are $[\Sigma ]_G\in s_G$ and a definable path $\Gamma $ in $U_s$ with
initial point $\inf \Gamma =u_s$ and $\Lambda =\Sigma \cdot \Gamma $. 

\begin{lem}\label{lem existence}
For all $s\in S$ and $k\in s_G$, $|s_G|$ and $V_{sk}$ are independent
of the choice of $u_s\in U_s$. Given $s'\in S$ then $|s'_G|=|s_G|$ and
$U_s\cap U_{s'}$ is non empty iff there are $k\in s_G$ and $l\in s'_G$
such that $V_{sk}\cap V_{s'l}$ is non empty.
If $k_1,k_2\in s_G$ are such that $k_1\neq k_2$ then $V_{sk_1}\cap
V_{sk_2}=\emptyset $.
\end{lem}

\pf
Let $u'_s\in U_s$, $k'=[\Sigma ']_G$$=[\Sigma \cdot \Theta ]_G$ where
$\inf \Theta =\sup \Sigma $ and $|\Theta |\subseteq U_s$, $V'_{sk'}$
the corresponding $V_{sk}$, $[\Lambda ]_G\in V_{sk}$,  and $k=[\Sigma
]_G$. Then $[\Lambda ]_G=[\Sigma \cdot \Gamma ]_G$ where $\inf \Gamma
=\sup \Sigma $ and $|\Gamma |\subseteq U_s$.  Hence, $[\Lambda
]_G=[\Sigma \cdot \Gamma ]_G$$=[(\Sigma \cdot \Theta )\cdot (\Theta
^{-1}\cdot \Gamma )]_G$$=[\Sigma '\cdot (\Theta ^{-1}\cdot \Gamma
)]_G$; since $\inf (\Theta ^{-1}\cdot \Gamma )$$=\sup \Sigma '$ and
$|(\Theta ^{-1}\cdot \Gamma )|\subseteq U_s$ we get $[\Lambda ]_G\in
V'_{sk'}$. The reverse inclusion is similar.   
 
If $U_s\cap U_{s'}$ is non empty then we have $|s_G|=|s'_G|$ since we
can take $u_s=u_{s'}\in U_s\cap U_{s'}$. We also have for the same
reason that there are $k\in s_G$ and $l\in s'_G$ such that $V_{sk}\cap
V_{s'l}$ is non empty. The general case follows because $X$ is
definably path connected. If there are $k\in s_G$ and $l\in s'_G$ such
that $[\Gamma ]_G\in V_{sk}\cap V_{s'l}$, then $p_G([\Gamma ]_G)\in
U_s\cap U_{s'}$.

Suppose that $[\Gamma ]_G\in V_{sk_1}\cap V_{sk_2}$. Let $k_i=[\Sigma
_i]_G$ and $[\Gamma ]_G=[\Sigma _i\cdot \Gamma _i]_G$ where $\inf
\Gamma _i=\sup \Sigma _i$ and $|\Gamma _i|\subseteq U_s$ for
$i=1,2$. Since, $\sup \Gamma =p_G([\Gamma ]_G)\in p_G(V_{sk})\cap
p_G(V_{sl})$ by the argument above, we can assume that $\sup \Sigma
_1=\sup \Sigma _2$. We have $[\Sigma _1\cdot \Sigma _2^{-1}]$$=[\Sigma
_1 \cdot (\Gamma _1\cdot \Gamma _2^{-1})\cdot \Sigma _2^{-1}]$
(since in $U_s$ every definable loop is definably homotopic to a
constant path) $=[(\Sigma _1\cdot \Gamma _1)\cdot (\Sigma _2\cdot \Gamma
_2)^{-1}]\in G$ and so $[\Sigma _1]_G=[\Sigma _2]_G.$
\qed

For each $s\in S$ let $I_s:=\{i\in I:U_s\cap X_i\neq \emptyset \}$ and
let $\{X_{Gskj}:s\in S,\, k\in s_G,\, j\in I_s\}$ be given by
$X_{Gskj}:=\{[\Lambda ]_G\in V_{sk}:\sup \Lambda \in U_s\cap
X_j\}$. Let $\phi _{Gskj}:X_{Gskj}\into \phi _j(X_j)$ be given by
$\phi _{Gskj}:=\phi _j\circ p_G$, and let ${\mathbf X}_G:=(X_G,
X_{Gskj},$$\phi _{Gskj})_{s\in S, k\in s_G, j\in I_s}$

\begin{thm}\label{thm existence}
For every subgroup $G$ of $\pi _1(X,x)$ with $|G|<\aleph _1$, ${\mathbf X}_G$ 
is a properly $\bigvee $-definably connected, properly $\bigvee $-definable
manifold. Moreover, $p_G:(X_G,x_G)\into (X,x)$ is a strictly properly $\bigvee
$-definable covering map such that $G=p_{G*}(\pi _1(X_G,x_G))$.
\end{thm}

\pf
$p_{G_{|V_{sk}}}:V_{sk}\into U_s$ is a bijection, because $U_s$ is
definably path connected and every definable loop in $U_s$ at $u_s$ is
definably homotopic to the constant path $\epsilon _{u_s}$. This
shows that ${\mathbf X}_G$ is a properly $\bigvee $-definable
manifold. Moreover, for each $s\in S$ and $j\in I_s$,
$p_G^{-1}(U_s\cap X_j)$ is the disjoint union of the definable sets
$X_{Gskj}$ with $k\in s_G$ and $p_{G_{|X_{Gskj}}}:X_{Gskj}\into U_s\cap
X_j$ is a bijection. Therefore, $p_G$ is an open, continuous and
surjective strictly properly $\bigvee $-definable map.

We now show that $X_G$ is definably path connected (and
therefore, by the above, $(X_G,p_G,X)$ is a strictly properly $\bigvee
$-definable covering space). Let $u=[\Sigma ]_G\in X_G$. We have
$\Sigma \simeq \Sigma _1\cdots \Sigma _n$ and there are for each $l\in
\{1, \dots , n\}$ there is $s(l)\in S$ such that $|\Sigma _l|\subseteq
U_{s(l)}$. Considering the definable bijection
$p_{G_{|V_{s(l)k(l)}}}:V_{s(l)k(l)}\into U_{s(l)}$ where $k(l)$ is such
that $x_G\in V_{s(0)k(0)}$, and for each $j=0,\dots , l-1$
$V_{s(j)k(j)}\cap V_{s(j+1)k(j+1)}\neq \emptyset $, it follows that
there is a definable path $\bar{\Sigma }\simeq \bar{\Sigma _1}\cdots
\bar{\Sigma _n}$ from $x_G$ to $u$ such that $p_G\circ \bar{\Sigma }=\Sigma $.

Finally, let us show that $G=p_{G*}(\pi _1(X_G,x_G)).$ Let $[\Sigma
]\in \pi _1(X,x)$. Then there is a unique definable lifting
$\bar{\Sigma }$ of $\Sigma $ in $X_G$ with $\inf \bar{\Sigma
}=x_G$. On the other hand, $[\Sigma ]\in p_{G*}(\pi _1(X_G,x_G))$ iff
$\inf \bar{\Sigma } =\sup \bar{\Sigma }=x_G$ iff $[\Sigma
]_G=[\epsilon _x]_G$ iff $[\Sigma ]=[\Sigma \cdot \epsilon _x^{-1}]\in
G$. Therefore, $G=p_{G*}(\pi _1(X_G,x_G)).$
\qed

\begin{cor}\label{cor universal covers}
Every strictly properly $\bigvee $-definable covering space
$((Y,y),p,$ $(X,x))$
is strictly properly $\bigvee $-definably isomorphic to a strictly properly
$\bigvee $-definable covering space of the form
$((X_G,x_G),p_G,(X,x))$. 
In particular, there is a universal strictly properly $\bigvee $-definable 
covering space $(\tilde{X},p,X)$.
\end{cor}

\end{subsection}

\begin{subsection}{Properly $\bigvee $-definable $G$-coverings spaces}
\label{subsection g coverings}

\begin{defn}\label{defn action}
{\em
Let $G$ be a group with $|G|<\aleph _1$. An action of $G$ on $X$ is a 
homomorphism $\rho :G\into Aut(X)$ induced by a 
map $\rho :G\times X\into X$, where $Aut(X)$ is
the group of all strictly properly $\bigvee $-definable homeomorphisms
of $X$. We often use the notation $g(x)=\rho (g,x)$.
The {\it orbit} $x/G$ of $x\in X$ is the subset $\{g(x):g\in G\}$ of
$X$. $X$ is the disjoint
union of the orbits. We denote by $X/G$ the set of orbits, and
$r:X\into X/G$ denotes the map that sends a point into its orbit (and
so, $r^{-1}(x/G)=x/G$). $X/G$ has a topology such that $r:X\into X/G$
is a continuous and open surjective map, but in general its not 
clear that $X/G$ can be made into a properly $\bigvee $-definable
manifold such that $r:X\into X/G$ is a strictly (properly) 
$\bigvee $-definable map.
}
\end{defn}

\begin{defn}\label{defn equivalent}
{\em
Two strictly properly $\bigvee $-definable covering spaces $(\bar{Y},p,Y)$ and
$(\bar{X},q,X)$ are {\it strictly properly $\bigvee $-definably
equivalent} if there are strictly properly $\bigvee $-definable
homeomorphisms $\phi $ and $\psi $ making the following diagram commutative:
\[ 
\begin{array}{clrc}
\bar{Y}\,\,\,\,\stackrel{\phi }{\rightarrow }\,\,\,\,\bar{X}\\
\downarrow ^p\,\,\,\,\,\,\,\,\,\,\,\,\,\,\,\downarrow ^q\\
Y\,\,\,\,\stackrel{\psi }{\rightarrow }\,\,\,\,X.
\end{array}
\] 
}
\end{defn}

\begin{lem}\label{lem actions and covers}
Consider the following commutative diagram of strictly properly 
$\bigvee $-definable covering maps
\[
\begin{array}{clrc}
Z\\
\,\,\,\,\,\,\,\,\,\,\,\,\,\,\,\,\,\,\searrow ^r\\
\,\,\,\,\,\,\,\,\,\,\,\,\,\,\,\,\,\,\downarrow ^p\,\,\,\,\,\,\,\,\,Y\\
\,\,\,\,\,\,\,\,\,\,\,\,\,\,\,\,\,\,\swarrow ^q\\
X
\end{array}
\]
where $(Z,p,X)$ and $(Z,r,Y)$ are regular; let $G=Cov(Z/Y)$ and
$H=Cov(Z/X)$. Then there is a commutative diagram 
\[
\begin{array}{clrc}
Z\\
\,\,\,\,\,\,\,\,\,\,\,\,\,\,\,\,\,\,\,\,\,\,\,\,\,\,\,\searrow ^{r'}\\
\,\,\,\,\,\,\,\,\,\,\,\,\,\,\,\,\,\,\,\,\,\,\,\,\,\,\,\downarrow ^{p'}
\,\,\,\,\,\,\,\,\,Z/G\\
\,\,\,\,\,\,\,\,\,\,\,\,\,\,\,\,\,\,\,\,\,\,\,\,\,\,\,\swarrow ^{q'}\\
Z/H
\end{array}
\]
of strictly properly $\bigvee $-definable covering spaces, each of which is
strictly properly $\bigvee $-definably equivalent to the corresponding
strictly properly $\bigvee $-definable covering space in the original diagram.
\end{lem}

\pf
We first show that $(Z,r',Z/G)$ is a strictly properly $\bigvee $-definable
covering space (here $r'$ is the natural map $r':Z\into Z/G$), strictly
properly $\bigvee $-definably equivalent to $(Z,r,Y)$. Note that since
$(Z,r,Y)$ is regular, $Cov(Z/Y)$ acts transitively on each fibre
$r^{-1}(y)$. This can be used to show that the map $\phi :Y\into Z/G$
given by $\phi (y)=r'(r^{-1}(y))$ is a well defined bijective map. 
Since both $r'$ and $r$ are continuous and open surjective maps, 
it follows that $\phi $ is a homeomorphism, from this we get the claim,
by putting on $Z/G$ the structure of a properly $\bigvee $-definable manifold
via $\phi $ which then becomes a well defined bijective strictly
properly $\bigvee $-definable map.

Since $G\subseteq H$ if follows that for every $z\in Z$, we have
$z/G\subseteq z/H$. Define $q':Z/G\into Z/H$ by $q'(z/G)=z/H.$ It is
now easy to verify that $q'$ is the map we are after.
\qed

\begin{defn}\label{defn G covering}
{\em
We say that a strictly properly $\bigvee $-definable covering space $(Z,r,Y)$
is a {\it strictly properly $\bigvee $-definable $G$-covering} if 
$(Z,r',Z/G)$ is a strictly properly $\bigvee $-definable covering
space strictly properly $\bigvee $-definably equivalent to $(Z,r,Y)$. 
}
\end{defn}

\begin{cor}\label{cor action and covers}
Let $(\tilde{X},p,X)$ be a strictly properly $\bigvee $-definable universal
covering space of $X$. Then every strictly properly  $\bigvee $-definable
covering space $(Y,q,X)$ is strictly properly $\bigvee $-definably
equivalent to $(\tilde{X}/G,r,X)$ for some subgroup $G$ of
$Cov(\tilde{X}/X)\simeq \pi _1(X,x).$
\end{cor}

\pf
There is a unique strictly properly $\bigvee $-definable covering space
$(\tilde{X}, s,Y)$ such that $p=q\circ s$. Since $\tilde{X}$ is
definably simply connected, both $(\tilde{X},p,X)$ and
$(\tilde{X},s,Y)$ are regular strictly properly $\bigvee $-definable
covering spaces. Therefore, by lemma \ref{lem actions and covers}, $(Y,q,X)$ is
strictly properly $\bigvee $-definably equivalent to
$(\tilde{X}/G,r,X)$ where $G=Cov(\tilde{X}/Y)$.
\qed

\begin{cor}\label{cor cov and covers}
Let $(\tilde{X},p,X)$ be a strictly properly $\bigvee $-definable universal
covering space of $X$. Denote the family of all strictly properly $\bigvee
$-definable covering spaces of $X$ of the form $(\tilde{X}/G,r,X)$,
where $G$ is a subgroup of $Cov(\tilde{X}/X)$, by $\mathcal{C}$ and
denote the family of all subgroups of $Cov(\tilde{X}/X)$ by
$\mathcal{G}$. Then $\Phi :\mathcal{C}\into \mathcal{G}$ defined by
$(Z,q,X)\into Cov(\tilde{X}/Z)$ and $\Psi :\mathcal{G}\into
\mathcal{C}$ defined by $G\into (\tilde{X}/G,r,X)$ are bijections
inverse to one another.
\end{cor}

\pf 
If $G\subseteq Cov(\tilde{X}/X)$, then $\Phi \Psi
(G)=Cov(\tilde{X}/(\tilde{X}/G))$. This group consists of all strictly
properly $\bigvee $-definable homeomorphisms $h:\tilde{X}\into
\tilde{X}$ such that 
\[
\begin{array}{clrc}
\tilde{X}\\
\,\,\,\,\,\,\,\,\,\,\,\,\,\,\,\,\,\,\,\,\,\,\,\,\,\,\,\searrow ^{h}\\
\,\,\,\,\,\,\,\,\,\,\,\,\,\,\,\,\,\,\,\,\downarrow ^{r}
\,\,\,\,\,\,\,\,\,\,\,\,\,\tilde{X}\\
\,\,\,\,\,\,\,\,\,\,\,\,\,\,\,\,\,\,\,\,\,\,\,\,\swarrow ^{r}\\
\tilde{X}/G
\end{array}
\] 

If $g\in G$ and $\tilde{x}\in \tilde{X}$, then
$\tilde{x}/G$$=g(\tilde{x})/G$,and so $rg=r$; hence $g\in \Phi \Psi
(G)$ and $G\subseteq \Phi \Psi (G)$. For the reverse inclusion, let
$h\in \Phi \Psi (G)$, then $rh=h$. If $\tilde{x}\in \tilde{X}$, then
$\tilde{x}/G=h(\tilde{x})/G$ and there exists $g\in G$ such that
$g(h(\tilde{x}))=\tilde{x}$. Since $g\in \Phi \Psi (G)$, it follows
that $gh\in \Phi \Psi (G)$ and by uniqueness, $gh=1_{\tilde{X}}$ and
$h=g^{-1}\in G$.

Similarly, its also easy to see that $\Psi \Phi $ is the identity.
\qed

\begin{cor}\label{cor fundamental group and covers}
Let $(\tilde{X},p,X)$ be a strictly properly 
$\bigvee$-definable universal covering
space of $X$. If $G$ is a subgroup of $Cov(\tilde{X}/X)\simeq \pi _1(X,x)$,
then $\pi _1(\tilde{X}/G,\tilde{x}/G)\simeq G$.
\end{cor}

\pf
We have $\pi _1(\tilde{X}/G,\tilde{x}/G)\simeq
Cov(\tilde{X}/(\tilde{X}/G))\simeq G$.  
\qed

\begin{defn}\label{defn proper action}
{\em
We say that $G$ {\it acts properly on $X$} if there is a cover 
$\{U_k:k\in K\}$ of $X$ with $|K|<\aleph _1$ such that each $U_k$ is 
an open definably connected definable subset and $gU_k\cap
U_k=\emptyset $  for all $g\in G\setminus \{1\}$. We call the 
$\{U_k:k\in K\}$ a {\it $G$-admissible family of definable open
subsets of $X$}.
}
\end{defn}

Its easy to see that if $(Z,r,Z/G)$ is a strictly properly 
$\bigvee $-definable $G$-covering space such that $G$ acts on $X$
without fixed points (e.g., $(Z,r,Z/G)$ is regular) then $G$ acts
properly on $Z$. 

\begin{cor}\label{cor proper action}
Suppose that $G$ acts properly on $X$. Then $(X,r,X/G)$ is a regular 
strictly properly $\bigvee $-definable covering space and 
$$G\simeq Cov(X/(X/G))\simeq \pi _1(X/G,x/G)/r_*(\pi _1(X,x)).$$
\end{cor}

\pf
The natural map $r:X\into X/G$ is a continuous, surjective open map
with $r^{-1}(r(U))=\bigcup _{g\in G}gU$ for any definable open subset 
$U$ of $X$. Let $\{U_k:k\in K\}$ be a $G$-admissible family of
definable open subsets of $X$. We claim that there is a structure of a 
$\bigvee $-definable manifold on $X/G$ such that $\{U_k:k\in K\}$ is
an $r$-admissible family definable open subsets of $X/G$. Let $k\in
K$, if $g,h$ are distinct elements of $G$ then $gU_k\cap hU_k=
\emptyset $. If $u\in U_k$, then $u/G=(gu)/G$ for every $g\in G$, and
so $r(gu)=r(u)$; hence $r_{|gU_k}$ is surjective. If $r(gu)=r(gv)$ for 
$u,v\in U_k$, then there is $h\in G$ with $gu=hgv$; hence 
$gU_k\cap hgU_k\neq \emptyset $, a contradiction. Therefore,
$r_{|gU_k}$ is a bijection. Let $Y:=X/G$, for each $i\in I$ and $k\in
K$, let $Y_{i,k}:=r_{|U_k\cap X_i}(U_k\cap X_i)$ and let 
$\psi _{i,k}:Y_{i,k}\into \phi _i(X_i)$ be given by 
$\psi _{i,k}(y):=\phi _i((r_{|U_k})^{-1}(y))$. Its now easy to see
that $(Y,(Y_{i,k},\psi _{i,k})_{i\in I,\,\, k\in K})$ is a (properly) 
$\bigvee $-definable manifold and $(X,r,Y)$ is a strictly properly 
$\bigvee $-definable $G$-covering space with $r$-admissible family of 
definable open subsets given by $\{Y_{i,k}:i\in I,\,\, k\in K\}$.

Its easy to see that $Cov(X/(X/G))$ acts transitively on a fiber, 
therefore $(X,r,X/G)$ is regular. The rest follows from previous results.   
\qed  

We denote by $G$-$COV(X)$ (resp., $G$-$COV(X,x)$) the set of equivalence
classes of strictly properly $\bigvee $-definable $G$-coverings
$(Y,p,X)$ of $X$ (resp., pointed $G$-coverings $((Y,y),p,(X,x))$ of
$(X,x)$) under strictly properly 
$\bigvee $-definable isomorphisms of strictly properly $\bigvee
$-definable $G$-coverings (resp., pointed $G$-coverings) of $X$ 
(resp., of $(X,x)$).

\begin{prop}\label{prop hom and g coverings}
There is a canonical bijection between $G$-$COV(X,x)$ and $Hom(\pi _1(X,x),G)$.
\end{prop}

\pf
We first define a map $A:Hom(\pi _1(X,x),G)\into G$-$COV(X,x)$ has 
follows: Let $\rho \in Hom(\pi _1(X,x),G)$ and let 
$((\tilde{X},\tilde{x}),p,(X,x))$ be the universal strictly properly 
$\bigvee $-definable covering space of $(X,x)$. Consider the action of 
$\pi _1(X,x)$ on $\tilde{X}\times G$ given by  $[\Sigma ](z,g)
:=(z*\Sigma ,g\rho ([\Sigma ]^{-1}))$. We claim that $\pi _1(X,x)$ 
acts properly on $\tilde{X}\times G$. In fact, if $\{U_l:l\in L\}$ is 
a $p$-admissible family of definable open neighbourhoods then for each 
$l\in L$, there is a strictly properly $\bigvee $-definable
homeomorphism $M_l:p^{-1}(U_l)\into $ $U_l\times \pi _1(X,x)$. Let 
$\{V_{l,m}:l\in L,\,\, m\in \pi _1(X,x)\}$$=\{V\subseteq p^{-1}(U_l):
l\in L,\,\, M_l(V)=U_l\times \{m\}\}$ and for $l\in L$, $g\in G$ and 
$m\in \pi _1(X,x)$, let $W_{l,m,g}:=V_{l,m}\times \{g\}$. Then, 
$\{W_{l,m,g}:l\in L,\,\, m\in \pi _1(X,x),\,\, g\in G\}$ is a cover of 
$\tilde{X}\times G$ by open definably connected, definable subsets
such that $[\Sigma ]W_{l,m,g}\cap W_{l,m,g}=\emptyset $ for all 
$[\Sigma ]\in \pi _1(X,x)\setminus \{1\}$. Let $(Y,y)
:=(\tilde{X}\times G/\pi _1(X,x), (\tilde{x},1)/\pi _1(X,x))$ and let 
$r:Y\into X$ be the strictly properly $\bigvee $-definable map induced
by $p$. There is a natural action of $G$ on $Y$ induced by the natural 
action of $G$ on $\tilde{X}\times G$. Its now easy to see, using an 
argument similar to the one above, that $G$ acts properly on $Y$ and 
$((Y,y),r,(X,x))$ is a strictly properly $\bigvee $-definable 
$G$-covering space.

Let $B:G$-$COV(X,x)\into Hom(\pi _1(X,x),G)$ be the map defined as
follows: given a strictly properly $\bigvee $-definable $G$ covering
space $((Y,y),p,(X,x))$, let $\rho \in Hom(\pi _1(X,x),G)$ be
determined by $\rho ([\Sigma ])y=y*\Sigma $. Now, some easy
computations, show that the map $A$ is the inverse of the map $B$. 
\qed

\end{subsection}

\begin{subsection}{The Seifert-van Kampen theorem}
\label{subsection seifert van kampem theorem}

We give here the proof of the o-minimal version of the Seifert-van
Kampen theorem for $X$ without assuming that $X$ is properly 
$\bigvee $-definably complete. The prove we present here is analogue to
Grothendieck proof in the classical case (see \cite{ro}).
 
\begin{lem}\label{lem patching covers}
Let $\{X^{\alpha }:\alpha \in A\}$ with $|A|<\aleph _1$ be a cover of
$X$ by properly $\bigvee $-definable open sets with properly 
$\bigvee $-definable covering maps $p_{\alpha }:Y^{\alpha }\into
X^{\alpha }$. Suppose that for each $\alpha ,\beta \in A$, we have strictly
properly $\bigvee $-definable isomorphisms
$$ \theta _{\alpha ,\beta }:p_{\alpha }^{-1}(X^{\alpha }\cap X^{\beta
})\into p_{\beta }^{-1}(X^{\alpha }\cap X^{\beta })$$
of strictly properly $\bigvee $-definable coverings of $X^{\alpha }\cap
X^{\beta }$, such that for every $\alpha ,\, \beta ,\, \gamma $$\in A$ we
have $\theta _{\alpha ,\alpha }=1_{Y^{\alpha }}$ and $\theta _{\gamma
,\alpha }=\theta _{\gamma ,\beta }\circ \theta _{\beta ,\alpha }$ on
$p_{\alpha }^{-1}(X^{\alpha }\cap X^{\beta }\cap X^{\gamma })$.
Then there is a strictly properly $\bigvee $-definable covering map $p:Y\into
X$ and there are strictly properly $\bigvee $-definable isomorphisms
$\lambda _{\alpha }:Y^{\alpha }\into p^{-1}(X^{\alpha })$ of strictly properly
$\bigvee $-definable coverings of $X^{\alpha }$ such that $\theta
_{\alpha ,\beta}=\lambda _{\beta }^{-1}\circ \lambda _{\alpha }$ on
$p_{\alpha }^{-1}(X^{\alpha }\cap X^{\beta})$ and $Y=\bigcup \{\lambda
_{\alpha }(Y^{\alpha }):\alpha \in A\}$. 

Moreover, if each $p_{\alpha }:Y^{\alpha }\into X^{\alpha }$ is a
strictly 
properly $\bigvee $-definable $G$-covering and each $\theta _{\alpha
,\beta }$ is a strictly properly $\bigvee $-definable isomorphism of
strictly properly $\bigvee $-definable $G$-coverings, then $p:Y\into X$ is a
strictly properly $\bigvee $-definable $G$-covering and each $\lambda _{\alpha
}:Y^{\alpha }\into p^{-1}(X^{\alpha })$ is a strictly properly
$\bigvee $-definable isomorphism of $\bigvee $-definable
$G$-coverings.
\end{lem}

\pf
We take $Y$ to be the disjoint union of copies $Y_{\alpha }$'s of the
$Y^{\alpha }$'s, we take on $Y$ the natural charts induced by those of
the $Y^{\alpha }$'s together with the charts induced by the $\theta
_{\alpha ,\beta }$'s. $\lambda ^{\alpha }:Y^{\alpha }\into Y_{\alpha
}$ is the natural strictly properly $\bigvee $-definable
homeomorphism, $p:Y\into X$ is defined by $p\circ \lambda _{\alpha
}=p_{\alpha }$ on $Y^{\alpha }$. Since $\lambda _{\alpha }$ is a
strictly properly $\bigvee $-definable homeomorphism, $p$ is in fact a
strictly properly $\bigvee $-definable covering map.

If each $p_{\alpha }:Y^{\alpha }\into X^{\alpha }$ is a strictly properly
$\bigvee $-definable covering map and each $\theta _{\alpha ,\beta }$
is a strictly properly $\bigvee $-definable isomorphism of strictly
properly $\bigvee $-definable $G$-coverings, then there is a unique
action of $G$ on $Y$
commuting with each $\lambda _{\alpha }$  i.e., $\lambda _{\alpha
}(gy^{\alpha })=g\lambda _{\alpha }(y^{\alpha })$ for $g\in G$ and
$y^{\alpha }\in Y^{\alpha }$. This gives the strictly properly 
$\bigvee $-definable covering the structure of a strictly 
properly $\bigvee $-definable $G$-covering, so that each 
$\lambda _{\alpha }$ is a strictly properly
$\bigvee $-definable isomorphism of strictly properly $\bigvee $-definable
$G$-coverings.
\qed

\begin{thm}\label{thm seifert-van kampen}
Let $\{X^{\alpha }:\alpha \in A\}$ with $|A|<\aleph _1$ be a cover of
$X$ by open properly $\bigvee $-definably connected properly $\bigvee
$-definable subsets such that for any $\alpha ,\beta \in A$ there is
$\gamma \in A$ with $X^{\gamma }=X^{\alpha }\cap X^{\beta }$ and for
all $\alpha \in A$, $x_0\in X^{\alpha }$ where $x_0$ is some point of
$X$. Then for any group $G$ and any homomorphisms $h_{\alpha }:\pi
_1(X^{\alpha },x_0)\into G$ such that $h_{\beta }=h_{\alpha }\circ
i_{\alpha \beta}$ whenever $X^{\beta }\subseteq X^{\alpha }$ and
$i_{\alpha \beta }:\pi _1(X^{\beta },x_0)\into \pi _1(X^{\alpha
},x_0)$ is the homomorphism induced by the inclusion, there is a
unique homomorphism $h:\pi _1(X,x_0)\into G$ such that $h_{\alpha
}=h\circ j_{\alpha }$ for all $\alpha \in A$, where $j_{\alpha }:\pi
_1(X^{\alpha },x_0)\into \pi _1(X,x_0)$ is the homomorphism induced by
the inclusion.
\end{thm}

\pf
By proposition \ref{prop hom and g coverings}, 
the homomorphisms $h_{\alpha }:\pi _1(X^{\alpha },x_0)\into G$
determine strictly properly $\bigvee $-definable $G$-coverings $p_{\alpha
}:Y^{\alpha }\into X^{\alpha }$ together with base points $y_0^{\alpha
}$ over $x_0$. The equality $h_{\beta }=h_{\alpha }\circ i_{\alpha
\beta}$ whenever $X^{\beta }\subseteq X^{\alpha }$ and $i_{\alpha
\beta }:\pi _1(X^{\beta },x_0)\into \pi _1(X^{\alpha },x_0)$ is the
homomorphism induced by the inclusion, allows by lemma \ref{lem
patching covers} the construction of a strictly properly $\bigvee $-definable
$G$-covering $p:Y\into X$ that restricts to the strictly properly 
$\bigvee $-definable $G$-coverings $p_{\alpha }:Y^{\alpha }\into X^{\alpha
}$. This strictly properly $\bigvee $-definable $G$-covering corresponds to a
homomorphism $h:\pi _1(X,x_0)\into G$, and the fact that the
restricted coverings agree means precisely that $h\circ j_{\alpha
}=h_{\alpha }$.
\qed

\begin{cor}\label{cor tietze2}
Suppose that $X$ is properly $\bigvee $-definably complete. Then
$\pi _1(X,x)$ is invariant under taking
elementary extensions, elementary substructures of $\N$ (containing
the parameters over which ${\mathbf X}$ is defined) and under taking
expansions of $\N$ and reducts of $\N$ on which ${\mathbf X}$ is
defined and $X$ has definable choice.  
\end{cor}

\pf
Since $X$ is properly $\bigvee $-definably complete, there is a cover 
$\{X^{\alpha }:\alpha \in A\}$ (with $|A|<\aleph _1$) of
$X$ by open properly $\bigvee $-definably connected, locally finite properly 
$\bigvee $-definable subsets such that for any $\alpha ,\beta \in A$ there is
$\gamma \in A$ with $X^{\gamma }=X^{\alpha }\cap X^{\beta }$ and for
all $\alpha \in A$, $x_0\in X^{\alpha }$ where $x_0$ is some point of
$X$. The result now follows from theorem \ref{thm tietze} and 
theorem \ref{thm seifert-van kampen}.
\qed

\end{subsection}

\begin{subsection}{Strictly properly $\bigvee $-definable 
$\mbox{\u C}$ech cohomology}
\label{subsection cech cohomology}

\begin{defn}\label{defn properly bigvee definable cover}
{\em
A {\it properly $\bigvee $-definable cover} of $X$ is a covering 
$\mathcal{U}$$=\{U_l:l\in L\}$ of $X$ with $|L|<\aleph _1$ such that
for each $l\in L$, $U_l$ is an open definably connected subset of $X$
and for each $i\in I$, there is a finite subset $L_i$ of $L$ such that 
$X_i\subseteq\cup \{U_l:l\in L_i\}$. We say that a properly 
$\bigvee $-definable cover ${\mathcal U}$ of $X$ is definably simply 
connected if for each $l\in L$, $U_l$ is definably simply connected.
}
\end{defn}

\begin{defn}\label{defn cech cohomology}
{\em
Let $G$ be a group with $|G|<\aleph _1$ and consider $G$ as a properly 
$\bigvee $-definable manifold of dimension zero. Let $\mathcal{U}$
$=\{U_l:l\in L\}$ be a properly $\bigvee $-definable cover of $X$. 
Let $S=\{(m,n)\in L\times L: U_m\cap U_n\neq \emptyset \}$. 
A {\it strictly properly $\bigvee $-definable \u{C}ech cochain of 
$\mathcal{U}$} if a collection $\{g_{m ,n}:(m,n)\in S\}$ (denoted 
$\{g_{m,n}\}$) of strictly properly $\bigvee $-definable continuous
maps $g_{m,n}:U_m\cap U_n\into G$. A strictly properly 
$\bigvee $-definable \u{C}ech cochain $\{g_{m,n}\}$ of $\mathcal{U}$ is 
a {\it strictly properly $\bigvee $-definable \u{C}ech cocycle of 
$\mathcal{U}$}, if for all $m,n,k\in L$ such that $(m,n), (m,k),
(k,n)\in S$ the following properties hold: $(1)$ $g_{m,m}=1_G$; $(2)$ 
$g_{m,n}=(g_{n,m})^{-1}$ and $(3)$ $g_{m,n}=g_{m,k}g_{k,n}$ on 
$U_m\cap U_k\cap U_m$.  Two strictly properly $\bigvee $-definable 
\u{C}ech cocycles $\{g_{m,n}\}$ and $\{h_{m,n}\}$ of $\mathcal{U}$ are 
said to be {\it strictly properly $\bigvee $-definably cohomologous} 
if there are strictly properly $\bigvee $-definable continuous maps 
$k_m:U_m\into G$ such that $h_{m,n}=(k_m)^{-1}g_{m,n}k_n$ on $U_m\cap
U_n$ for all $m,n\in L$ such that $(m,n)\in S$.  This is an
equivalence relation, the equivalence classes are called {\it strictly 
properly $\bigvee $-definable \u{C}ech cohomology classes on 
${\mathcal U}$ with coeficents on $G$} and the set of equivalence
classes is denoted by $H^1(X,{\mathcal U};G)$.  
}
\end{defn}
 
\begin{prop}\label{prop covering and cech cohomology}
Suppose that ${\mathcal U}$ is a definably simply connected, properly 
$\bigvee $-definable cover of $X$. Then there are canonical bijections 
between the following sets: $H^1(X,{\mathcal U};G)$, $G$-$COV(X)$ and 
$Hom(\pi _1(X,x),G)/conjugancy$. 
\end{prop}

\pf
A canonical bijection between $Hom(\pi _1(X,x),G)/conjugancy$ and 
$G$-$COV(X)$ is obtained from the proof of
proposition \ref{prop hom and g coverings}. We now construct a
bijection between $H^1(X,{\mathcal U};G)$ and $G$-$COV(X)$.

Let $(Y,p,X)$ be a strictly properly $\bigvee $-definable
$G$-covering. Since each $U_l$ is definably simply connected, there
are strictly properly $\bigvee $-definable isomorphisms 
$\alpha _l:U_l\times G\into p^{-1}(U_l)$ of strictly properly 
$\bigvee $-definable $G$-coverings. For each $(m,n)\in S$, let 
$\alpha _{m,n}:U_m\cap U_n\times G\into U_m\cap U_n\times G$ be given
by $\alpha _{m,n}:=(\alpha _{m})^{-1}\circ \alpha _n$. Its easy that, 
for each $(m,n)\in S$, there is a unique strictly properly 
$\bigvee $-definable continuous map $g_{m,n}:U_m\cap U_n\into G$ such
that $\alpha _{m,n}(x,g)=(x,gg_{m,n}(x))$.  It folows from the
definitions that $\{g_{m,n}\}$ is a strictly properly 
$\bigvee $-definable \u{C}ech cocycle on ${\mathcal U}$. Moreover, an
argument similar to the one above, shows that the class of
$\{g_{m,n}\}$ in $H^1(X,{\mathcal U};G)$ does not depend on the choice
of $\{\alpha _l:l\in L\}$ and depends only on  the class of $(Y,p,X)$
in $G$-$COV(X)$.

Conversely, given a strictly properly $\bigvee $-definable \u{C}ech
cocycle $\{g_{m,n}\}$ on ${\mathcal U}$ with coeficients in $G$, the 
result follows from lemma \ref{lem patching covers} 
if we take, for each $l,m,n\in L$ 
such that $(m,n)\in S$, $X^l=U_l$, $Y^l=U_l\times G$, $p_l:Y^l\into
X^l$ given by $p_l(x,g):=x$ and $\theta _{m,n}:p_m^{-1}(X^m\cap X^n)
\into p_m^{-1}(X^m\cap X^n)$ given by $\theta
_{m,n}(x,g):=(x,gg_{m,n}(x))$.  Moreover, its easy to verify that, the 
strictly properly $\bigvee $-definable $G$-coverings constructed in
this way from strictly properly $\bigvee $-definable cohomologous
strictly properly $\bigvee $-definable \u{C}ech cocycles are strictly
properly $\bigvee $-definably isomorphic $G$-coverings.
\qed

Let $H^1((X,x),{\mathcal U};G)$ be the set of equivalence classes of 
strictly properly $\bigvee $-definable \u{C}ech cocycles $\{g_{m,n}\}$ on 
${\mathcal U}$ with coefficients in $G$ such that $g_{m,n}(x)=e$ for
all $(m,n)\in S$ with $x\in U_m\cap U_n$, under the equivalence
relation given by $\{g_{m,n}\}\thicksim \{h_{m,n}\}$ iff there are
strictly properly $\bigvee $-definable continuous maps 
$\{k_l:U_l\into G:l\in L\}$ such that for all $l,m,n\in L$ with $x\in
U_l$ and $(m,n)\in S$, we have $k_l(x)=e$ and $h_{m,n}=(k_m)^{-1}g_{m,n}k_m$.

\begin{cor}\label{cor covering and cech cohomology}
Let ${\mathcal U}$ be a definably simply connected properly 
$\bigvee $-definable cover of $X$. Then there are canonical bijections 
between the following sets: $H^1((X,x),{\mathcal U};G)$, $G$-$COV(X,x)$
and $Hom(\pi _1(X,x),G)$. 
\end{cor}

\end{subsection}

\begin{subsection}{Strictly properly $\bigvee $-definable $H$-groups}
\label{subsection bigvee definable h groups}

\begin{defn}\label{defn H manifold}
{\em
$(X,\mu ,$$x_0)$ is a 
{\it strictly properly $\bigvee $-definable $H$-manifold} if 
the {\it strictly properly $\bigvee $-definable
multiplication} $\mu :(X\times X,(x_0,x_0)\into (X,x_0)$ and the {\it
unit} $x_0$ are continuous and satisfy $[\mu \circ i_1]\simeq [1_X]=[\mu
\circ i_2]$, where $i_1,i_2:X\into X\times X$ are the continuous
strictly properly $\bigvee $-definable maps $i_1(x)=(x,x_0)$ and
$i_2(x)=(x_0,x)$. A strictly properly $\bigvee $-definable $H$-group
is $(X,\mu ,\iota ,x_0)$ where $(X,\mu ,x_0)$ is a strictly properly
$\bigvee $-definable $H$-manifold with $[\mu \circ (\mu \times
1_X)]=[\mu \circ (1_X\times \mu)]$ (i.e., $\mu $ is strictly properly
$\bigvee $-definably $H$-associative) and the {\it strictly properly
$\bigvee $-definable $H$-inverse} $\iota :X\into X$ is continuous and
satisfies $[\mu \circ (\iota \times 1_x)\circ \Delta _X]=$$[\epsilon
_{x_0}]$$=[\mu \circ (1_x\times \iota )\circ \Delta _X]$ where $\Delta
_X:X\into X\times X$ is the diagonal map. A strictly properly $\bigvee
$-definable $H$-group $(X,\mu ,\iota ,x_0)$ is strictly properly
$\bigvee $-definably $H$-abelian if $[\mu ]=[\mu \circ \tau ]$ where
$\tau :X\times X\into X\times X$ is given by $\tau (x,y)=(y,x).$
}
\end{defn}

\begin{defn}\label{defn h homomorphisms}
{\em
A pointed, continuous strictly properly $\bigvee $-definable map
$h:(X,x_0)\into (Y,y_0)$ between strictly properly $\bigvee
$-definable $H$-manifolds $(X,\mu ,x_0)$ and $(Y,\gamma ,y_0)$ (resp.,
$H$-groups $(X,\mu ,\iota ,x_0)$ and $(Y,\gamma ,\zeta ,y_0)$) is
called a {\it strictly properly $\bigvee $-definable $H$-map} (resp.,
{\it $H$-homomorphism}) if $[h\circ \mu]=[\gamma \circ (h\times h)]$
(resp., also $[h\circ \iota ]=[\zeta \circ h]$).
}
\end{defn}

Let $\Sigma $ 
and $\Gamma $ be  definable path in $X$ with $\Sigma =\Sigma  _1\cdots 
\Sigma _k$, $\Gamma =\Gamma _1\cdots \Gamma _l$ where the
$\Sigma _i$'s (resp., the $\Gamma _j$'s) are definable basic paths
parametrised by $\sigma _i :(a_i,b_i)\into X$ (resp., $\gamma
_j:(c_j,d_j)\into X$). For $t\in [a_i,b_i]$ (resp., $t\in [c_j,d_j]$)
we define $(\Sigma ,\Gamma )(t)$ $:=(\Sigma (t), \Gamma (t))$ to be
$(\sigma _i(t), \inf \Gamma )$ (resp., $(\sup \Sigma ,\gamma _j(t))$)
where, if $\Sigma _i=\epsilon _x$ and so $(a_i,b_i)=\emptyset $
(resp., $\Gamma _j=\epsilon _y$ and so $(c_j,d_j)=\emptyset $) we
set $\sigma _i(t)=x$ (resp., $\gamma _j(t)=y$).
 
\begin{lem}\label{lem pi of an h manifold}
If $(X,\mu ,x_0)$ is a strictly properly $\bigvee $-definable
$H$-manifold, then $\pi _1(X,x_0)$ is an abelian group.
\end{lem}

\pf
For definable paths $\Sigma $ and $\Gamma $ in $X$ let $\Sigma \Gamma
$ be the definable path such that $\Sigma \Gamma (t):=\mu (\Sigma
(t),\Gamma (t)).$ It follows from the definition 
that: $(1)$ $[\epsilon _{x_0}\Sigma]=[\Sigma \epsilon _{x_0}]=[\Sigma
]$; $(2)$ if $[\Sigma ]=[\Sigma ']$ and $[\Gamma ]=[\Gamma ']$ then
$[\Sigma \Gamma ]=[\Sigma '\Gamma ']$; and $(3)$ $(\Sigma \cdot \Gamma
)(\Sigma '\cdot \Gamma ')=(\Sigma \Sigma ')\cdot (\Gamma \Gamma ').$
Now the lemma follows from the definable homotopies:
$[\Sigma \cdot \Gamma ]=[(\Sigma \epsilon _{x_0})\cdot (\epsilon
_{x_0}\Gamma )]=$$[(\Sigma \cdot \epsilon _{x_0})(\epsilon _{x_0}\cdot
\Gamma )]$$=[\Sigma \Gamma ],$
$[\Gamma \cdot \Sigma ]=[(\epsilon _{x_0}\Gamma )\cdot (\Sigma
\epsilon _{x_0})]=$
$[(\epsilon _{x_0}\cdot \Sigma )(\Gamma \cdot \epsilon
_{x_0})]$$=[\Sigma \Gamma ].$ 
\qed

\begin{prop}\label{prop lifting homotopies}
Suppose that $(\bar{X},p,X)$ is a strictly 
properly $\bigvee $-definable covering space. 
Consider the diagram of continuous  strictly properly $\bigvee $-definable maps
\[
\begin{array}{clcr}
Y\,\,\,\,\,\stackrel{\bar{f}}
{\rightarrow }\,\,\,\,\,\,\bar{X}\\
\downarrow ^{j}\,\,\,\,\,\,\,\nearrow ^{\bar{H}}\,\,\,\,\downarrow ^p\\
Y\times |\Gamma |\,\stackrel{H}{\rightarrow }X,\,\,\,
\end{array}
\]
where $j(y):=(y,\inf \Gamma )$ for all $y\in Y$. Then there exists a
unique continuous strictly properly $\bigvee $-definable map $\bar{H}:
Y\times |\Gamma |\into Y$ making the diagram commutative.
\end{prop}

\pf
Since $Y$ is properly $\bigvee $-definably connected then so is
$Y\times |\Gamma |$, and because $\bar{H}(y,\inf \Gamma )=\bar{f}(y)$
for all $y\in Y$, by lemma \ref{lem unique lift}, $\bar{H}$ is unique. 

Note that its clearly enough to assume that $\Gamma $ is a definable
basic path parametrised by say $\gamma :(a,b)\into Z$. And therefore
we can assume as well that $|\Gamma |=[a,b]$.
Let $\{U_l:l\in L\}$ be a $p$-admissible family of open definable
subsets of $X$. 
Since $H$ is continuous strictly properly $\bigvee $-definable, for
each $l\in L$, there is a family $\{V^j_l:j\in J_l\}$ with
$|J_l|<\aleph _1$ of open definable subsets of $Y\times I$ such that
$\{V^j_l:j\in J_l\}=H^{-1}(U_l)$. Let $\pi _1:Y\times [a,b]\into Y$
and $\pi _2:Y\times [a,b]\into [a,b]$ be the natural projections. 

We will first construct $\bar{H}^j_l:V^j_l\into \bar{X}$, for each
$l\in L$ and $j\in J_l$ such that $(1)$ $p\circ
\bar{H}^j_l=H_{|V^j_l}$, $(2)$ for all $(y,a)\in V^j_l$,
$\bar{H}^j_l(y,a)=\bar{f}(y)$ and $(3)$ for all $(y,t)\in
V^{j'}_{l'}\cap V^j_l$, $\bar{H}^j_l(y,t)=\bar{H}^{j'}_{l'}(y,t)$. 
Its then clear that under these conditions,
the collection $\{\bar{H}^j_l:l\in L,\,\,\,j\in J_l\}$ determines the
strictly properly $\bigvee $-definable map $\bar{H}$ satisfying the
proposition.

Let $\{O^n_l:n\in N_l\}$ be the open definable sheets in $\bar{X}$ over
$U_l$, and let $A^n_l=\bar{f}^{-1}(O^n_l)$ in $Y$. For each $n\in
N_l$, let $V^{j,n}_l:=V^j_l\cap \pi _1^{-1}(A^n_l)$. Then $V^{j,n}_l$
is a disjoint cover of $V^j_l$ by open definable sets. Define
$\bar{H}^{j,n}_l:V^{j,n}_l\into O^n_l$ by
$\bar{H}^{j,n}_l:=(p_{|O^n_l})^{-1}\circ H$. 
We have thus constructed $\bar{H}^{j,n}_l:V^{j,n}_l\into \bar{X}$, for
each $l\in L$, $j\in J_l$ and $n\in N_l$ such that $(1)$ $p\circ
\bar{H}^{j,n}_l=H_{|V^{j,n}_l}$, $(2)$ for all $(y,a)\in V^{j,n}_l$, 
$\bar{H}^{j,n}_l(y,a)=\bar{f}(y)$  and $(3)$ for all $(y,t)\in
V^{j,n}_{l}\cap V^{j',n'}_l$, 
$\bar{H}^{j,n}_l(y,t)=$ $\bar{H}^{j',n'}_{l'}(y,t)$. 

The collection $\{\bar{H}^{j,n}_l:n\in N_l\}$ clearly determines
$\bar{H}^j_l$ satisfying $(1)$, $(2)$ and $(3)$.
\qed

Given a strictly properly $\bigvee $-definable $H$-manifold $(X,\mu
,x_0)$ (resp., given a strictly properly $\bigvee $-definable $H$-group
$(X,\mu ,\iota ,x_0)$) we define in the obvious way the notion of
strictly properly $\bigvee $-definable $H$-submanifold (resp.,
strictly properly $\bigvee $-definable $H$-subgroup, strictly properly
$\bigvee $-definable $H$-centre $Z(X)$ of $X$, etc.,). 

\begin{thm}\label{thm cover of h manifold}
Suppose that $((Y,y_0),p,(X,x_0))$ is a strictly 
properly $\bigvee $-definable covering
space and $(X,\mu ,x_0)$ is a strictly properly $\bigvee $-definable
$H$-manifold. 

(1) There is a unique structure $(Y,\gamma ,y_0)$ of a strictly
properly $\bigvee $-definable $H$-manifold on $Y$ such that the
diagram below commutes.
\[ 
\begin{array}{clrc}
Y\times Y\,\,\,\,\stackrel{\gamma }{\rightarrow
}\,\,\,\,Y\\
\,\,\,\,\,\,\downarrow ^{p\times
p}\,\,\,\,\,\,\,\,\,\,\,\,\,\,\,\,\,\,\downarrow ^p\\
X\times X\,\,\,\,\stackrel{\mu }{\rightarrow }\,\,\,\,X.
\end{array}
\]
Moreover, if $(X,\mu ,\iota ,x_0)$ is a (strictly properly $\bigvee
$-definably abelian) strictly properly $\bigvee $-definable $H$-group
then there is a unique structure $(Y,\gamma ,\zeta ,y_0)$ of a
(strictly properly $\bigvee $-definably abelian) strictly properly
$\bigvee $-definable $H$-group such that the following diagram
\[
\begin{array}{clrc}
Y\,\,\,\,\stackrel{\zeta}{\rightarrow }\,\,\,\,Y\\
\,\,\downarrow ^p\,\,\,\,\,\,\,\,\,\,\,\,\,\,\,\downarrow ^p\\
X\,\,\,\,\stackrel{\iota }{\rightarrow }\,\,\,\,X
\end{array}
\]
is commutative.

(2) With respect to $\gamma $, $p^{-1}(x_0)$ is an abelian group
isomorphic with 
$$\pi _1(X,x_0)/p_*(\pi _1(Y,y_0))\simeq Cov(Y/X).$$
Moreover, if $(X,\mu ,\iota ,x_0)$ is a strictly properly 
$\bigvee $-definable $H$-group then 
$$(p^{-1}(x_0),\gamma _{|p^{-1}(x_0)}, \zeta _{|p^{-1}(x_0)}, y_0)$$ 
is a 
strictly properly $\bigvee $-definable $H$-subgroup contained in the 
strictly properly $\bigvee $-definable $H$-centre $Z(Y)$ of $Y$.
\end{thm}

\pf
(1) Let $f=\mu \circ (p\times p)$ then by the proof of lemma \ref{lem
pi of an h manifold} we see that $$f_*(\pi _1(Y\times
Y,(y_0,y_0))\subseteq p_*(\pi _1(Y,y_0)).$$
Therefore, by proposition \ref{prop lifting functions} there is a
unique strictly properly $\bigvee $-definable $\gamma $ such that
$\gamma(y_0,y_0)=y_0$ and the diagram in (1) is commutative. To see
that $\gamma $ satisfies the condition of a multiplication of a
strictly properly $\bigvee $-definable $H$-manifold, we consider by
proposition \ref{prop lifting homotopies} the strictly
properly $\bigvee $-definable liftings $\tilde{H}_ j$ ($j=1,2$) of the
corresponding strictly properly $\bigvee $-definable homotopies $H_j$
for $[\mu \circ i_j]=[1_X]$ in the strictly properly $\bigvee
$-definable $H$-manifold $X$. The rest of (1) is proved in a similar way.

(2) By lemma \ref{lem pi of an h manifold}, $\pi _1(X,x_0)$ is
abelian. Therefore, $(Y,p,X)$ is regular by corollary \ref{cor
regular} and $Cov(Y/X)\simeq \pi _1(X,x_0)/p_*(\pi _1(Y,y_0).$ We define a
bijection $\psi :\pi _1(X,x_0)\into p^{-1}(x_0)$ by $\psi ([\Sigma
])=y_0*\Sigma .$ By the diagram in (1), $(p^{-1}(x_0),
\gamma_{|p^{-1}(x_0)}, y_0)$ is a strictly properly $\bigvee
$-definable $H$-submanifold. By lemma \ref{lem pi of an h manifold} we
have, $[\Sigma \cdot \Gamma ]$$=[\Sigma \Gamma ]$ and so $y_0*(\Sigma
\cdot \Gamma)=y_0*(\Sigma \Gamma ).$ Therefore, $\psi $ is a
homomorphism and from previuos results we see that $ker \psi =p_*(\pi
_1(Y,y_0))$.
For $y\in Kerp$, the strictly properly $\bigvee $-definable left and
right $H$-translations $L_y,R_y:Y\into Y$ given by $L_y(z)=\gamma
(z,y)$ and $R_y(z)=\gamma (y,z)$ respectively satisfy $[p\circ
L_y]=[p]$$=[p\circ R_y]$ and $[L_y(y_0)]$$=[\epsilon _y]=[R_y(y_0)]$,
hence $[L_y]=[R_y]$ by proposition \ref{prop lifting homotopies}.
\qed

Similarly we get the following result:

\begin{cor}\label{cor h manifolds}
Suppose that $((Y_i,y_{i,0}),p_i,(X_i,x_{i,0}))$ for $i=1,2$ 
are properly $\bigvee $-definable covering
spaces and $(X_i,\mu _i,x_{i,0})$ (resp., $(X_i,\mu _i,\iota _i,x_{i,0})$)
are strictly properly $\bigvee $-definable $H$-manifolds (resp.,
$H$-groups) with a strictly properly 
$\bigvee $-definable $H$-map (resp., $H$-homomorphism)
$h:(X_1,x_{1,0})\into (X_2,x_{2,0})$. Then there are 
unique structures $(Y_i,\gamma _i,y_{i,0})$ (resp., $(Y_i,\gamma
_i,\zeta _i,y_{i,0})$ of strictly
properly $\bigvee $-definable $H$-manifolds (resp., $H$-groups) 
on $Y_i$ and a unique strictly properly 
$\bigvee $-definable $H$-map (resp., $H$-homomorphism)
$l:(Y_1,y_{1,0})\into (Y_2,y_{2,0})$ such that the diagram below
commutes
\[ 
\begin{array}{clrc}
(Y_1,y_{1,0})\,\,\,\,\stackrel{l }{\rightarrow }\,\,\,\,(Y_2,y_{2,0})\\
\downarrow ^{p_1}\,\,\,\,\,\,\,\,\,\,\,\,\,\,\,\,\,\,\,\,\,\,\,\,\,
\downarrow ^{p_2}\\
(X_1,x_{1,0})\,\,\,\,\stackrel{h }{\rightarrow }\,\,\,\,(X_2,x_{2,0}).
\end{array}
\] 
\end{cor}

\end{subsection}

\end{section}

\begin{section}{Strictly properly $\bigvee $-definable groups}
\label{section strictly properly bigvee definable groups}

In this section, ${\mathbf X}$$=(X, (X_i,\phi _i)_{i\in I})$ and 
${\mathbf Y}$$=(Y, (Y_j,\psi _j)_{j\in J})$ will be (properly) 
$\bigvee $-definable manifolds and $Z\subseteq X$ will be a (properly) 
$\bigvee $-definable subset of $X$.  Note that we do not assume that
$Z$, $X$ or $Y$ have definable choice or are (properly) $\bigvee $-definably 
connected.
 
\begin{subsection}{Strictly properly $\bigvee $-definable groups}
\label{subsection strictly properly bigvee definable groups}

\begin{defn}\label{defn bigvee definable groups}
{\em
A {\it strictly (properly) $\bigvee $-definable group} is a group
$(Z,\mu ,\iota ,z_0)$ on $Z$ with identity $z_0$ and such that the
product map $\mu $ and the inverse map $\iota $ are 
strictly (properly) $\bigvee $-definable maps. {\it Strictly (properly)
$\bigvee $-definable rings} are defined in a similar way. 
}
\end{defn}

\begin{lem}\label{lem large in a group}
Suppose that $Z$ is a strictly (properly) $\bigvee $-definable group
and let $V$ be a large (properly) $\bigvee $-definable subset of $Z$. 
Then countably many translates of $V$ cover $Z$.
\end{lem} 

\pf
This is just like in lemma 2.4 in \cite{p1}: Let $\M$$\prec \N$ be a
small model over which $(Z,\mu ,\iota ,z_0)$ and $V$ are defined, and
assume without loss of generality that $\N$ is $\aleph
_1$-saturated. Let $i\in I^Z$, $a\in Z_i$ and let $c\in Z_i$ be a
generic point of $Z$ over $M$ such that $tp(c/Ma)$ is finitely
satisfiable in $M$. Then $c$ is a generic point of $Z$ over $Ma$ 
(see proof of lemma 2.4 in \cite{p1}). Since $V$ is a large properly 
$\bigvee $-definable subset of $Z$, so is $\mu (V,\iota (a))$ and
therefore, $c\in \mu (V_k,\iota (a))$ and $a\in \mu (\iota (c),V_k)$
for some $k\in K_i$ where $K_i$ a finite subset of $I^V$ such that 
$Z_i\subseteq \cup \{\mu (\iota (Z_i),V_k):k\in K_i\}$ (this exists 
because $\mu $ and $\iota $ are strictly (properly) $\bigvee
$-definable). Since $tp(c/Ma)$ is finitely satisfiable over $M$, there
is $b\in Z_i(M)$ such that $a\in  \mu (b,V_k)$ for some $k\in K_i$. 
Therefore, by compactness theorem, for each $i\in I^Z$, there are 
$b_1,\dots b_{r_i}\in Z_i(M)$ such that for every $a\in Z_i$, $a\in 
\mu (b_j,V_k)$ for some $j=1,\dots ,r_i$ and $k\in K_i$.
\qed

Lemma \ref{lem large in a group} together with the properly 
$\bigvee $-definable cell decomposition theorem, gives similarly to
what happens in the definable case (see \cite{p1} and \cite{pps1}) the 
following result for strictly (properly) $\bigvee $-definable groups and 
rings. The version of this result (also included in the statement of
theorem \ref{thm pvd groups}) for strictly $\bigvee $-definable groups
and rings, which shows that these groups and rings are strictly 
$\bigvee $-definable topological groups and rings appears in \cite{pst2}.

\begin{thm}\label{thm pvd groups}
Let $Z$ be a strictly (properly) $\bigvee $-definable group 
(resp., ring) and $W$  a strictly (properly) $\bigvee $-definable 
subgroup (resp., left or right ideal). Then there are unique
structures of properly $\bigvee $-definable manifolds on $Z$ and $W$ such that
$Z$ (resp., $W$) is a strictly (properly) submanifold of $X$ (resp.,
$Z$) and  the group (resp., ring) operations are continuous (in fact
$C^p$) strictly properly $\bigvee $-definable maps. Any strictly
(properly) $\bigvee $-definable homomorphism between strictly
(properly) $\bigvee $-definable groups (resp., rings) is a continuous
(in fact $C^p$) strictly (properly) $\bigvee $-definable
homomorphism. Moreover, $W$ is closed in $Z$, $W$ is open in $Z$ iff 
$dimW=dimZ$ and, when both $Z$ and $W$ are strictly properly 
$\bigvee $-definable groups then $dimW=dimZ$ iff $W$ has countable
index in $Z$. 
\end{thm}

\pf
For each $i\in I^Z$ we have $Z_i=\cup \{Z^s_i:s\in S_i\}$ where 
$|S_i|<\aleph _1$ (in fact $|S_i|=1$ if $X$ is properly $\bigvee
$-definable) and each $Z_i^s$ is a definable subset of $X_i$. Suppose
that $dimZ=k$. Then by cell decomposition theorem, each $Z^s_i$ is a
finite union of cells. For each $i\in I^Z$ and $s\in S_i$, let $V^s_i$
be the union of the cells of $Z^s_i$ in the cell decomposition of
$Z^s_i$, of dimension $k$. Note that each $V^s_i$ is a disjoint union
of finitely many  definable subsets $U^s_{i,1}, \dots , U^s_{i,n_i}$ 
definably homeomorphic to an open definable subset of $N^k$ and 
$V=\cup \{V^s_i:i\in I^Z,\,\,\, s\in S_i\}$ is a large open (properly) 
$\bigvee $-definable subset of $Z$. Using the same argument as in the
proof of proposition 2.5 in \cite{p1}, we can further assume that: the 
inverse map is a continuous strictly (properly) $\bigvee $-definable
map from $V$ into $V$; there is a large (properly) $\bigvee
$-definable subset $U$ of $Z\times Z$ such that $U$ is open and dense
in $V\times V$, multiplication is a continuous strictly (properly) 
$\bigvee $-definable map from $U$ into $V$ and for any $a\in V$, if
$b$ is a generic of $V$ over $a$, then $(b,a)\in U$ and $(\iota (b),
\mu (b,a))\in U$. The rest of the arguments in \cite{p1}
show that $U^s_{i,j}$'s give $Z$ the structure of a properly 
$\bigvee $-definable manifold such that the group (resp., ring) 
operations are continuous (in fact $C^p$) strictly properly 
$\bigvee $-definable maps. Lemma 2.6 in \cite{pst2} shows that the
structures of properly $\bigvee $-definable manifolds on $Z$ and $W$
such that the group (resp., ring) operations are continuous (in fact
$C^p$) strictly properly $\bigvee $-definable maps and $Z$ (resp.,
$W$) is a strictly (properly) submanifold of $X$ (resp., $Z$) are
unique (i.e., the inclusion maps are strictly (properly) 
$\bigvee $-definable homeomorphisms onto their image). The result about
strictly (properly) $\bigvee $-definable homomorphism is proved in
lemma 2.8 in \cite{pst2} and lemma 2.6 in \cite{pst2} shows that $W$
is closed in $Z$ and $W$ is open in $Z$ iff $dimW=dimZ$.

Suppose that both $W$ and $Z$ are properly $\bigvee $-definable. Then 
clearly, if $W$ has countable index in $Z$ the $dimW=dimZ$. Suppose
that $dimW=dimZ$. Looking at $Z$ as a properly $\bigvee $-definable
manifold and $W$ as an open properly $\bigvee $-definable subset, we
see that $Z$ is a countable (disjoint) union of cosets of $W$ in $Z$. 
\qed

For general strictly $\bigvee $-definable groups, even though we have
that, if $W$ has countable index in $Z$ then $dimW=dimZ$ the
reciprocal does not hold: take $\N$$=(N,0,+,<)$ an $\aleph
_1$-saturated extension of $(\RR$$,0,+,<)$, $Z=(N,0,+)$ and $W$ the
convex hull of $\RR$ in $Z$.

\begin{nrmk}
{\em
We will from on  assume that if $(Z,\mu ,\iota ,z_0)$ is a strictly
(properly) $\bigvee $-definable group (resp., ring), then ${\mathbf Z}$$=(Z,
(Z_k,\tau _k)_{k\in K})$ is the corresponding unique properly
$\bigvee $-definable manifold on $Z$ given by theorem \ref{thm pvd
groups}. Since $Z$ is then a properly $\bigvee $-definable subset of
$Z$ (relative to ${\mathbf Z}$) we will simply say that $Z$ be a strictly 
properly $\bigvee $-definable group (resp., ring). 
By a definable (resp., properly $\bigvee $-definable,
$\bigvee $-definable) subset (resp., subgroup, ideal, etc.,) of $Z$, we will
mean a definable (resp., properly $\bigvee $-definable,
$\bigvee $-definable) subset (resp., subgroup, ideal, etc.,) of $Z$
relative to ${\mathbf Z}$. Finally, as usual we will some times write
$xy$ for $\mu (x,y)$ and $x^{-1}$ for $\iota (x)$.
}
\end{nrmk}

\begin{cor}\label{cor dcc}
Let $Z$ be a strictly properly $\bigvee $-definable group (resp.,
ring). Then the properly $\bigvee $-definable connected component
$Z^0$ of $Z$ is the smallest strictly properly $\bigvee $-definable 
subgroup (resp., ideal) of $Z$ of countable index. If $\{Z^s:s\in S\}$
is a decreasing sequence of strictly properly $\bigvee $-definable 
subgroups (resp., left or right ideals) of $Z$ then $\cap \{Z^s:s\in
S\}$$=\cap \{Z^s:s\in S_0\}$ for some $S_0\subseteq S$ with
$|S_0|<\aleph _1$ and is a strictly properly $\bigvee $-definable
subgroup (resp., left or right ideal) of $Z$.
\end{cor}

\pf
The first part is clear. So suppose that $\{Z^s:s\in S\}$ is a
decreasing sequence of strictly properly $\bigvee $-definable
subgroups of $Z$. For each $s\in S$, let $k_s:=dimZ^s$. Since 
$\{k_s:s\in S\}\subseteq \{0,\dots ,dimZ\}$, there are $k_1<\dots
<k_m$ in $\{0,\dots ,dimZ\}$ and there are disjoint subsets $S_1,\dots
,S_m$ of $S$ such that $S=S_1\cup \cdots \cup S_m$ and for each 
$l\in \{0,\dots ,m\}$, if $s\in S_l$ then $dimZ^s=k_l$. Therefore, 
since we want to determine $\cap \{Z^s:s\in S\}$, we may assume
without loss of generality that, for all $s\in S$, $dimZ^s=r$. It
follows from the first part, that for all $s\in S$, the properly 
$\bigvee $-definable connected component of $Z^s$ is the same properly 
$\bigvee $-definable subgroup $V$. Let $s_0$ be the first element of
$S$ (we can assume, without loss of generality that $s_0$
exists). Then there is a decreasing sequence $\{U_s:s\in S\}$ of
properly $\bigvee $-definable countable subsets of $Z^{s_0}$
containing the identity element and such that for each $s\in S$, 
$Z^s=\cup \{uV:u\in U_s\}$. Note that, there is $S_0\subseteq S$ such
that $|S_0|<\aleph _1$ and $\{U_s:s\in S\}$$=\{U_s:s\in S_0\}$. Let 
$U=\cap \{U_s:s\in S\}$. Then $U$ is a countable nonempty (contains
the identity) properly $\bigvee $-definable subset of $Z^{s_0}$.  Let 
$W:=\cup \{uV:u\in U\}$. Then $W$ is a strictly properly 
$\bigvee $-definable subgroup of $Z$ such that $W=\cap \{Z^s:s\in S\}$. 

The results for the strictly properly $\bigvee $-definable rings
follows from the corresponding results for strictly properly 
$\bigvee $-definable groups.
\qed

From (the proof of) corollary \ref{cor dcc} we  easily get the
following result.

\begin{cor}\label{cor centralizer}
Let $Z$ be a strictly properly $\bigvee $-definable group and let 
$S\subseteq Z$. Then $C_Z(S)=\{z\in Z:\forall s\in S,\,\, zs=sz\}$, 
the centraliser of $S$ in $Z$, is a strictly properly 
$\bigvee $-definable subgroup. In fact there is $S_0\subseteq S$ such
that $|S_0|<\aleph _1$ and $C_Z(S)=C_Z(S_0)$, and  for each 
$k\in K^{C_Z(S)}$, there is a finite subset $A_k$ of $S$ such that 
$C_Z(S)\cap Z_k=C_Z(A_k)\cap Z_k$.
In particular, if $A$ is a subgroup of $Z$, then the centre of
$C_Z(A)$ is a strictly properly $\bigvee $-definable subgroup of $Z$ 
containing $A$, which is normal if $A$ is normal.
\end{cor}

Similarly to corollary 2.15 in \cite{p1} (see also proposition 5.6 in 
\cite{p2}) we get:

\begin{cor}\label{cor infinite abelian subgroup}
Let $Z$ be an infinite strictly properly $\bigvee $-definable group. 
Then $Z$ has an infinite strictly properly $\bigvee
$-definable abelian subgroup.  
\end{cor}

We finish this subsection with a result that generalises a theorem
from \cite{e1} on definable groups. 

\begin{thm}\label{thm quotient}
Let $X$ be a strictly properly $\bigvee $-definable group and let $Z$
be a normal strictly properly $\bigvee $-definable subgroup of
$X$. Then we have  strictly properly $\bigvee $-definable extension 
$1\rightarrow Z\rightarrow X\stackrel{j}{\rightarrow }Y\rightarrow 1$
of strictly properly $\bigvee $-definable groups, with strictly
properly $\bigvee $-definable section $s:Y\into X$.
\end{thm}

\pf
This is proved by adapting the corresponding result from \cite{e1} for 
definable groups. We first find $Y$. The argument in \cite{e1} shows
that for each $i\in I$, there is a large definable subset $U_i$ of
$X_i$ and there are definable functions $l_{i,1},\dots
,l_{i,m}:U_i\into N$ such that for each $x\in U_i$, there is 
$z\in xV\cap X_i$ with $\phi _i(z)=(l_{i,1}(x),\dots ,l_{i,m}(x))$ and
for all $y\in U_i$, if $xV=yV$ then $(l_{i,1}(x),\dots ,l_{i,m}(x))=$
$(l_{i,1}(y),\dots ,l_{i,m}(y))$.
Let $U=\cup \{U_i:i\in I\}$. Then $U$ is a large properly 
$\bigvee $-definable subset of $X$. Let $I=\{0,1,\dots \}$ be an
enumeration of $I$ and let $l:U\into X$ the strictly properly 
$\bigvee $-definably map defined inductively in the following way: 
for $x\in U_0$, $l(x):=\phi _i^{-1}(l_0(x))$; suppose that $l$ has 
been defined on $U_0\cup \cdots \cup U_k$, then we define $l$ on 
$U_{k+1}$ by $l(x):= \phi _{k+1}^{-1}(l_{k+1}(x))$ if 
$xV\cap (U_0\cup \cdots \cup U_k)=\emptyset $ or $l(x):=l(y)$ for some 
(for all) $y\in xV\cap (U_0\cup \cdots \cup U_k)$. Clearly, $l:U\into
X$ is well defined. Since $U$ is a large properly $\bigvee $-definable 
subset of $X$, there is a subset $\{x_s:s\in S\}$ of $X$ with 
$|S|<\aleph _1$ such that $X=\cup \{x_sU: s\in S\}$. Let 
$S=\{0,1,\dots \}$ be an enumeration of $S$ and let $L:X\into X$ be
the strictly properly $\bigvee $-definable map defined inductively in 
the following way: for $x\in x_0U$,  $L(x):=x_0l(x_0^{-1}x)$; suppose 
that $L$ has been defined on $x_0U\cup \cdots \cup x_kU$, then we
define $L$ on $x_{k+1}U$ by $L(x):=x_{k+1}l(x_{k+1}^{-1}x)$ if 
$xV\cap (x_0U\cup \cdots \cup x_kU)=\emptyset $ or $L(x):=L(y)$ 
for some (for all) $y\in xV\cap (x_0U\cup \cdots \cup x_kU)$. Clearly, 
$L:X\into X$ is well defined and moreover, for $x,y\in X$, $L(x)=L(y)$
iff $xV=yV$. Take $Y:=L(X)$ and $j=L:X\into Y$. It turns out that $Y$
is a properly $\bigvee $-definable subset of $X$ with a strictly
properly $\bigvee $-definable group structure given by 
$xy=L(L^{-1}(x)L^{-1}(y))$. 

To find the strictly properly $\bigvee $-definable section $s:Y\into
X$ we use the same argument, just like in the definable case (see \cite{e1}).  
\qed

\end{subsection}

\begin{subsection}{The centerless case}\label{subsection the centerless case}

\begin{defn}\label{defn bigvee semisimple}
{\em
We say that a strictly properly $\bigvee $-definable group $X$ is 
{\it properly $\bigvee $-definably semisimple} if $X$ has no 
strictly properly $\bigvee $-definable normal abelian subgroup of 
dimension bigger than zero.
}
\end{defn}

In particular, a strictly properly $\bigvee $-definable group has
centre of dimension zero.  The following lemma will be useful later.

\begin{lem}\label{lem discrete centre}
Let $X$ be a properly $\bigvee $-definably connected strictly properly 
$\bigvee $-definable group. Then every strictly properly 
$\bigvee $-definable normal subgroup of $X$ of dimension zero is 
contained in $Z(X)$ and if $Z(X)$ has dimension zero then, $X/Z(X)$ is 
a centerless strictly properly $\bigvee $-definable group.
\end{lem}

\pf
This is proved in the same way as in the definable case (see \cite{e1}).
\qed

\begin{thm}\label{thm bigvee semisimple}
If $X$ is a centerless, properly $\bigvee $-definably
semisimple, properly $\bigvee $-definably connected strictly properly 
$\bigvee $-definable group, then $X=X_1\times \cdots \times X_l$ and 
for each $k\in \{1,\dots ,l\}$, there is a definable real closed field 
$R_k$ such that there is no definable bijection between a distinct
pair among the $R_k$'s, and  there is an 
$R_k$-semialgebraically connected, $R_k$-semialgebraic subgroup $G_k$
of $GL(n_k,R_k)$ which is a direct product of $R_k$-semialgebraically
simple, $R_k$-semialgebraic subgroups of $GL(n_k,R_k)$ such that $X_k$
is strictly properly $\bigvee $-definably isomorphic to a $\bigvee
$-definable open and closed subgroup of $G_k$.
\end{thm}

\pf
The proof is a modification of the corresponding result in \cite{pps1}
for definably semisimple groups. We will therefore assume the readers 
familiarity with the terminology of \cite{pps1}.

Arguing as in the proof of theorem 3.1 in \cite{pps1} and using the
fact that the centraliser $C_X(U)$ in $X$ of any subset $U$ of $X$ is
a strictly properly $\bigvee $-definable subgroup of $X$ (see
corollary \ref{cor dcc}), we see that $X$ is a direct product of
strictly properly $\bigvee $-definable unidimensional subgroups. So we
may assume that $X$ is unidimensional.  Further, we may also assume
that $\N$ is $\aleph _1$-saturated and just like in \cite{pps1}, there
is an open transitive interval $M$ such that $e=(d,\dots ,d)$ for some 
$d\in M$, where $e=\phi _{i_0}(x_0)$, $x_0$ is the identity of $X$ and 
$x_0\in X_{i_0}$ and moreover, if $B= M^n$ where $n=dimX$ then 
$\phi _{i_0}^{-1}(B)$ is an open definable neighbourhood of $x_0$. Let 
$\rho :M\into B$ be the continuous injection defined as 
$\rho (x)=(x,d,\dots ,d)$. Let $M^{+}:=\{b\in M:b>d\}$, for 
$b\in M^{+}$ let $M_b :=\{c\in M: d<c<b\}$, and let
$Y_b:=C_X(\bar{M}_b)$ where $\bar{M}_b:=\phi _{i_0}^{-1}(\rho
(M_b)))$. Clearly, $\{Y_b:b\in M^{+}\}$ is a sequence of strictly
properly $\bigvee $-definable subgroups of $X$ such that if $b'<b$
then $Y_b\subseteq Y_{b'}$ and therefore, by corollary \ref{cor dcc}, 
$\{C_X(Y_b):b\in M^{+}\}$  is a sequence of strictly properly 
$\bigvee $-definable subgroups of $X$ such that if $b'<b$ then
$C_X(Y_{b'})\subseteq C_X(Y_b)$. Let $Y=\cup \{Y_b:b\in M^{+}\}$. Then 
$C_X(Y)=\cap \{C_X(Y_b):b\in M^{+}\}$. Hence by corollary \ref{cor
dcc}, there is subset $\{b_s:s\in S\}$ of $M^{+}$ with $|S|<\aleph _1$
and such that $C_X(Y)=\cap \{C_X(Y_{b_s}):s\in S\}$. By saturation,
there is $b\in M^{+}$ such that for all $s\in S$, $b>b_s$. Therefore, 
$C_X(Y)=C_X(Y_b)$ and $Y=Y_b$ (since $\bar{M}_b\subseteq C_X(Y_b)$
$=C_X(Y)$, we have $Y\subseteq C_X(\bar{M}_b)=Y_b$).
Since $X$ is centerless and properly $\bigvee $-definably connected, 
$dimY<dimX$ (otherwise, $Y=X$ and $\bar{M}_b\subseteq Z(X)$). Hence, 
$B$ cannot be covered by finitely many left cosets of $Y$, and arguing 
as in \cite{pps1}, there is a definable real closed field $R$ on some 
open subinterval of $M$. Furthermore, just like in the proof of 
theorem 3.2 in \cite{pps1} we see that $Ad:X\into GL(n,R)$, where 
$Ad(x)$ is the differential at $x_0$ of the strictly properly 
$\bigvee $-definable automorphism $a(x):X\into X$ given by 
$a(x)(z):=xzx^{-1}$, is a strictly properly $\bigvee $-definable 
injective homomorphism.  Let $G:=Aut({\mathbf x})$$<GL(n,R)$ where 
${\mathbf x}$ is the Lie algebra of $X$. By theorem 2.37 in
\cite{pps1}, $G$ and $G^0$ are $R$-semialgebraically semisimple 
$R$-semialgebraic groups of dimension $dimX$, and clearly, by theorem
\ref{thm pvd groups} $Ad(X)$ 
is an open and closed $\bigvee $-definable subgroup of $G$ and since $X$ is
properly $\bigvee $-definably connected, $Ad(X)$ 
is an open and closed strictly properly $\bigvee $-definable subgroup
of $G^0$.  
\qed

\begin{cor}\label{cor dim 1 case}
Let $X$ be a properly $\bigvee $-definably connected, strictly
properly $\bigvee $-definable group such that $dimX=1$. Then, $X$ is 
abelian, divisible and either $X$ is torsion-free and properly 
$\bigvee $-definably ordered or $X$ is a definably compact definable group. 
\end{cor}

\pf
Suppose that $X$ is not abelian. Then every strictly properly 
$\bigvee $-definable subgroup of $X$ has dimension zero. In
particular, 
$dimZ(X)=0$ and by lemma \ref{lem discrete centre}, $X/Z(X)$ is a 
centerless, strictly properly $\bigvee $-definably semisimple strictly 
properly $\bigvee $-definable group of dimension one. But by theorem 
\ref{thm bigvee semisimple} we get a contradiction. The rest, follows 
by adapting the proofs of the corresponding results for definable
groups (see \cite{r} and \cite{PiS1}).
\qed

\end{subsection}

\begin{subsection}{The solvable case}\label{subsection the solvable case}

Recall from \cite{e1} that a definable abelian group $U$ {\it has no 
definably compact parts} if there are definable subgroups 
$1=U_0<U_1<\cdots <U_n=U$ such that for each $j\in \{1,\dots ,n\}$, 
$U_j/U_{j-1}$ is a one-dimensional definably connected, torsion-free 
definable group; and a definable solvable group $U$ {\it has no 
definably compact parts} if there are definable subgroups 
$1=U_0\trianglelefteq U_1\trianglelefteq \cdots \trianglelefteq U_n=U$ 
such that for each $j\in \{1,\dots ,n\}$, $U_j/U_{j-1}$ is a definable 
abelian group with no definably compact parts. Definable solvable
groups with no definably compact parts are classified in
\cite{e1}. The next result (theorem \ref{thm solvable case}), uses
this fact to reduce the classification of strictly properly 
$\bigvee $-definable solvable groups to the the classification of
properly $\bigvee $-definably complete such groups. 

Note that, theorem \ref{thm quotient} makes possible to develop group 
extension theory and group
cohomology theory in the category of strictly properly 
$\bigvee $-definable groups with strictly properly $\bigvee $-definable 
homomorphisms, just like in the category of definable groups with
definable homomorphisms treated in \cite{e1}. The proof of the next
theorem will use this theory, we therefore assume the readers
familiarity the corresponding results from \cite{e1}.

\begin{thm}\label{thm solvable case}
If $X$ is a strictly properly $\bigvee $-definable solvable
group. Then we have a strictly properly $\bigvee $-definable extension 
$1{\rightarrow }Z{\rightarrow }X{\rightarrow }Y{\rightarrow }1$, where
$Z$ is a definable solvable group with no definably compact parts and
$Y$ is a properly $\bigvee $-definably complete, strictly properly 
$\bigvee $-definable solvable group.
\end{thm}

\pf
The proof is just like in the definable case (\cite{e1}) and is based
in the main result of \cite{ps}. We therefore, need to show the
analogue of  the main result of \cite{ps} in the our more general
context. Suppose that $X$ is not properly $\bigvee $-definably
complete, and let $\sigma :(a,b)\into X$ be a definable injective map
such that $\lim _{t\into b^-}\sigma (t)$ does not exist in $X$. Let 
$I:=\sigma (a,b)$ with the natural order $<$, for $b\in I$ let 
$I^{>b}:=\{x\in I:x>b\}$, and for each $x\in X$, let $xI:=\{xt:t\in
I\}$. As in \cite{ps}, define a properly $\bigvee $-definable
relation $\prec _I$ on $X$ by $x\prec _Iy$ iff for all $t\in I$, for
all $V$ definable open neighbourhood of $y$, there $s\in I$ and $U$ a 
definable open neighbourhood of $x$ such that $UI^{>s}\subseteq
VI^{>t}$. And let $\thicksim _I$ be defined by $x\thicksim _Iy$ iff 
$x\prec _Iy$ and $y\prec _Ix$. Arguing just like in \cite{ps}, we see 
that: $(i)$ $\thicksim _I$ is a properly $\bigvee $-definable
equivalence relation on $X$; $(ii)$ the class $X_1$ of the identity $1$
of $X$ is a strictly properly $\bigvee $-definable subgroup and
$(iii)$ the equivalence classes of $\thicksim _I$ are exactly the left 
cosets of $X_1$.  Moreover, lemma 3.7 \cite{ps} shows that $X_1$ is
in fact definable with dimension less than or equal to one, lemma 3.8 
\cite{ps} shows that $X_1$ has dimension one and lemma 3.9 \cite{ps} 
shows that $X_1$ is torsion-free.
\qed

\end{subsection}

\begin{subsection}{Covers of strictly properly $\bigvee $-definable groups}
\label{subsection covers of groups}

We are now ready to prove the main results of the paper, theorem
\ref{thm main result1} and theorem \ref{thm main result2} below. But
first we need the following lemma. 

\begin{lem}\label{lem definable choice for groups}
Let $X$ be a strictly properly $\bigvee $-definable group (resp.,
ring). Then $X$ has strong definable choice.
\end{lem}

\pf
Let $R(X)$ be the maximal, properly $\bigvee $-definably connected, 
strictly properly $\bigvee $-definable solvable normal subgroup of
$X$. Then, since we have a strictly properly $\bigvee $-definable
extension $1\rightarrow R(X)\rightarrow X\rightarrow Y\rightarrow 1$
where $Y$ is properly $\bigvee $-definably semisimple, which by
theorem \ref{thm bigvee semisimple} has strong definable choice, its
enough to show that $R(X)$ has strong definable choice. One the
another hand, by theorem \ref{thm solvable case}, and a similar
argument, its enough to show that a properly $\bigvee $-definably
complete, strictly properly $\bigvee $-definable group has strong
definable choice. But this can be proved using the same argument used
in \cite{e1} to prove that a definably compact definable group has
strong definable choice.
\qed

\begin{thm}\label{thm main result1}
Let $(X,\mu ,\iota ,x_0)$ be a  
strictly properly $\bigvee $-definable group and suppose that 
$((Y,y_0),p,(X,x_0))$ is a strictly properly $\bigvee $-definable covering
space. 

(1) Then there is a unique structure $(Y,\gamma ,\zeta ,y_0)$ of a
strictly properly $\bigvee $-definable group on $Y$ such that the diagram
\[ 
\begin{array}{clrc}
Y\times Y\,\,\,\,\stackrel{\gamma }{\rightarrow
}\,\,\,\,Y\\
\,\,\,\,\,\,\downarrow ^{p\times
p}\,\,\,\,\,\,\,\,\,\,\,\,\,\,\,\,\,\,\downarrow ^p\\
X\times X\,\,\,\,\stackrel{\mu }{\rightarrow }\,\,\,\,X.
\end{array}
\]
and the diagram
\[
\begin{array}{clrc}
Y\,\,\,\,\stackrel{\zeta}{\rightarrow }\,\,\,\,Y\\
\,\,\downarrow ^p\,\,\,\,\,\,\,\,\,\,\,\,\,\,\,\downarrow ^p\\
X\,\,\,\,\stackrel{\iota }{\rightarrow }\,\,\,\,X
\end{array}
\]
are commutative.

(2) Moreover, 
$(p^{-1}(x_0), \gamma _{|p^{-1}(x_0)}, \zeta _{|p^{-1}(x_0)}, y_0)$
is an abelian strictly properly $\bigvee $-definable subgroup
isomorphic with $$\pi _1(X,x_0)/p_*(\pi _1(Y,y_0))\simeq Cov(Y/X)$$ and
contained in the centre $Z(Y)$ of $Y$.
\end{thm}

\pf
This result follows immediately
from theorem \ref{thm cover of h manifold}, but it can also be proved
directly using the same argument and proposition \ref{prop lifting
functions} instead of proposition \ref{prop lifting homotopies}. 
\qed

\begin{thm}\label{thm main result2}
If $((Y_i,y_{i,0}),p_i,(X_i,x_{i,0}))$ are strictly properly 
$\bigvee $-definable covering spaces, $(X_i,\mu _i,\iota
_i,x_{i,0})$ are strictly properly 
$\bigvee $-definable groups and $h:(X_1,x_{1,0})\into (X_2,x_{2,0})$
is a strictly properly $\bigvee $-definable homomorphism, then there are 
unique structures $(Y_i,\gamma _i,\zeta _i,y_{i,0})$ of strictly
properly $\bigvee $-definable groups on $Y_i$ and there is 
a unique strictly properly 
$\bigvee $-definable homomorphism
$l:(Y_1,y_{1,0})\into (Y_2,y_{2,0})$ such that the diagram below
commutes
\[ 
\begin{array}{clrc}
(Y_1,y_{1,0})\,\,\,\,\stackrel{l }{\rightarrow }\,\,\,\,(Y_2,y_{2,0})\\
\downarrow ^{p_1}\,\,\,\,\,\,\,\,\,\,\,\,\,\,\,\,\,\,\,\,\,\,\,\,\,
\downarrow ^{p_2}\\
(X_1,x_{1,0})\,\,\,\,\stackrel{h }{\rightarrow }\,\,\,\,(X_2,x_{2,0}).
\end{array}
\] 
\end{thm}

\pf
This result is proved in a similar way to theorem \ref{thm main result1}.
\qed

Clearly, there are results analogue to theorem  \ref{thm main
result1} and theorem \ref{thm main result2} for strictly properly
$\bigvee $-definable rings. 

\begin{cor}\label{cor from main}
If $(X,p,Z)$ be a strictly properly $\bigvee $-definable covering
space, then $X$ is properly $\bigvee $-definably complete (resp.,
abelian, nilpotent or solvable) iff $Z$ is properly $\bigvee $-definably 
complete (resp., abelian, nilpotent or solvable).
\end{cor}

\pf
This is a consequence
of theorem \ref{thm main result1} and corollary \ref{cor dcc}.
\qed

\begin{thm}\label{thm quotient by dim0}
Let $X$ be properly $\bigvee $-definably connected, strictly properly 
$\bigvee $-definable group and let $Z$ be a strictly properly 
$\bigvee $-definable normal subgroup of $X$ with $dimZ=0$. Then 
$(X,j,X/Z)$ is a regular strictly properly $\bigvee $-definable 
$Z$-covering space, and 
$$Z\simeq Cov(X/(X/Z))\simeq \pi _1(X/Z,x/Z)/j_*(\pi _1(X,x)).$$ 
In particular, if $h:Y\into X$ is a strictly properly 
$\bigvee $-definable surjective homomorphism such that $dim(Kerh)=0$,
then $(Y,h,X)$ is a regular strictly properly $\bigvee $-definable 
$Kerh$-covering space.
\end{thm}

\pf
By corollary \ref{cor proper action}, its enough to show that $Z$ acts 
properly on $X$. 
For $i\in I$, let $\thicksim _i$ be the definable equivalence relation
on $X_i$ given by $x\thicksim _iy$ iff $xZ\cap X_i=yZ\cap
X_i$ i.e., iff there is $z\in Z_i$ such that $y=zx$ or $x=zy$ (where
$Z_i$ is some finite subset of $Z$, which exists by o-minimality since 
every thing is strictly properly $\bigvee $-definable). Clearly, for
all $x\in X_i$, $|x/\thicksim _i|\leq |Z_i|$ and there are: a positive 
natural number $m_i$, disjoint definable subsets $A_{i,l}$ 
($l=1,\dots ,m_i$) of $X_i$ whose union is $X_i$, and positive natural 
numbers $k_{i,1}>k_{i,2}>\cdots >k_{i,m_i}$ such that for all 
$x\in A_{i,l}$, we have: $|x/\thicksim _i|=k_{i,l}$ and 
$x/\thicksim _i\subseteq A_{i,l}$. Therefore $\thicksim _i$ induces a 
definable equivalence relation on $A_{i,l}$ by restriction and there 
are subsets $Z(i,l)$ of $Z_i$ such that for all $x,y\in A_{i,l}$, 
$x\thicksim _iy$ iff there is $z\in Z(i,l)$ such that $y=zx$ or
$x=zy$. The definable set $\{(x,y)\in A_{i,l}\times X_i: y\in x/
\thicksim _i\}$ is a disjoint union of the diagonal $\Delta _{i,l}$ of 
$A_{i,l}$ together with the disjoint definable sets $R_{i,l,z}$ 
(where $z\in Z(i,l)$). Let $\alpha _i,\beta _i:X_i\times X_i\into X_i$ 
be the definable maps given by $\alpha _i(x,y)=x$ and $\beta
_i(x,y)=y$. Then $R_{i,l,z}$  is the graph of the definable
homeomorphism $\gamma _{i,l,z}:\alpha _i(R_{i,l,z})\into \beta
_i(R_{i,l,z})$ given by $\gamma _{i,l,z}(x):=\beta _i(R_{i,l,z}\cap 
\{(x,y): y=zx\})$ or of the definable homeomorphism $\gamma
_{i,l,z^{-1}}:\alpha _i(R_{i,l,z})\into \beta _i(R_{i,l,z})$ given by 
$\gamma _{i,l,z^{-1}}(x):=\beta _i(R_{i,l,z}\cap \{(x,y): zy=x\})$ . 
Let $A_{i,l,z} :=\alpha _i(R_{i,l,z})$ and let $B_{i,l,z}:=
\beta _i(R_{i,l,z})$. Its clear that there is an open definable subset 
$V_i$ of $X_i$ which is the interior in $X_i$ of a definable set of
the form $A_{i,l_1,z_1}\cup \cdots \cup A_{i,l_{r},z_{r}}$ and such
that: $(1)$ for every $x\in X_i\setminus V_i$ there is $z\in Z_i$ such
that $x\in \bar{zV_i}$; $(2)$ $dim(X_i\setminus \{zV_i:z\in Z_i\})$
$<dimX_i$; and $(3)$  for every $z\in Z_i$, if $z\neq 1$ then
$zV_i\cap V_i$$=\emptyset $. 

Since the properly $\bigvee $-definable subset $\cup \{zV_i:i\in
I,\,\, z\in Z_i\}$ is large in $X$, there is a subset $\{x_s:s\in S\}$
of $X$ with $|S|<\aleph _1$ and such that $X=\cup \{zx_sV_i:i\in I,
\,\, s\in S, \,\, z\in Z_i\}$ (note that by lemma \ref{lem discrete centre}, 
$Z\subseteq Z(X)$). To finish we show that $\{zx_sV_i:i\in I,\,\, s\in
S,\,\, z\in Z_i\}$ is a $Z$-admissible family of definable open
subsets of $X$. In fact, if $u\in Z$ and $uzx_sV_i\cap zx_sV_i\neq
\emptyset $, then $uV_i\cap V_i\neq \emptyset $. Therefore, $u\in Z_i$
or $u^{-1}\in Z_i$ and so $u=1$. 
\qed

\end{subsection}
 
\begin{subsection}{Local strictly properly $\bigvee $-definable isomorphism}
\label{subsection local isomorphisms}

\begin{defn}\label{defn local isomorphism}
{\em
Let $X$ and $Y$ be strictly properly $\bigvee $-definable groups, $U$
an open properly $\bigvee $-definable neighbourhood of the identity in
$X$, and $f:U\into Y$ a strictly properly $\bigvee $-definable map. We
say that $f$ is a {\it locally strictly properly $\bigvee $-definable 
homomorphism} if for all $x,y\in U$ such that $xy\in U$, we have
$f(xy)=f(x)f(y)$ and there is a open properly $\bigvee $-definable
neighbourhood $V$ of the identity of $X$ such that $V^{-1}V\subseteq
U$ and $\{x_lV:l\in L\}$ is an open cover of $X$ with $|L|<\aleph _1$, 
such that for each $i\in I$, there is a finite subset $L_i$ of $L$
with $X_i\subseteq \cup \{x_lV:l\in L_i\}$.
}
\end{defn}
 
\begin{thm}\label{thm local isomorphism}
Let $X$ and $Y$ be strictly properly $\bigvee $-definable groups, $U$
an open properly $\bigvee $-definable neighbourhood of the identity in
$X$, and let $f:U\into Y$ be a locally strictly properly $\bigvee
$-definable homomorphism. If $X$ is definably simply connected, then
$f$ is uniquely extendible to a strictly properly $\bigvee $-definable 
homomorphism $\bar{f}:X\into Y$.  
\end{thm}

\pf
Let $x\in X$ and consider a definable path $\Gamma $ in $X$ from the 
identity $1$ to $x$.
Then, there is a finite subset $L'$ of $L$ such that 
$|\Gamma |\subseteq \cup \{x_lV:l\in L'\}$ and therefore, there are
definable paths $\Gamma _1,\dots ,\Gamma _m$ in $X$ such that 
$\Gamma =\Gamma _1\cdot \cdots \cdot \Gamma _m$ and for each 
$j\in \{1,\dots ,m\}$, $|\Gamma _j|\subseteq x_{l_j}V$ for some 
$l_j\in L'$. From this, it follows that for each $j\in \{1,\dots
,m\}$, $|\Gamma _j|^{-1}|\Gamma _j|\subseteq V\subseteq U$. Now define
$$\bar{f}_{\Gamma }(x)=f((\inf \Gamma _1)^{-1}\sup \Gamma _1)f((\inf 
\Gamma _2)^{-1}\sup \Gamma _2)\cdots f((\inf \Gamma _m)^{-1}\sup
\Gamma _m).$$
The property of $f$ in $U$ shows implies that $\bar{f}_{\Gamma }(x)$
$=\bar{f}_{\Sigma }(x)$ for any definable path in $X$ from $1$ to $x$ 
such that $\Gamma \simeq \Sigma $.

We now show that $\bar{f}_{\Gamma }(x)$ is determined independently of
the choice of the definable path $\Gamma $. Let $\Sigma $ be another 
definable path from $1$ to $x$. Since $X$ is definably simply
connected, there is a definable homotopy $\Gamma \thicksim _H\Sigma $ 
between $\Gamma $ and $\Sigma $. Similarly as before, there is a
finite subset $L''$ of $L$ such that $|H|\subseteq \cup \{x_lV:l\in
L''\}$ and for every $k$-cell $K$ ($k=0,1,2$) of $H$, 
$|K|\subseteq x_{l_K}V$ for some $l_K\in L''$. By the property of $f$
in $U$ and by the construction of $\bar{f}_{\Gamma }(x)$ its enough to
show the claim when $H$ is a $k$-cell ($k=0,1,2$) such that 
$|H|\subseteq x_lV$ for some $l\in L$ and $\Gamma \thicksim _H^0\Sigma
$. But under these assumptions, the claim is clear. 

Now define $\bar{f}:X\into Y$, by $\bar{f}(x):=\bar{f}_{\Gamma }(x)$
for some (for every) definable path $\Gamma $ in $X$ from $1$ to $x$. 
By construction, $\bar{f}$ is an extension of $f$ and, if we consider a
properly $\bigvee $-definable system of definable paths in $X$, we see
that $\bar{f}$ is a strictly properly $\bigvee $-definable
map. Moreover, $\bar{f}$ is continuous  because $f$, the
multiplication and inverse map on $X$ are continuous.

Let $x,y\in X$ and let $\Gamma $ and $\Sigma $ be definable paths in
$X$ from $1$ to $x$ and $y$ respectively. Then, $\Gamma \Sigma $ 
(notation from the proof of lemma \ref{lem pi of an h manifold}) is a 
definable path in $X$ from $1$ to $xy$ and we have $\bar{f}_{\Gamma
}(x)\bar{f}_{\Sigma }(y)=\bar{f}_{\Gamma \Sigma }(xy)$. This shows
that $\bar{f}(x)\bar{f}(y)=\bar{f}(xy)$.    
\qed

Two strictly properly $\bigvee $-definable groups $X$ and $Y$ are
called {\it locally strictly properly $\bigvee $-definably
isomorphic}, if there are locally strictly properly $\bigvee $-definable 
homomorphisms $f:U\subseteq X\into Y$ and $g:V\subseteq Y\into X$ such
that $g\circ f_{|f^{-1}(f(U)\cap V)}=1_{f^{-1}(f(U)\cap V)}$ and 
$f\circ g_{|g^{-1}(g(V)\cap U)}=1_{g^{-1}(g(V)\cap U)}$. Note that, in
a strictly properly $\bigvee $-definable covering space $(Y,p,X)$, $Y$
and $X$ are locally strictly properly $\bigvee $-definably isomorphic.

\begin{cor}\label{cor bigvee local isomorphism1}
If the strictly properly $\bigvee $-definable groups $X$ and $Y$ are 
definably simply connected, then $X$ and $Y$ are locally strictly
properly $\bigvee $-definably isomorphic iff $X$ and $Y$ are strictly 
properly $\bigvee $-definably isomorphic.
\end{cor}

\pf
Let $f:U\subseteq X\into Y$ and $g:V\subseteq Y\into X$ be as in the
definition of locally strictly properly $\bigvee $-definably
isomorphic. By theorem \ref{thm local isomorphism}, they can be
uniquely extended to strictly properly $\bigvee $-definable
homomorphisms $\bar{f}:X\into Y$ and $\bar{g}:Y\into X$. Both 
$\bar{g}\circ \bar{f}:X\into X$ and $1_X$ are extensions of the
inclusion $U\into X$. By uniqueness of the extension, we have 
$\bar{g}\circ \bar{f}=1_X$. Similarly, $\bar{f}\circ \bar{g}=1_Y$. 
Thus $X$ and $Y$ are strictly properly $\bigvee $-definably
isomorphic. The converse is clear.
\qed

\begin{cor}\label{cor bigvee local isomorphism2}
Let $X$ and $Y$ be properly $\bigvee $-definably connected, strictly 
properly $\bigvee $-definable groups and let $\tilde{X}$ and
$\tilde{Y}$ be their universal strictly properly $\bigvee $-definable 
covering spaces. Then $X$ and $Y$ are locally strictly properly 
$\bigvee $-definably isomorphic iff $\tilde{X}$ and $\tilde{Y}$ are strictly 
properly $\bigvee $-definably isomorphic.
\end{cor}

\end{subsection}

\begin{subsection}{The $m$-torsion points 
of a definable abelian group}
\label{subsection m torsion points of a definable abelian group}

In this subsection we describe $\pi _1(X)$ and the subgroup $X[m]$ of
$m$-torsion points of a strictly properly $\bigvee $-definable abelian
group $X$ for which there is a definable group
$Z$ and a strictly properly $\bigvee $-definable covering space
$(X,p,Z)$.

\begin{lem}\label{lem finite generators def case}
Let $X$ be properly $\bigvee $-definably connected, strictly properly 
$\bigvee $-definable group for which there is a definable group
$Z$ and a strictly properly $\bigvee $-definable covering space
$(X,p,Z)$. Then $\pi _1(X)$ is a finitely generated abelian group.
\end{lem}

\pf
We have $\pi _1(X)\simeq p_*(\pi _1(X))\leq \pi _1(Z)$ and by 
lemma \ref{lem pi of an h manifold}, both $\pi _1(X)$ and $\pi _1(Z)$
are abelian groups. Therefore, its enough to show that $\pi _1(Z)$ is
finitely generated. By \cite{e1} there definable groups $V\trianglelefteq
Z$ and $W\trianglelefteq U=Z/V$ such that $V$ is the maximal definable,
solvable normal subgroup of $Z$ with no definably compact parts, $W$
is  the maximal definable, definably compact, abelian normal subgroup
of $U$ and $U/W$ is a definably semisimple definable group. Moreover,
$\pi _1(Z)\simeq \pi _1(V)\times \pi _1(W)\times \pi _1(U/W)$,
$\pi _1(V)=0$ and by \cite{e2}, $\pi _1(U/W)$ is a finite group. On
the other hand, by theorem \ref{thm tietze} $\pi _1(W)\simeq G(K,T)$
for some cell decomposition $K$ of $W$ with a maximal tree $T$ 
and therefore, $\pi _1(W)$ is also
finitely generated and the result follows.
\qed

\begin{lem}\label{lem counting torsions}
Let $X$ be properly $\bigvee $-definably connected, strictly properly 
$\bigvee $-definable group and suppose that there is a definable group
$Z$ and a strictly properly $\bigvee $-definable covering space
$(X,p,Z)$. Then $X$ has unbounded exponent, the 
subgroup $Tor(X)$ of torsion points of $X$ is countable (in
particular, if $\N$ is $\aleph _0$-saturated, then $X$ has elements of 
infinite order) and, if $X$ is properly $\bigvee $-definably
complete and solvable then $X$ is abelian.
\end{lem}

\pf
This follows from similar results for definable groups 
(see \cite{e1} and \cite{s}) together with 
theorem \ref{thm main result1} and corollary \ref{cor from main}.
\qed

\begin{lem}\label{multiplication by m}
Let $X$ be a properly $\bigvee $-definably connected, strictly
properly $\bigvee $-definable abelian group.
For $m\in \NN$, let $m:X\into X$ be the multiplication by $m$
homomorphism. Then, 
$m_*:\pi _1(X)\into \pi _1(X)$ is the homomorphism defined by 
$m_*([\Gamma ])$$=m[\Gamma ]$.
\end{lem}

\pf
This is by induction on $m$. For $m=1$ the result is clear, and if it
is true for $m>1$, then $(m+1)_*([\Gamma ])=[\Gamma (m\circ \Gamma )]$
(notation of lemma \ref{lem pi of an h manifold})$=[\Gamma \cdot
(m\circ \Gamma )]$ (by lemma \ref{lem pi of an h manifold}) $=[\Gamma
][m\circ \Gamma ]$$=(m+1)[\Gamma ]$.
\qed

\begin{thm}\label{thm counting torsions def case}
Let $X$ be a properly $\bigvee $-definably connected, strictly
properly $\bigvee $-definable abelian group. Suppose that there is a 
definable group $Z$ and a strictly properly $\bigvee $-definable covering space
$(X,p,Z)$. Then $X$ is divisible, $\pi _1(X)$ is a finitely generated
torsion-free abelian group and
for each $m\in \NN$, the subgroup $X[m]$ of $m$-torsion points of $X$
is a finite group isomorphic to $\pi _1(X)/m\pi _1(X)$. In particular, $X$ is
definably simply connected iff $X$ is torsion-free. 
\end{thm}

\pf
For $m\in \NN$, let $m:X\into X$ be the multiplication by $m$
homomorphism. By lemma \ref{lem counting torsions},
$dim(m^{-1}(0))=0$. Since $mX$ is a strictly properly $\bigvee
$-definable subgroup of $X$ with $dimX=dim(mX)$, we have $mX=X$ and
$X$ is divisible. 
By theorem \ref{thm quotient by dim0}, $m$ is a strictly properly 
$\bigvee $-definable covering map and $X[m]\simeq \pi _1(X)/m_* \pi _1(X)$.
By lemma \ref{multiplication by m}, 
$X[m]\simeq \pi _1(X)/m\pi _1(X)$ and so, $X$ is definably simply 
connected iff $X$ is torsion-free. Since $\pi _1(X)$ is a strictly
properly $\bigvee $-definable subgroup of the universal strictly
properly $\bigvee $-definable covering $\tilde{X}$ of $X$, which
is torsion-free by the above, $\pi _1(X)$ is torsion-free. 
By lemma \ref{lem finite generators def case}, 
$\pi _1(X)$ is finitely generated and the result
follows.
\qed

In particular, assuming (as we conjecture),
that for a definably compact, definably connected, definable abelian
group $X$, $\pi _1(X)\simeq \ZZ$$^{dimX}$ it
follows from theorem \ref{thm counting torsions def case}, that for every 
definably compact, definably connected, definable abelian group $X$, 
$X[m]\simeq (\ZZ$$/m\ZZ$$)^{dimX}$ for every $m\in \NN$.

\end{subsection}
\end{section}

\end{document}